\documentclass[a4paper,11pt]{article}
\usepackage[square,sort,comma,numbers]{natbib}
\usepackage{amsmath}
\usepackage{accents}
\newcommand{\ubar}[1]{\underaccent{\bar}{#1}}
\usepackage{amssymb}
\usepackage{float}

\let\emptyset\varnothing
\usepackage{amsthm}
\usepackage[utf8]{inputenc}
\usepackage{bm}
\usepackage{bbm}
\usepackage{amsmath}
\usepackage{algorithm}
\usepackage{algpseudocode}
\usepackage{etoolbox}

\usepackage{varwidth}
\usepackage{framed}
\usepackage{dcolumn}
\usepackage{float}
\usepackage[T1]{fontenc}
\usepackage{graphicx}
\usepackage{longtable}
\usepackage{supertabular}
\usepackage{geometry}
\usepackage[utf8]{inputenc}
\usepackage{lipsum}
\usepackage{mathrsfs}
\usepackage{newtxmath}
\usepackage{amsmath}
\usepackage{bbm}
\usepackage{setspace}
\usepackage{subfigure}
\usepackage[raggedright ]{titlesec}
\usepackage{natbib}
\usepackage{mathtools}
\DeclarePairedDelimiter\ceil{\lceil}{\rceil}
\usepackage{graphicx}
\usepackage{multicol}
\usepackage{booktabs}
\usepackage{multirow}
\usepackage{amsmath}   
\usepackage{algorithm,algpseudocode,float}
\usepackage{lipsum}
{\theoremstyle{plain}
\usepackage{hyperref}
\hypersetup{colorlinks, linkcolor=blue}

\usepackage[table]{xcolor}

\setcitestyle{square}
\allowdisplaybreaks 
\usepackage{url}
\newcommand{\propnumber}{} 

\usepackage{amsthm}
\theoremstyle{definition}
\newtheorem{theorem}{Theorem}
\newtheorem{prop}[theorem]{Proposition}
\newtheorem{remark}{Remark}
\newtheorem{example}{Example}
\newtheorem{assumption}{Assumption}
\newtheorem{cor}[theorem]{Corollary}
\newtheorem{lem}[theorem]{Lemma}

\usepackage{nomencl}
\usepackage{pgfplots}
\usepackage{etoolbox}
\renewcommand\nomgroup[1]{%
  \item[\bfseries
  \ifstrequal{#1}{A}{Indices and Sets}{%
  \ifstrequal{#1}{B}{Parameters}{%
  \ifstrequal{#1}{C}{Decision variables}{%
  \ifstrequal{#1}{D}{Random Parameter}{}}}}%
]}
\usepackage{amsfonts}
\usepackage{xpatch}
\usepackage{mathtools}
\xpatchcmd{\thenomenclature}{%
  \section*{\nonauthor}
}{
}{\typeout{Success}}{\typeout{Failure}}
\makeatletter

\makeatother
\usepackage{tikz}
\usepackage{pgfplots}
\usepgfplotslibrary{units}
\usepackage{filecontents}
\usetikzlibrary{arrows.meta,calc,fit,shapes.geometric,positioning,backgrounds}
\usepackage{authblk}
\usepackage{blindtext}
\usepackage{caption}
\usepackage{pgfplotstable}
\usepackage{lscape}
\DeclareMathOperator*{\argmax}{argmax}

\usetikzlibrary{plotmarks}
\usetikzlibrary{patterns}
\usepackage{pgfplotstable}
\usepgfplotslibrary{groupplots}
\usepackage{boldline,multirow}

\makenomenclature 
\RequirePackage{ifthen}
\renewcommand{\nomname}

\def\bf{\mathbf}
\def\it{\mathit}

\def\cal{\mathcal}
\def\bsm{\boldsymbol}
\newcounter{alg}

\begin{document}
\title{Two-Stage Distributionally Robust Optimization: Intuitive Understanding and Algorithm Development from the Primal Perspective}
\author {Zhengsong Lu and Bo Zeng}
\affil{Department of Industrial Engineering, University of Pittsburgh}
\date{}
\maketitle

\renewcommand{\baselinestretch}{1.5}
\noindent
\abstract{
	In this paper, we study the two-stage distributionally robust optimization (DRO) problem from the primal perspective. Unlike existing approaches, this perspective allows us to build a deeper and more  intuitive  understanding on DRO, to leverage classical and well-established solution methods and to develop a general and fast decomposition algorithm (and its variants), and to address a couple of unsolved issues that are critical for modeling and computation. 
	Theoretical analyses regarding the strength, convergence, and iteration complexity of the developed algorithm are also presented. 
	A numerical study on different types of instances of the distributionally robust facility location problem demonstrates that the proposed solution algorithm (and its variants) significantly outperforms existing methods. It solves instances up to several orders of magnitude faster, and successfully addresses new types of practical instances that previously could not be handled. We believe these results will significantly enhance the accessibility of DRO, break down barriers, and unleash its potential to solve real world challenges.
}

\section{Introduction}
Distributionally robust optimization (DRO) has emerged as a powerful optimization paradigm to address uncertainties in decision making. As a flexible unification of stochastic programming (SP) and robust optimization (RO), it adopts an \textit{ambiguity set}, which leverages available 
distributional information  while accounting for unforeseen perturbations on top of that distributional information, to describe randomness. Hence, rather than assuming a single fixed probability distribution, which is the central idea behind SP, DRO takes an RO perspective to hedge against all (often infinitely many) possible distributions within that ambiguity set to ensure robustness. Obviously, the derived solution is more robust than that of SP.  
Also, by restricting the perturbations in parameters of the underlying distribution, the solution is less conservative compared to that of RO, which just focuses on worst-case scenarios without considering any distributional information. Actually, if those parameters are exact and complete that defines a single distribution, DRO reduces to SP. And if the perturbations are allowed to be arbitrarily large, DRO reduces to RO.

Because of its flexibility in handling data and strength in taking advantage of the associated information, DRO is believed to be an ideal data-driven decision making tool. On one hand, it attracts a large amount of methodological studies in the literature, including those on the design of ambiguity sets and their mathematical and statistical properties, new variants and customizations for different modeling purposes, and strong solution methods to handle complex data/structures, including exact and approximate mixed integer reformulations and specific algorithms.  On the other hand,  it has found applications in various fields, including finance, energy, supply chain management, and machine learning, where some data on random parameters are available but the true distribution may not be known precisely. 

Regardless the fact that the first recognized study on DRO \cite{scarf1958min} appeared in 1958, a long time ago, DRO is notably less popular compared to SP or RO, particularly in the context of two-stage decision making. We observe that a couple of issues contribute to this situation. The first one is that present solution methods, which are reviewed in the next section, are, not sufficiently strong or lack scalability to handle practical-scale problems. Another one is that some fundamental issues have not been sufficiently addressed, e.g., the feasibility requirement of the recourse problem and the non-convex ambiguity sets, and hence restrictive assumptions are often imposed to circumvent those situations.  Actually, existing studies tend to be theoretically sophisticated, providing less intuitive understanding to appreciate DRO. This high theoretical threshold renders DRO less accessible to practitioners, further hindering its applications in the real world.  

In this paper, instead of following the mainstream that adopts the dual perspective to study DRO and to develop solution methods, we take the primal perspective to analyze its structure. This ``new'' angle offers prominent advantages that help us address aforementioned barriers. The reminder of this paper will demonstrate those advantages, with an emphasis on intuitive understanding and developing general and fast algorithms that drastically outperform the state-of-the-art. We first review the relevant literature in Section \ref{sect_literature}. Then, Section \ref{sect_primal_DRO} analyzes and computes the worst-case expected value from the  primal perspective that is more accessible and general.     
 In Section \ref{sect_CCG_DRO}, we present the whole computational algorithm in its basic form for two-stage DRO, with a couple of variants described in Section \ref{sect_CCGvariants} to achieve a stronger solution capacity. Section \ref{sect_computation} reports and analyzes the results of numerical experiments, while Section \ref{sect_conclusion} provides the conclusions of this paper.     

\textbf{Notation}:
We generally denote matrices in upper-case, vectors in bold lower-case and scalars in lower-case, unless explicitly noted otherwise. Special notations include that: an upper-case letter in calligraphic denotes a major set, $\hat{}$ \ generally denotes a subset and \ $\widehat{}$ \ denotes a union of subsets,  $\bf M$ represents a sufficiently large positive number, and $[n]$ indicates the set of positive integer numbers from $1$ to $n$.

\section{Related Research in the Literature}
\label{sect_literature}
In this section, we review the literature that is most relevant to the work presented in this paper. Note that over the last 10 years, many studies on two-stage DRO (and those with more stages) have been published, which involve making recourse decisions after the randomness is cleared  \cite{rahimian2022frameworks}.  Consider the following tri-level two-stage DRO formulation
\begin{align}
\label{eq_2stgDRO}
\mathbf{2-Stg \ DRO}: \ \ \ \ w^* =& \min_{\bf x\in \mathcal{X}} f_1(\bf x)+\sup_{P\in \mathcal{P}} E_P[Q(\bf x, \xi)],
\end{align}
where $\mathcal X$ denotes the feasible set of the first-stage (also known as ``here-and-now'') decision variables $\bf x$, and $\xi$ the random vector representing the uncertain parameter, and $Q(\bf x, \xi)$ the value function capturing the optimal value of the second-stage recourse problem  
\begin{align}
\label{eq_recourse}
\mathbf{Reco\ Prob}: & \ \min\Big\{f_2(\bf x,\bsm \xi,\bf y): \bf y\in\mathcal{Y}(\bf x, \xi)\Big\}.   
\end{align} 
Variables $\bf y$ denote the recourse (also known as ``wait-and-see'') decision variables, which may contain discrete ones, and $\mathcal Y(\bf x,\xi)$ is their feasible set. We highlight that $Q(\bf x, \xi)$ returns the optimal value of the recourse problem for given $\bf x$ and $\xi$, not representing the whole recourse problem.  
When the recourse problem is infeasible, i.e., $\mathcal{Y}(\bf x, \xi)=\emptyset$, $Q(\bf x,\xi)$ is set to $+\infty$ by convention.  Random vector $\xi$ is defined on a measurable space $(\Xi, \mathcal F)$ with $\Xi\subseteq\mathbb{R}^{n_\xi}$ being the support and $\mathcal{F}$ a $\sigma$-algebra that contains all singletons in $\Xi$. It follows probability distribution $P\in \cal{P}$ where $\cal P\subseteq \mathcal{M}(\Xi, \mathcal F)$ is a family of probability distributions and  $\mathcal{M}(\Xi, \mathcal F)$ is the collection of all probability distributions on $(\Xi,\mathcal{F})$. Note that $\cal P$, assumed to be non-empty in this paper, is commonly referred to as the ambiguity set in the DRO literature.

It can be seen that the fundamental feature that distinguishes DRO from SP or RO is its ambiguity set. Based on the current literature, ambiguity sets generally can be grouped into two main categories: moment-based ones and statistical discrepancy-based ones. For the moment-based ambiguity, the set is usually defined by the moments of $\xi$ across the entire support $\Xi$, e.g.,  \cite{delage2010distributionally,wiesemann2014distributionally}. 
For discrepancy-based ambiguity sets,  examples include the $\phi$-divergence ambiguity sets \cite{ben2013robust,lam2019recovering} and Wasserstein ambiguity sets \cite{gao2023distributionally,mohajerin2018data,zhao2018data}. The latter ones are becoming more popular as they are directly defined with respect to empirical data, a natural demonstration of the data-driven concept.  A few less-investigated ambiguity sets also do exist in the literature (\cite{rahimian2022frameworks}). We note that, regardless of the structure of $\Xi$, $\cal{P}$ is always represented as a convex set in the existing literature. For sophisticated  $\cal{P}$ that requires a mixed integer representation, we have not been aware of any study available yet.

Theoretically, two-stage DRO has been proven to be  NP-hard in general \cite{ben2004adjustable,bertsimas2010models}. To solve this challenging problem exactly,  there exist two types of popular solution strategies. One is to reformulate it into single-level (nonlinear)  mixed integer programs (MIPs) that can be directly handled by existing solvers or computational algorithms, e.g.,\cite{hanasusanto2018conic,mohajerin2018data,zhao2018data,xie2020tractable}. Nevertheless, those reformulations are generally with very large (often infinite) number of variables and/or constraints. Directly computing them turns out to be practically infeasible. Hence, another strategy is to develop decomposition algorithms so that variables and/or constraints are introduced on the fly to reduce the computational burden, e.g., \cite{bagheri2018resilient,bansal2018decomposition,saif2021data,gamboa2021decomposition,gangammanavar2022stochastic,duque2022distributionally}. It is worth highlighting that the fundamental idea underlying those two strategies is to leverage the  duality results of the integration-based convex problems (i.e., the $\sup$ part  in \eqref{eq_2stgDRO}) to reduce it to a monolithic formulation \cite{shapiro2001duality,mohajerin2018data,zhao2018data,blanchet2019quantifying,gao2023distributionally}. 

Regarding decomposition algorithms, they are generally designed with a master-subproblem framework that runs iteratively before convergence. Similar to the case of RO~\cite{ZENG2013457}, those algorithms can be classified into two groups based on the representation of information feedback from the subproblem to cut the current solution of the master problem.   One group is Benders-dual (BD) type  of algorithms  \cite{bansal2018decomposition,duque2022distributionally,gamboa2021decomposition,gangammanavar2022stochastic,luo2022decomposition}, which use the dual information of the recourse problem and define cutting planes in the form of Benders (optimality) cut. Another group uses the column-and-constraint (C\&CG) method to generate cut represented by a new replicate of the whole recourse problem \cite{bagheri2018resilient,saif2021data,gamboa2021decomposition,el2024distributionally}. Common to those decomposition algorithms is that their master problems are augmented over iterations to derive stronger relaxations and therefore tighter lower bounds for \eqref{eq_2stgDRO}. On the other hand, the best feasible solution found so far provides an upper bound. Once those two bounds meet or are with a sufficiently small difference, the optimal value and an optimal first stage solution to $\mathbf{2-Stg \ DRO}$ in \eqref{eq_2stgDRO} are obtained.

While those decomposition algorithms are iterative, they, by making use of professional MIP solvers to compute master and subproblems, often demonstrate significantly better computational performance than directly computing DRO's monolithic reformulations. Also, the emerging features of those solvers, such as lazy constraints and the ability to handle bilinear constraints, can be leveraged to increase our solution capacity to address more complex practical problems. In addition to those general and exact methods, several algorithms that handle special structures (e.g., \cite{long2024supermodularity}) or compute approximate solutions (e.g., \cite{jiang2018risk,bertsimas2019adaptive,georghiou2021optimality}) have also been investigated in the literature, which are beyond the scope of this paper. Regardless the aforementioned many efforts on developing  algorithms for general $\mathbf{2-Stg \ DRO}$, we note that the following fundamental issues remain unsolved or seriously underinvestigated. 

\noindent ($i$) No existing algorithm has addressed the feasibility  (equivalently the infeasibility) issue of the recourse problem yet, noting that in general some $\bf x$ may render \eqref{eq_recourse} infeasible under some scenario $\xi\in \Xi$. This challenge has always been circumvented by assuming that $Q(\bf x, \xi)$ is finite for any $\bf x\in \mathcal{X}$ and $\xi \in \Xi$ combination. Nevertheless, such a strategy is not realistic or practically acceptable for many real systems (e.g., critical infrastructure systems), as decision makers might be required to  identify an $\bf x$ that guarantees the feasibility of \eqref{eq_recourse} with respect to $\mathcal P$. Moreover, we note from our numerical studies that addressing the feasibility issue in the context of DRO could be computationally much more demanding.  \\    
$(ii)$ No study has been performed on complex ambiguity sets that are non-convex. Such type of ambiguity sets may happen if data come from different sources and they are rather not consistent, which nevertheless can be represented by employing mixed integer sets.  
Actually, due to their non-convexity, we believe that simply taking the dual perspective, which is behind current reformulations or solution algorithms,  to handle this type of ambiguity sets is not a valid or viable direction. Hence, it would be necessary to extend our solution capacity to solve $\mathbf{2-Stg \ DRO}$ with mixed integer ambiguity sets.

Also, as noted, existing studies on DRO heavily rely on sophisticated mathematical concepts and derivations with different assumptions or conditions. It is desired to develop simpler, more intuitive, and general reasoning and analyses, which will undoubtedly inspire more scholars and practitioners to study and apply DRO. Indeed, a recent study~\cite{zhang2022simple} has made an effort to simplify complex mathematical derivations in the context of Wasserstein ambiguity set based DRO, aiming to provide a better accessibility and applicability.  The last and probably the most important issue is our rather limited computational capacity to handle general, large-scale and practical DRO instances, which is significantly inferior to what we currently have for instances of SP and RO.

In this paper, we make an effort to address those critical challenges from the primal perspective. Rather than relying on the  duality of convex programs constructed on $\cal{P}$ to develop analytical insights and computational methods, we adopt a rather simple approach to directly compute the worst-case expected value (WCEV) and leverage it to solve $\mathbf{2-Stg \ DRO}$. It is worth noting that this strategy offers us three key advantages. First, it helps us establish a clear and intuitive understanding of DRO that facilitates its applications in practice. Second, it allows us to directly investigate and attack more challenging structures that might not be approachable by existing methods. Third, it provides new opportunities to develop stronger and more general solution algorithms. The remaining of this paper will showcase those advantages in our analyses and algorithm development.

\section{Computing the WCEV from the Primal Perspective}
\label{sect_primal_DRO}
	
In this section, we consider a core structure and challenge embedded in solving $\mathbf{2-Stg \ DRO}$, i.e., how to compute the worst-case expected value with respect to $\cal P$. Since we do not assume any particular value of $\bf x$, we ignore $\bf x$ in all derivations and analyses of this section to simplify our exposition.

The basic intuitions behind our study and algorithm design  are: $(i)$ Generally any probability distribution can be represented or approximated with an arbitrary accuracy by a discrete probability distribution; $(ii)$ The WCEV over an ambiguity set can be obtained by solving a typical finite mathematical program or by customizing and implementing some well-known optimization algorithms.   
	
\subsection{Expected Value and Worst-Case Expected Value}
\label{sect_EVWCEV}
In this subsection, we present some results regarding the  computations of the expected value for any distribution and the WCEV over a set of probability distributions.  Different from the current mainstream understandings obtained from a dual perspective, those results are derived by directly taking the primal perspective to appreciate those two expected values. They also provide theoretical support for us to develop corresponding algorithms. Recall that	$\mathcal M(\Xi, \mathcal F)$ denotes the collection of  probability distributions on $(\Xi, \mathcal F)$, and consider a real-valued function $Q(\xi)$ defined on $\Xi$.     
For the remainder of this paper, we make a couple of very mild assumptions. 
 \begin{assumption}
 	($i$) Support $\Xi$ is a closed and bounded set; $(ii)$ $Q(\xi)>-\infty$ for $\xi \in \Xi$.
 \end{assumption}
 We mention that after $\bf x$ is introduced in Section \ref{sect_CCG_DRO}, the second part of this assumption is revised to $Q(\bf x, \xi)>-\infty$ for $\bf x\in \mathcal X$ and $\xi \in \Xi$. Next, we consider the general integration based expected value for an arbitrary continuous distribution. 
\begin{theorem}
\label{thm_EV_sum}
 Assume that   $Q(\Xi)$  is Lebesgue integrable over $\underline\Xi\equiv\{\xi\in \Xi: Q(\xi) <+\infty\}$. Let $\Big(\bsm\xi(n),P(\bsm\xi(n))\Big)$ denote a discrete probability distribution over $\bsm\xi(n)\equiv\big(\xi_1(n),\dots,\xi_n(n)\big)$ with $P(\bsm\xi(n))\equiv\Big(p(\xi_1(n)),\dots,p(\xi_n(n))\Big)$ being the associated probability. For any probability distribution  $P \in \mathcal M(\Xi, \mathcal F)$,  there exists a sequence of discrete probability distributions in the form of
 $\Bigg\{\Big(\bsm\xi(n), P(\bsm\xi(n))\Big)\Bigg\}_{n=1}^{+\infty}$
 such that 
 \begin{equation}
E_{P}[Q(\xi)]=\int_{\xi\in \Xi} Q(\xi)P(d\xi) = \lim\limits_{n \rightarrow + \infty} \sum\limits_{j=1}^n Q\big(\xi_j(n)\big)p(\xi_j(n)). \tag*{\qed}
\end{equation}
\end{theorem}

We note that Theorem \ref{thm_EV_sum}  is rather mathematically intuitive. It indicates that we may be able to replace the integration operation by a weighted sum  over a set of discrete scenarios, which should be computationally much more friendly. Actually, if we aim to compute the WCEV over a family of distributions, i.e., the ambiguity set underlying DRO, it has been reported that the worst-case one may reduce to a discrete probability distribution. Next, we consider the following ambiguity set that generalizes all known moment-based  and Wasserstein ambiguity sets adopted in the DRO literature.    
\begin{align}
	\label{eq_ambiguity}
	\mathcal{P} = \Big\{P \in \mathcal{M}(\Xi,\mathcal{F}): \int_{\Xi}  P(d\mathbf \xi) = 1, \, \int_{\Xi} \mathbf \psi_i(\xi)P(d\xi) \leq \gamma_i\ \ \forall i\!\in [m]\Big\}.
\end{align}
Function $\psi_i(\cdot)$ is real valued and bounded on $\Xi$ for $i=1,\dots, m$. With this ambiguity set, we define the following optimization problem to compute the WCEV. 
\begin{align}
	\label{eq_WCEV}	\mathbf{WCEV}: \sup_{ P\in \mathcal P} E_P[Q(\bf \xi)]\!\equiv\!\sup\Big\{\int_{\Xi} Q(\xi)P(d\xi): P\in \cal P\Big\}.
\end{align}
When $\Xi$ is a finite discrete set, we have $P\equiv(p_1,\dots,p_n)$ with $n=|\Xi|$. And the integration in \eqref{eq_ambiguity} and \eqref{eq_WCEV} reduces to summation. As the analysis for this case is simpler, we, unless noted otherwise, mainly present our derivations and analyses for the case where $|\Xi|$ is not finite. We mention that those results, with little or  minor changes, are naturally applicable when $|\Xi|$ is finite.    
\begin{remark}
	It is clear that an ambiguity set defined by moment inequalities is a special case of \eqref{eq_ambiguity}. Also, an ambiguity set defined by Wasserstein metric, as explained in Section \ref{subsect_Wasserstein}, belongs to  \eqref{eq_ambiguity}. Regarding  a $\phi$-divergence-based ambiguity set, we can approximate $\phi$, a non-negative and convex function,  by a piece-wise linear one with an arbitrary accuracy. As a result, its mathematical representation is also in the form of \eqref{eq_ambiguity}. 
	\hfill \qed
\end{remark}

  Consider an arbitrary $P\in \cal P$ and let $\varepsilon$ be a sufficiently small positive constant. Following the proof  of Theorem  \ref{thm_EV_sum} in Appendix A.1, it is easy to see that we can always identify a partition of $\Xi$ for $\psi_i$ and derive a discrete distribution that approximates $P$ with a violation of $\gamma_i$ up to $\varepsilon$. By taking the intersection of those partitions for $Q$ and $\psi_i$ for $i\in [m]$, a finer partition yielding a stronger approximation to $P$ with less violation of all $\gamma_i$'s will be obtained. Such a construction procedure leads to the next result readily that allows us to leverage a sequence of discrete distributions to compute $\mathbf{WCEV}$ formulation with an arbitrary accuracy.  
  Before that, we define that a  discrete probability distribution  $\Big(\bsm\xi(n), P(\bsm\xi(n))\Big)$ is $\varepsilon$-feasible to $\cal P$ if $\displaystyle\sum_{j=1}^n p(\xi_j(n))\psi_i(\xi_j(n))\leq \gamma_i+\varepsilon$ for $i\in [m].$ Moreover, it is an  $\varepsilon$-approximation to $P\in \cal P$ if $\big|\int_{\Xi} Q(\xi)P(d\xi)-\displaystyle\sum_{j=1}^n p(\xi_j(n))Q(\xi_j(n))\big|\leq \varepsilon.$

\begin{cor}
	\label{cor_eps_optimal_feasible}
Suppose that $\displaystyle{\sup_{ P\in \mathcal P} E_P[Q(\bf \xi)]}$ is finite.    $(i)$ For every $P\in \cal P$, there exists an $\varepsilon$-feasible discrete distribution, denoted by  $\Big(\bsm\xi(n), P(\bsm\xi(n))\Big)$, that is an $\varepsilon$-approximation to $P$. $(ii)$ There exists an $\varepsilon$-feasible discrete distribution, denoted by  $\Big(\bsm\xi'(n), P'(\bsm\xi'(n))\Big)$, that is $\varepsilon$-optimal to the $\mathbf{WCEV}$ formulation, i.e.,   
\begin{equation}
\big|\sup_{ P\in \mathcal P} E_P[Q(\bf \xi)]-\displaystyle\sum_{j=1}^n p'(\xi'_j(n))Q(\xi'_j(n))\big|\leq \varepsilon. \tag*{\qed}
\end{equation}
\end{cor}

Note that the result in Corollary \ref{cor_eps_optimal_feasible} holds in general as it does not require any restrictive condition. Yet, if an actual solution  to  the $\mathbf{WCEV}$ problem is desired, Corollary \ref{cor_eps_optimal_feasible} does not help to ensure the existence of an optimal $P$ or to obtain a discrete distribution that is feasible and ($\varepsilon$-)optimal. Next, we present some sufficient conditions that substantially support us on those issues.

\begin{theorem}
	\label{thm_feasible_EV_sum}
	Suppose $Q(\cdot)$ is an upper semicontinuous function and $\psi_i(\cdot)$ are lower semicontinuous functions over $\Xi$ for $i\in[m]$. The optimal value of $\mathbf{WCEV}$ can be attained, i.e., $\displaystyle{\sup_{ P\in \mathcal P} E_P[Q(\bf \xi)] = \max_{ P\in \mathcal P} E_P[Q(\bf \xi)]}$. Moreover, if $\cal P$ has an interior point $P^0$ and function $g(t)$,  defined as:
	\[g(t) = \sup_{P \in \mathcal{M}(\Xi, \mathcal F)}\Big\{\int_{\Xi} Q(\xi)P(d\xi): \int_{\Xi}  P(d\mathbf \xi) = 1, \, \int_{\Xi} \psi_i(\xi)P(d\xi) \leq \gamma_i - t\ \ \forall i\!\in [m]\Big\},\]
	is continuous at $t=0$, there exists a sequence of discrete probability distributions such that their associated expected values of $Q(\cdot)$'s converge to $\displaystyle{\max_{ P\in \mathcal P} E_P[Q(\bf \xi)]}$. \hfill \qed
\end{theorem}

\begin{proof}
		We start with the proof for the first statement. Since $\Xi$ is a compact metric space and $Q(\cdot)$ is upper semicontinuous, $Q(\cdot)$ is bounded and the set of all probability measures $\mathcal P(\Xi)$ on $(\Xi, \mathcal F)$ is weakly compact. Let $\{P^k\}_k \subseteq \mathcal P(\Xi)$ be a sequence converging weakly to $P^{\infty} \in \mathcal P(\Xi)$. According to Portmanteau Theorem, given that $\psi_i(\cdot)$ for $i\!\in [m]$, is lower semicontinuous and bounded from below on $\Xi$, we have
	\[\liminf_{k \to \infty} \int_{\Xi} \psi_i(\xi)P^k(d\xi) \geq \int_{\Xi} \psi_i(\xi)P^{\infty}(d\xi),\]
	indicating that $P \mapsto \int_{\Xi} \mathbf \psi_i(\xi)P(d\xi)$ is weakly lower semicontinuous on $\mathcal P(\Xi)$.  Given that the pre-image of $(-\infty, \gamma_i]$ for $i\!\in [m]$, under a lower semicontinuous mapping is closed, $\mathcal P_{\gamma_i} = \Big\{P \in \mathcal P(\Xi):\int_{\Xi} \mathbf \psi_i(\xi)P(d\xi) \leq \gamma_i \Big\}$ for $i\!\in [m]$, is weakly closed. As $\mathcal P = \displaystyle\bigcap_{i \in [m]} \mathcal{P}_{\gamma_i} \subseteq \mathcal{P}(\Xi)$, $\mathcal P$ is weakly closed and tight, i.e., $\mathcal P$ is weakly compact. Applying Portmanteau Theorem again, we can conclude that $P \mapsto \int_{\Xi} Q(\xi)P(d\xi)$ is weakly upper semicontinuous on $\mathcal P$. Therefore, by Weierstrass' Theorem, the optimal value of $\mathbf{WCEV}$ problem can be attained.
	
   We next prove the second statement. Our proof depends on the construction of an auxiliary sequence of measures $\{P^k\}_k$:
	\begin{align}
	\label{eq_proof_P*_exist}
	P^k = \argmax_{P \in \mathcal{M}(\Xi, \mathcal F)}\Big\{\int_{\Xi} Q(\xi)P(d\xi): \int_{\Xi}  P(d\mathbf \xi) = 1, \, \int_{\Xi} \psi_i(\xi)P(d\xi) \leq \gamma_i - \frac{\Delta}{k}\ \ \forall i\!\in [m]\Big\}
	\end{align}
	where $\Delta$ is a fixed value that ensures the feasible set of  \eqref{eq_proof_P*_exist} is consistent for all $k$, i.e., \eqref{eq_proof_P*_exist} is feasible,  and thus $P^k$'s are attained  for any $k$. Note that $\Delta$ can always be found given the existence of interior point $P^0$. Clearly, the sequence,  $\Big\{\int_{\Xi} Q(\xi)P^k(d\xi)\Big\}_k$, is monotonically increasing and bounded above by the optimal value of $\mathbf{WCEV}$ problem. Because of the continuity of $g(t)$ at $t=0$ and the  monotone convergence theorem, we have
	\[\lim_{k \rightarrow \infty} \int_{\Xi} Q(\xi)P^{k}(d\xi) = \max_{ P\in \mathcal P} E_P[Q(\bf \xi)].\]

	Using the argument presented before Corollary \ref{cor_eps_optimal_feasible}, for $P^k$, we can obtain a  discrete distribution based on a partition $\{\Xi_j(n_k)\}_j$ of $\Xi$ that satisfies
    $$\mathrm{both} \ \ \ \lim\limits_{n_k \rightarrow + \infty} \sum\limits_{j=1}^{n_k} Q\big(\xi_j(n_k)\big)p\big(\xi_j(n_k)\big) = \int_{\Xi} Q(\xi)P^k(d\xi)$$
     $$\mathrm{and } \  \ \ \ \lim\limits_{n_k \rightarrow + \infty}\sum\limits_{j=1}^{n_k} \psi_i\big(\xi_j(n_k)\big)p\big(\xi_j(n_k)\big) = \int_{\Xi} \psi_i(\xi)P^k(d\xi).$$

     For the discrete distribution associated with $P^k$, we let \begin{align*}
	 N_k = \min\Bigg\{n_k: & \Big|\sum\limits_{j=1}^{n_k} Q\big(\xi_j(n_k)\big)p\big(\xi_j(n_k)\big) - \int_{\Xi} Q(\xi)P^k(d\xi) \Big|\leq \frac{\Delta}{k}, \\  & \Big|\sum\limits_{j=1}^{n_k} \psi_i\big(\xi_j(n_k)\big)p\big(\xi_j(n_k)\big) - \int_{\Xi} \psi_i(\xi)P^k(d\xi) \Big| \leq \frac{\Delta}{k}\ \ \forall i\in [m] \Bigg\}.
   \end{align*}
	\noindent Claim:  All discrete distributions  in  $\Bigg\{\Big(\bsm\xi(N_k), P(\bsm\xi(N_k))\Big)\Bigg\}_{k=1}^{+\infty}$ are feasible and $$\lim_{k \rightarrow \infty}\sum_{j=1}^{N_k} Q\big(\xi_j(N_k)\big)p\big(\xi_j(N_k)\big)=\max_{ P\in \mathcal P} E_P[Q(\bf \xi)].$$\\
	\textit{Proof of Claim:} Obviously, any discrete distribution in $\Bigg\{\Big(\bsm\xi(N_k), P(\bsm\xi(N_k))\Big)\Bigg\}_{k=1}^{+\infty}$ is feasible as it satisfies all constraints in $\mathcal P$. For every $\varepsilon>0$, there is $k(\varepsilon)$ such that for all $k\geq k(\varepsilon)$, we have
	\begin{align*}
		\Big|\max_{ P\in \mathcal P} E_P[Q(\bf \xi)] -\sum\limits_{j=1}^{N_k} Q\big(\xi_j(N_k)\big)p\big(\xi_j(N_k)\big) \Big| \ \leq \ & \Big|\max_{ P\in \mathcal P} E_P[Q(\bf \xi)] - \int_{\Xi} Q(\xi)P^k(d\xi)\Big| \ \\
		+\ \Big|\int_{\Xi} & Q(\xi)P^k(d\xi) -\sum\limits_{j=1}^{N_k} Q\big(\xi_j(N_k)\big)p\big(\xi_j(N_k)\big) \Big| \leq \varepsilon.
	\end{align*}
	The first inequality follows from the triangle inequality.  As for the second inequality, the convergence of $\Big\{\int_{\Xi} Q(\xi)P^k(d\xi)\Big\}_k$ guarantees the difference between two integrals is not more than $\frac\varepsilon{2}$. Furthermore, by the definition of $N_k$ and $\Bigg\{\Big(\bsm\xi(N_k), P(\bsm\xi(N_k))\Big)\Bigg\}_{k=1}^{+\infty}$, the difference between the sum and the integral can be made no more than  $\frac\varepsilon{2}$ when $k(\varepsilon) \geq \frac{2\Delta}{\varepsilon}$, which results in the second inequality. Hence, the claim holds when $k\rightarrow + \infty$.

	With this claim proved, the second statement of this theorem follows directly.
\end{proof}

By simply making use of Theorems \ref{thm_EV_sum} and \ref{thm_feasible_EV_sum}, we have the next  result. 
\begin{cor}
\label{cor_WCEV_limit}
Assuming that all the sufficient conditions in Theorem \ref{thm_feasible_EV_sum} hold, the next formulation is equivalent to $\mathbf{WCEV}$ problem in \eqref{eq_WCEV}.
\begin{equation}
\label{eq_limit_WCEV}
\max_{ P\in \mathcal P} E_P[Q(\bf \xi)]= \lim_{n\rightarrow\infty}\max_{\Big(\bsm\xi(n), P(\bsm\xi(n))\Big)}\Big\{\sum_{j=1}^nQ\big(\xi_j(n)\big)p(\xi_j(n)): P\in \cal P, \xi_j(n)\in \Xi\ \ \forall j\in [n]\Big\}
\end{equation}
\end{cor}
\begin{proof} 
	By the definition of $\mathbf{WCEV}$ problem, it is clear that 
\[\max_{ P\in \mathcal P} E_P[Q(\bf \xi)] \geq \lim_{n\rightarrow\infty}\max_{\{(\xi_j(n), p(\xi_j(n)))\}_{j=1}^n}\Big\{\sum_{j=1}^nQ(\xi_j(n))p(\xi_j(n)): P\in \cal P, \xi_j(n)\in \Xi\ \ \forall j \in [n]\Big\}.\]
According to Theorem \ref{thm_feasible_EV_sum}, there exists a sequence of discrete distributions in $\mathcal P$ such that
$\displaystyle\max_{ P\in \mathcal P} E_P[Q(\bf \xi)] =\lim_{n\rightarrow\infty} \sum_{j=1}^nQ(\tilde \xi_j(n))\tilde p(\tilde \xi_j(n)).$ 
Given that $$\displaystyle\sum_{j=1}^nQ(\tilde \xi_j(n))\tilde p(\tilde \xi_j(n)) \leq \displaystyle\max_{\{(\xi_j(n), p(\xi_j(n)))\}_{j=1}^n}\sum_{j=1}^nQ(\xi_j(n))p(\xi_j(n)),$$ we have
\begin{align*}
		\max_{ P\in \mathcal P} E_P[Q(\bf \xi)]\leq \lim_{n\rightarrow\infty}\max_{\{(\xi_j(n), p(\xi_j(n)))\}_{j=1}^n}\Big\{\sum_{j=1}^nQ(\xi_j(n))p(\xi_j(n)): P\in \cal P, \xi_j(n)\in \Xi \ \ \forall j\in [n]\Big\}.
\end{align*}
Hence, the desired result simply follows.
\end{proof}

Actually, it is worth highlighting that $\mathbf{WCEV}$ problem in \eqref{eq_WCEV} can be seen as an infinite-column linear program. Under some little bit stronger conditions, a more insightful result has been derived \cite{shapiro2001duality}: $\mathbf{WCEV}$ problem, if subject to equality constraints in \eqref{eq_ambiguity} and with a finite optimal value, has an optimal solution that has at most $(m+1)$ scenarios with non-zero probabilities. This insight actually is verified in our  numerical result presented in Section~\ref{sect_computation}. 
 Extending and making use of this result,  we can build a finite mathematical program (FMP) that is concise and helps to solve $\mathbf{WCEV}$ problem exactly.

\begin{prop}
\label{prop_WCEV_OPT}
Suppose $Q(\cdot)$ is upper semicontinuous and $\psi_i(\cdot)$ are continuous over $\Xi$ for $i\in [m]$. The WCEV, i.e., the optimal value of $\mathbf{WCEV}$ problem, is attainable and can be obtained by solving the following FMP
\begin{align}
\label{eq_WCEV_opt}
\begin{split}
\mathbf{WCEV-FMP}: \max\Big\{\sum_{j=1}^{m+1} Q(\xi_j)p_j: & \sum_{j=1}^{m+1}p_j=1, \  \sum_{j=1}^{m+1}\psi_i(\xi_j)p_j\leq \gamma_i\ \ \forall i\in[m],\\ 
&\xi_j\in \Xi\ \ \forall j\in [m+1],\ p_j\geq 0\ \ \forall j\in [m+1]\Big\}. 
\end{split}
\end{align}
Its optimal solution, denoted by $(P^*,\xi^*)$ with $P^*\equiv(p^*_1,\dots,p^*_{m+1})$ and $\xi^*\equiv(\xi^*_1,\dots,\xi^*_{m+1})$, is then feasible and optimal to  $\mathbf{WCEV}$ problem. \hfill \qed
\end{prop} 

\begin{proof}
		Given that $\psi_i(\cdot)$ for $i\in [m]$ are continuous, according to Theorem \ref{thm_feasible_EV_sum}, the optimal value of the original $\mathbf{WCEV}$ in \eqref{eq_WCEV} can be attained. Let $P'(\Xi)$ be an optimal solution and compute $\gamma_i' = \int_{\Xi} \psi_i(\xi)P'(d\xi)$ for $i\in [m]$. Then, the original $\mathbf{WCEV}$ problem is equivalent to the following one with $(m+1)$ equality constraints.
	\begin{align}
	\label{eq_proof_geq_eq}
		\max_{P \in \mathcal M(\Xi, \mathcal F)}\Big\{\int_{\Xi} Q(\xi)P(d\xi): & \int_{\Xi} P(d\xi)=1, \  \int_{\Xi} \psi_i(\xi)P(d\xi) = \gamma_i'\ \ \forall i\in[m]\Big\}.
	\end{align}
	According to Richter-Rogosinski Theorem \cite{shapiro2001duality}, \eqref{eq_proof_geq_eq} has an optimal distribution with a support of at most  $(m+1)$ scenarios. Clearly, it is also feasible and optimal to the original $\mathbf{WCEV}$ problem. Consequently, $\mathbf{WCEV}$ reduces to $\mathbf{WCEV-FMP}$, and an optimal solution to $\mathbf{WCEV-FMP}$ solves  $\mathbf{WCEV}$.  	
\end{proof}

\begin{remark}
$(i)$ We mention that, compared to the integration-based $\mathbf{WCEV}$ formulation, both \eqref{eq_limit_WCEV} and \eqref{eq_WCEV_opt} are more accessible. The finite nonlinear program $\mathbf{WCEV-FMP}$ certainly allows us to take advantage of many existing mathematical programming tools or results. For example, a strong nonlinear programming solver or algorithm can be readily used as an oracle to compute $\mathbf{WCEV-OPT}$ if the incorporation of $Q(\bf x)$ and $\psi_i$ is computationally friendly.  Indeed, even if the oracle is a fast approximation or heuristic one, i.e., the exactness of the derived $(P^*, \xi^*)$ cannot be guaranteed, it provides a basis for applying more sophisticated procedures for refinements. On the other hand, \eqref{eq_limit_WCEV} indicates that in general we can gradually augment $\bsm\xi(n)$ and associated $P(\bsm\xi(n))$ to approach $\mathbf{WCEV}$ arbitrarily.

\noindent $(ii)$ We would like to highlight the critical advantage of the primal representation demonstrated by \eqref{eq_limit_WCEV} and \eqref{eq_WCEV_opt}  in handling much more complex and general ambiguity sets. Note that Corollary \ref{cor_WCEV_limit} does not require $\mathcal{P}$ to be convex in $P$.  
It is very different from the existing duality-based methodologies to compute $\mathbf{WCEV}$ problem. Indeed, if $\mathcal{P}$ is a linear mixed integer set and an upper bound on the number of its constraints is known, it is viable to construct the corresponding $\mathbf{WCEV-FMP}$ to compute $\mathbf{WCEV}$. Certainly, solution procedures for mixed integer nonlinear programs are needed. Such a general situation is addressed in Section \ref{subsect_MIPAmbiguity}.  \qed 
\end{remark}

Nevertheless, we mention that $\mathbf{WCEV-FMP}$ is a challenging non-convex program, given the bilinear terms between $P$ and $Q(\xi)$ or $\psi_i(\xi)$ in \eqref{eq_WCEV_opt}. Actually, in the context of $\mathbf{2-Stg \ DRO}$, value function $Q(\xi)$ represents the recourse cost, which is a complex function of $\xi$ and renders $\mathbf{WCEV-FMP}$ computationally very intractable. Even just with first moment constraints, we observe that it could take an extremely long time for a state-of-the-art professional solver to solve small-scale instances. Hence, regardless of its simplicity and compactness, $\mathbf{WCEV-FMP}$ is still difficult. On the other hand, \eqref{eq_limit_WCEV} inspires us to develop computationally effective procedures to compute $\mathbf{WCEV}$ problem (including its variants).

\subsection{A Decomposition Algorithm for WCEV Problem}
\label{subsect_CG}
In the remainder of this paper, we assume, unless noted otherwise, the sufficient conditions presented in Proposition~\ref{prop_WCEV_OPT} to ensure the attainability of the WCEV. Rather than treating ambiguity set $\mathcal P$ as a whole set, we consider its sample space and the probability distribution in a separate fashion. Recall that  $\bsm \xi_n\equiv (\xi_1,\dots,\xi_n)$ represents a pool of discrete scenarios, and  $P(\bsm\xi_n)\equiv\big(p(\xi_1), \dots,p(\xi_n)\big)\in \mathcal{P}(\bsm\xi_n)\subseteq \cal P$ a discrete probability distribution and the subset of $\cal P$ defined on $\bsm \xi_n$. Also, the optimal value of an infeasible maximization problem is conventionally set to $-\infty$. The next result directly follows from the $\sigma$-algebra of $\Xi$ and Corollary \ref{cor_WCEV_limit}.
\begin{cor}
	For $\mathbf{WCEV}$ problem, i.e., $\displaystyle{\max_{P\in \mathcal P} E_P[Q(\bf \xi)]}$, we have
	\begin{equation}
	\label{eq_WCEV_reduction}
	\max_{ P\in \mathcal P} E_P[Q(\bf \xi)]\geq  \max_{\bsm\xi_n\subseteq \Xi}\max_{P(\bsm\xi_n)\in \mathcal{P}}\sum_{j=1}^n Q(\xi_j)p(\xi_j)   \geq \max_{P(\bsm\xi^0_n)\in \mathcal{P}(\bsm\xi^0_n)}\sum_{j=1}^nQ(\xi^0_j) p(\xi^0_j),
   \end{equation}
where $\bsm\xi^0_n\equiv\{\xi^0_1,\dots,\xi^0_n\}$ is a set of fixed scenarios. \hfill \qed
\end{cor} 
The non-convex program in the middle of \eqref{eq_WCEV_reduction} actually equals the WCEV if $n\geq m+1$. The rightmost optimization problem, although just providing a lower bound, is a very simple linear program (LP) in $P$. 
\begin{remark}
It is worth highlighting that the lower bound from that LP in \eqref{eq_WCEV_reduction}  is the strongest one we can have for given $\bsm\xi^0_n$. That is, if its strength with respect to $\displaystyle\max_{ P\in \mathcal P} E_P[Q(\bf \xi)]$ is weak, $\bsm\xi^0_n$ should be expanded by including additional nontrivial scenarios. Computationally, we can start with a small-sized $\bsm\xi^0_n$ and then gradually expand it for stronger lower bounds.   \qed
\end{remark}

 Actually, if the expansion of $\bsm\xi^0_n$ can be managed appropriately, the optimal value of that LP approaches  $\displaystyle\max_{P\in \mathcal P} E_P[Q(\bf \xi)]$ exactly or with an arbitrary accuracy, which corresponds to the limit operation presented in Corollary \ref{cor_WCEV_limit}. We mention that such an expansion process can be realized by customizing the  well-known column generation (CG) algorithm, a classical decomposition method proposed to solve large-scale LPs \cite{KanZal.1951,ford1958suggested,dantzig1960decomposition,gilmore1961linear}. Next, we first present the explicit form of the LP in \eqref{eq_WCEV_reduction}, referred to as the \textit{pricing master problem} (PMP) in the remainder of this paper.  
	\begin{align}
		 \label{eq_PMP_basic}
	\begin{split}
		\mathbf{PMP}: \ \underline \eta^*(\bsm\xi^0_n)= \max\Big\{&\sum_{j=1}^{n} Q(\xi^0_j) p(\xi^0_j): \sum_{j=1}^{n} p(\xi^0_j) = 1, \ 
		\sum_{j=1}^{n} \psi_i(\xi^0_j)p(\xi^0_j)\leq \gamma_i \ \ \forall i \in [m],\\
		& \ p(\xi^0_j)\geq 0\ \ \forall j \in [n]\Big\}.
	\end{split}
	\end{align}
Let $\alpha$ and $\bsm\beta\equiv[\beta_1,\dots,\beta_m]$ be  dual variables of its constraints, respectively. Supposing that $\mathbf{PMP}$ is feasible and its shadow prices are 
$(\alpha^*,\bsm\beta^*)$, the corresponding \textit{pricing subproblem} (PSP) to derive a new scenario with the largest reduced cost is 
\begin{align}
	\label{eq_PSP}
	\mathbf{PSP}:  \ v^*(\bsm\xi^0_n) = \max_{\xi \in \Xi} Q(\xi) - \alpha^* -\sum_{i=1}^m\psi_i(\xi)\beta^*_i.
\end{align}
Next, we provide an estimation on the strength of the lower bound derived from $\bsm\xi^0_n$. 
\begin{prop}
\label{prop_CG_basic}
Suppose that $\mathbf{PMP}$ is feasible and both $\mathbf{PMP}$ and $\mathbf{PSP}$ are solved to optimality. We have 
\begin{equation}
\underline{\eta}^*(\bsm \xi^0_n)  \leq \max_{ P\in \mathcal P} E_P[Q(\bf \xi)]\leq \underline{\eta}^*(\bsm\xi^0_n) +v^*(\bsm\xi^0_n). \tag*{\qed}
\end{equation}
\end{prop}
\begin{proof} Note that it is sufficient to prove the second inequality. 
	From $\mathbf{PSP}$ in \eqref{eq_PSP}, it can be seen that $\displaystyle{v^*(\bsm\xi^0_n)+\alpha^* \geq \max\limits_{\xi \in \Xi} Q(\xi) - \sum_{i=1}^m \psi_i(\xi) \beta_i^*}$, or equivalently
    \[v^*(\bsm\xi^0_n)+\alpha^* \geq Q(\xi) - \sum_{i=1}^m \psi_i(\xi) \beta_i^* \ \ \forall \xi\in\Xi.\]
	On the other hand, noting that \eqref{eq_WCEV} is an infinite-column linear program with the strong duality \cite{shapiro2001duality}, we have its dual problem 
\begin{align*}
    \min\Big\{\alpha + \sum_{i=1}^m \beta_i \gamma_i: \alpha +\sum_{i=1}^m\psi_i(\xi)\beta_i \geq Q(\xi) \ \ \forall \xi\in \Xi, \ \alpha \ \text{free}, \ \beta_i\geq 0\  \ \forall i \in [m]\Big\}.
\end{align*}
Clearly, by setting $\alpha=v^*(\bsm\xi^0_n)+\alpha^*$ and $\beta_i=\beta_i^*$ for all $i$, we obtain a feasible solution to the dual problem. Hence, we have 
\[\max_{ P\in \mathcal P} E_P[Q(\bf \xi)] \leq v^*(\bsm\xi^0_n)+\alpha^*+\sum_{i}\beta^*\gamma_i=v^*(\bsm\xi^0_n)+\underline{\eta}^*(\bsm\xi^0_n),\]
where the last equality follows from the strong duality of $\mathbf{PMP}$.
\end{proof}

According to Proposition \ref{prop_CG_basic}, for a given $\bsm\xi^0_n$ (and hence $\mathbf{PMP}$ and $\mathbf{PSP}$ are given), if the optimal value of $\mathbf{PSP}$ is $0$, $\underline{\eta}^*(\bsm\xi^0_n)$ equals the WCEV. Otherwise, denoting $\mathbf{PSP}$'s optimal solution by $\xi^*$, we can  augment $\bsm\xi^0_n$ by including $\xi^*$ and recompute $\mathbf{PMP}$. This process is repeated until the reduced price becomes sufficiently small, which is the basic idea of CG algorithm.
For simplicity, we refer to  the approach directly constructing $\mathbf{WCEV-FMP}$ and solving it by some stand-alone method(s)  as \textit{Oracle-1}, and one using an iterative procedure as \textit{Oracle-2}. In the context of this paper, \textit{Oracle-2 } is just the following customized CG algorithm. Note that we do not include subscript $n$, unless we need to track the number of scenarios in $\bsm\xi^0$.

\refstepcounter{alg} 
\noindent\begin{minipage}{\textwidth}
    \makebox[0.92\textwidth]{\hrulefill}\\[-0.8em]
    \noindent\textit{Oracle-2:} The CG Algorithm to Compute the WCEV \label{alg_CG}
\end{minipage}
\begin{algorithmic}
\State \textbf{Step 1}  Given initial $\bsm\xi^0$ and optimality tolerance $\varepsilon$, set the iteration counter $k=1$. 
\State\textbf{Step 2} Solve $\mathbf{PMP}$  to derive optimal value $\underline \eta^{k*}(\bsm\xi^0)$ and shadow price $(\alpha^*, \bsm\beta^*)$.
\State\textbf{Step 3} Solve \textbf{PSP}  to derive its optimal solution $\bf \xi^*$ and optimal value $v^*(\bsm\xi^0)$.
\State\textbf{Step 4} \parbox[t]{\dimexpr\linewidth-\algorithmicindent}{If $v^*(\bsm\xi^0_n) \leq \varepsilon$, report  $\underline \eta^{k*}(\bsm\xi^0)$ as the optimal value of \eqref{eq_WCEV} and terminate. Otherwise, update $\bsm\xi^0 = \bsm\xi^0 \cup \{\bf \xi^*\}$ and $k=k+1$, and go to \textbf{Step 2}.}
\end{algorithmic}
\noindent\makebox[0.92\textwidth]{\hrulefill}

We mention that, as an algorithm from the primal perspective to compute $\mathbf{WCEV}$ problem, it generates a set of discrete scenarios and their respective probabilities.  As discussed next, this algorithm brings us several important features that can be further explored.

\begin{remark}    
$(i)$ Similar to $\mathbf{WCEV-FMP}$, this CG algorithm mainly works in the primal space and is general in handling different ambiguity sets. Even if $\mathcal P$ is a mixed integer set, for which the strong duality-based approaches do not work, it is feasible to compute the WCEV by the extension of \textit{Oracle-2}, i.e., the well-established Brand-and-Price (B\&P) algorithm. Alternatively, we propose a new CG variant  in Section \ref{subsect_MIPAmbiguity} to handle such ambiguity sets without  developing B\&P procedures.  \\
\noindent$(ii)$ At the termination, we also have 
	\begin{align}
	\label{eq_WCEV_primal_opt}
	\begin{split}
	\underline{\eta}^{k*}(\bsm\xi^0) = \max\Big\{\sum_{\xi\in \bsm{\xi}^0} Q(\xi) p(\xi): (p_{\xi})_{\xi\in \bsm\xi^0}\in \mathcal{P}\Big\}.
	\end{split}
\end{align}
By Proposition \ref{prop_CG_basic}, the WCEV is bounded by $\ubar \eta^{k*}(\bsm\xi^0)+\varepsilon$ at termination. Instead of utilizing a time-consuming oracle to solve $\mathbf{PSP}$ exactly, this observation allows us to  develop a fast approximation algorithm for complex $\mathbf{PSP}$ and hence for the WCEV, as long as its approximation bound is available.\\
$(iii)$ Compared to directly computing the nonlinear program $\mathbf{WCEV-FMP}$ by a professional solver, \textit{Oracle-2} demonstrates a superior capacity that is generally faster by multiple orders of magnitude. Actually, it  is rather a vanilla version of CG. More advanced implementation techniques, after customization to fit into the  DRO context, should be able to help us accelerate the computation or the convergence. 	\\
$(iv)$ Although \textit{Oracle-2} is currently faster by several orders of magnitude, we anticipate that this situation may change with the development of specialized techniques to improve \textit{Oracle-1}. Yet, the applicability of \textit{Oracle-2} is less restrictive, as it does not depends on the number of scenarios that are with non-zero probabilities as represented  in Proposition \ref{prop_WCEV_OPT}.   \qed
\end{remark}

 According to the CG literature, it is not restrictive to assume that $\mathbf{PMP}$ is feasible, as a modified CG based on Farka's Lemma can generate new scenarios (i.e., columns) to ensure it to be feasible. Alternatively, as showed in Section \ref{subsect_CCGDRO}, the algorithm framework for $\mathbf{2-Stg \ DRO}$ can be used to address this issue of infeasible $\mathbf{PMP}$. In the next subsection, we consider the convergence issue of \textit{Oracle-2} and its computational complexity. 
\subsection{Convergence and Complexity of \textit{Oracle-2}}
\label{subsect_CG_convergence}
Without loss of generality, we assume $\bsm\xi_0=\emptyset$ at the initialization.  Since \textit{Oracle-2} derives optimal $\xi^*$ from $\mathbf{PSP}$ and updates $\bsm\xi^0 = \bsm\xi^0 \cup \xi^*$ in each iteration, a sequence $\big\{\underline \eta^{k*}(\bsm\xi^0)\big\}_k$, consisting of optimal values of $\mathbf{PMP}$s, can be found. As it can be easily recovered, we omit $\bsm\xi_0$ when appropriate to simplify our arguments, unless otherwise stated. The next theorem reveals the relationship between the limit of $\big\{\underline \eta^{k*}\big\}_k$ and the WCEV.

\begin{lem}
Assuming that $\displaystyle{\max_{ P\in \mathcal P} E_P[Q(\bf \xi)]}$ exists and an $\varepsilon$-optimal solution of $\mathbf{WCEV}$ problem can be derived in finite iterations for $\varepsilon>0$, $\big\{\underline \eta^{k*}\big\}_k$  converges to $\displaystyle{\max_{ P\in \mathcal P} E_P[Q(\bf \xi)]}$.\hfill \qed
\end{lem}
\begin{proof}
 Consider the simple case where \textit{Oracle-2} terminates with $v^*=0$ in the $k_t$-th iteration and hence it involves $k_t-1$ different $\xi$'s. According to Proposition \ref{prop_CG_basic}, we obtain the optimal value of $\mathbf{WCEV}$ problem, which leads the conclusion. Next, we consider the other case where an $\varepsilon-$optimal solution of $\mathbf{WCEV}$ problem can be retrieved in finite iterations.

According to Corollary \ref{cor_WCEV_limit}, for any given $\varepsilon > 0$, there exists an integer $K$ such that
\[\max_{ P\in \mathcal P} E_P[Q(\bf \xi)] - \varepsilon < \underline{\eta}^{K*} \leq \max_{ P\in \mathcal P} E_P[Q(\bf \xi)].\] Otherwise $\max_{ P\in \mathcal P} E_P[Q(\bf \xi)] - \varepsilon$ would be an upper bound of $\{\underline{\eta}^{k*}\}_{k}$, which contradicts the derivation of an $\varepsilon-$optimal solution by \textit{Oracle-2}. Because $\varepsilon$ can be arbitrarily small,  $\{\underline{\eta}^{k*}\}_{k}$ is increasing with $k$, and because of the monotone convergence theorem, $$\lim\limits_{k\rightarrow +\infty}\underline{\eta}^{k*} = \max_{ P\in \mathcal P} E_P[Q(\bf \xi)],$$
which is the expected conclusion. 
\end{proof}

The next theorem shows that we actually can always find an optimal or  $\varepsilon$-optimal solution of $\mathbf{WCEV}$ problem  in finite iterations under some mild conditions. Before that, consider a continuous function $f:\ \cal X \times \cal Y \to \mathbb R$, for $\bf y \in \cal Y$ and for any $\varepsilon > 0$, if there exist a $\delta >0$, such that 
\[|f(\bf x_1, \bf y) - f(\bf x_2, \bf y)| \leq \varepsilon\ \ \mbox{if} \ \ \Vert \bf x_1 - \bf x_2\Vert \leq \delta \ \ \mbox{for} \ \bf x_1, \bf x_2\in \cal X,\]
we say that $f(\bf x, \bf y)$ is continuous, uniformly with respect to $\bf x$ for $\bf y\in \cal Y$.
\begin{theorem}
\label{CG-epsilon}
Assume that the reduced cost function, i.e., $\displaystyle{r(\xi; \alpha,\bsm \beta)\equiv Q(\xi) - \alpha - \sum_{i=1}^m \psi_i(\xi) \beta_i}$, is continuous, uniformly with respect to $\xi$ over $\Xi$ for $\bsm \beta\geq \bf 0$. Then,  \textit{Oracle-2} returns  an optimal or $\varepsilon$-optimal solution of $\mathbf{WCEV}$ problem in finite iterations.\hfill \qed
\end{theorem}
\begin{proof}
 According to the assumption, for any $\varepsilon > 0$, there exists a $\delta > 0$ such that \[
	|r(\xi_1;\alpha,\bsm\beta) - r(\xi_2;\alpha,\bsm\beta)| \leq \varepsilon \ \ \mbox{if} \ \ \Vert\xi_1 - \xi_2\Vert \leq \delta \ \ \mbox{for} \ \xi_1, \xi_2\in \Xi, \bsm \beta\geq \bf 0.\]
Let $B(\xi_1, \delta)$ denote the closed ball with radius $\delta$ at $\xi_1\in \Xi$, i.e.,  $B(\xi_1,\delta)=\{\xi_2\in \Xi, \Vert\xi_1 - \xi_2\Vert \leq \delta\}$. We consider the following claim. 

\noindent \textit{Claim}: In any particular iteration with associated $\bsm \xi^0$, \textit{Oracle-2} either solves $\mathbf{WCEV}$ problem to an optimal or $\varepsilon$-optimal solution, or produces a new scenario $\xi^*$ that is not contained in ball $B(\xi, \delta)$ for any $\xi \in \bsm \xi^0$.\\
\textit{Proof of Claim:} Recall $(\alpha^*,\bsm\beta^*)$ denotes  the shadow price obtained from computing $\mathbf{PMP}$. We prove this claim by contradiction. 

Suppose that scenario $\xi^*$ with $r(\xi^*; \alpha^*,\bsm\beta^*)$ is identified by $\mathbf{PSP}$. Note that it is sufficient to assume that 
 $r(\xi^*; \alpha^*,\bsm\beta^*)> \varepsilon$, since otherwise  it is an optimal or $\varepsilon$-optimal solution of $\mathbf{WCEV}$ problem by Proposition \ref{prop_CG_basic}. Assume further that it is contained in ball $B(\xi', \delta)$ for some $\xi' \in \bsm \xi^0$. We have   
\begin{align*}
	v^*(\bsm\xi^0) = &  r(\xi^*;\alpha^*, \bsm\beta^*)\\
	= &  r(\xi^*;\alpha^*, \bsm\beta^*) - r(\xi';\alpha^*, \bsm\beta^*)+r(\xi';\alpha^*, \bsm\beta^*) \\
	= & |r(\xi^*;\alpha^*, \bsm\beta^*) - r(\xi';\alpha^*, \bsm\beta^*)|+ r(\xi';\alpha^*, \bsm\beta^*) \\
	\leq & \varepsilon + r(\xi';\alpha^*, \bsm\beta^*)\\
	\leq & \varepsilon.
\end{align*}
The second equation holds due to an identity transformation, the third follows from the fact  that $r(\xi';\alpha^*,\bsm\beta^*)\leq  0$ and the first inequality follows the definition of  uniformly continuous function. The last inequality, which is valid due to the fact that $r(\xi;\alpha^*,\bsm\beta^*)\leq  0$ for  $\xi \in \bsm\xi^0$, clearly contradicts to our first assumption.  

With the aforementioned contradiction, we can conclude that either $\bsm \xi^0$ (and the associated probabilities) is an  optimal or $\varepsilon$-optimal solution of $\mathbf{WCEV}$ problem, or $\xi^*$ is not contained in ball $B(\xi, \delta)$ for any $\xi \in \bsm \xi^0$. \qed

It simply follows from the claim that the distance between any two scenarios in $\Xi^0$ is more than $\delta$, indicating the two balls centered at them with $\frac{\delta}{2}$-radius are disjoint. Let $V(\frac{\delta}{2})$ denote the volume of a such ball. 
Given that $\Xi$ is compact, it is clear that $\Xi$ is contained in a ball centered at some $\xi^0\in \Xi$ with radius $\hat d$. Consequently, $B(\xi^0, \hat d+\frac{\delta}2)$ is compact, and its volume,  denoted by $V( \hat d+\frac{\delta}2)$, is finite.  Because  $V( \hat d+\frac{\delta}2)/V(\frac{\delta}{2})$ is finite, it follows that an optimal or  $\varepsilon$-optimal solution of $\mathbf{WCEV}$ problem can be found by \textit{Oracle-2} within a finite number of iterations, bounded by $\left\lceil V( \hat d+\frac{\delta}2)/V(\frac{\delta}{2})\right\rceil $.
\end{proof}

Next, we consider some special case where an optimal solution is guaranteed to be obtained within finite iterations.  
\begin{cor}
\label{CG-exact}
Assume that  $\Xi$ is a polytope, $Q(\xi)$ is convex and $\psi_i(\xi)$ are concave over $\Xi$ for $i\in [m]$. Let $\mathsf{X\!V}(\Xi)$ be the set of extreme points of $\Xi$. Then,  \textit{Oracle-2} terminates with an optimal solution to  $\mathbf{WCEV}$ problem within $|\mathsf{X\!V}(\Xi)|$ iterations.\hfill \qed
\end{cor}
\begin{proof}
We note that, under the aforementioned assumptions,  
an optimal solution  to $\mathbf{PSP}$ in \eqref{eq_PSP} is an extreme point of $\Xi$, i.e., it belongs to $\mathsf{X\!V}(\Xi)$. Moreover, it can be seen from the proof of Theorem \ref{CG-epsilon} that a new extreme point of $\Xi$, which is not contained in $\bsm\xi^0$, will be identified and then used to update  $\bsm\xi^0$ in each iteration before \textit{Oracle-2} terminates.

 Given that $\mathsf{X\!V}(\Xi)$ is a finite set, the expected conclusion follows. 
\end{proof}

\begin{remark}
	$(i)$ Actually, results of Theorem \ref{CG-epsilon} can be extended to handle some discontinuous  $r(\xi; \alpha, \bsm \beta)$. For example,  as long as $r(\xi; \alpha, \bsm \beta)$ is uniformly continuous within each subset of a finite partition of $\Xi$ for $\bsm \beta \geq \bf 0$,  where Proposition \ref{prop_WCEV_OPT} does not hold, an $\varepsilon$-optimal solution for $\mathbf{WCEV}$ problem can still be found within finite iterations.\\
	$(ii)$ As discussed,  $\psi_i(\cdot)$ functions are continuous over $\Xi$ for ambiguity sets defined by moment inequalities and Wasserstein metric, which render Theorem  \ref{CG-epsilon} to hold when $Q(\cdot)$ guarantees $r(\xi; \alpha, \bsm \beta)$ uniformly continuous over $\Xi$ for $\bsm \beta \geq \bf 0$. Moreover, in the case where $\Xi$ is a polytope and $\cal P$ involves the first moment inequalities, the second moment inequalities subject to lower bounds,  or the Wasserstein metric restriction defined with $L_1$ norm, optimal $\xi^*$ of $\mathbf{PSP}$, which is to maximize a convex function, can always be obtained at some extreme point of $\Xi$. According to Corollary \ref{CG-exact}, an optimal solution for $\mathbf{WCEV}$ problem can be obtained within finite iterations. \qed
\end{remark}

In the next section, we seek to compute the whole $\mathbf{2-Stg \ DRO}$, and present a decomposition algorithm that incorporates those primal oracles for $\mathbf{WCEV}$ problem.

\section{From Computing the WCEV to Solving $\bf{2-Stg \ DRO}$}
\label{sect_CCG_DRO}
In this section, we present a generally applicable and computationally strong solution scheme to solve $\bf{2-Stg \ DRO}$, along with theoretical analyses on its strength,  convergence and computational complexity.  Compared to all existing methods in the literature, it enhances our understanding on DRO and significantly expands and improves our ability to solve practical instances, particularly those with the challenging infeasible recourse issue.     

We highlight that the basic intuitions behind our algorithm design are: $(i)$ For a set of given scenarios, the worst-case probability distribution can be readily determined by solving an LP; $(ii)$ Together with the first stage decision making, a simple bilevel optimization problem can be constructed to obtain a lower-bound approximation of the original $\bf{2-Stg \ DRO}$; $(iii)$ If such lower-bound approximation is not satisfactory, additional scenarios needed to determine the WCEV can be identified by solving $\bf{WCEV}$ problem and then are employed to expand that scenario set through the framework of C\&CG.

\subsection{Integrating Primal Oracles within C\&CG}
\label{subsect_CCGDRO}
As noted, different from the current mainstream strategies on solving DRO, our algorithm development takes the primal perspective: it integrates the oracles for $\bf{WCEV}$ problem discussed in the previous section within the C\&CG framework, resulting a simple and intuitive algorithmic  structure that helps  its adoption among practitioners.  

We say that two formulations are \textit{equivalent} to each other if they share the same optimal value, and one's optimal first-stage solution is also optimal to the other one.  
With this definition  and by the results in Corollary \ref{cor_WCEV_limit} and Proposition \ref{prop_WCEV_OPT}, two equivalent reformulations can be obtained  for $\bf{2-Stg \ DRO}$.

\begin{prop}
 For $\bf{2-Stg \ DRO}$ in \eqref{eq_2stgDRO}, we have 
\begin{subequations}
\label{eq_WCEV_2RO}
\begin{align*}
	w^* =&\min\Big\{f_1(\bf x)+\eta:  \bf x\in \cal X, \ \eta \geq \lim_{n \rightarrow + \infty} \max_{\xi_j(n)\in \Xi, j\in [n]}  \max\{\sum_{j=1}^n Q(\bf x, \xi_j(n))p_j: P\in \cal{P}\} \Big\} \\
	=&\min\Big\{f_1(\bf x)+\eta:  \bf x\in \cal X, \ \eta\geq \max\big\{\sum_{j=1}^{m+1} Q(\bf x, \xi_j)p_j: P\in \cal{P}, \xi_j\in \Xi \ \ \forall j\in[m+1] \big\}\Big\}.\tag*{\qed}
\end{align*}
\end{subequations}
\end{prop}
Both equivalent reformulations can be treated as bilevel optimization problems with a non-convex lower-level maximization problem. Note that no closed-form is available to characterize the optimal solution set of such lower-level problem to build a single-level reformulation further. Indeed, they can easily yield lower bounds to $w^*$ if we reduce the lower-level problems to simpler ones, e.g., LPs. 
Let $\bf x^*$ denote a fixed first-stage solution. The resulting lower bound we obtain might not be strong enough to substantiate $\bf x^*$'s quality, which may demand us to revise our simplification.   Next, we show that, by employing  the oracles developed for $\bf{WCEV}$ problem and the C\&CG scheme, a complete algorithmic scheme for $\bf{2-Stg \ DRO}$ is developed, which ensures  an optimal and feasible $\bf x^*$ can be produced.

\subsubsection{Addressing the Challenge of Infeasible Recourse}
By convention, we set $Q(\bf x,\xi)=+\infty$ if $\cal Y(\bf x, \xi)=\emptyset$, i.e., the recourse problem in \eqref{eq_recourse} is infeasible. Also, let $\cal E$ denote the event when $\cal Y(\bf x, \xi)=\emptyset$ (or equivalently $Q(\bf x,\cdot)=+\infty$). Then, by extending the feasibility definition from the deterministic context, we say a first-stage decision $\bf x$ is \textit{almost surely} feasible to $\bf{2-Stg \ DRO}$ if $\displaystyle\sup_{P\in \cal P}P(\cal E)=0$, i.e., the probability of this event $\cal E$ is 0 for $P\in \cal P$. To implement this concept in computation, we introduce a supporting problem, whose optimal value corresponds to the status of $\cal Y(\bf x, \xi)$. Specifically, we assume that $\cal Y(\bf x, \xi)\equiv\big\{\bf y: \bf g(\bf x,\xi, \bf y)\geq \bf 0\big\}$ and define 
\begin{align}
	\cal{\tilde Y}(\bf x, \xi)\equiv\big\{(\bf y, \bf{\tilde y}): \bf g(\bf x,\xi, \bf y)+\bf{\tilde y}\geq \bf 0,\ \bf{\tilde y}\geq \bf 0\big\},
\end{align}
with artificial variable $\bf{\tilde y}$. Note that $\cal{\tilde Y}(\bf x, \xi)\neq \emptyset$ regardless of $\bf x$ and $\xi$.
	With $\lVert \cdot \rVert_{1}$ denoting the  $L_1$ norm, the supporting problem is 
\begin{align}
	\label{eq_support_feasibility}
\tilde Q_f(\bf x,\xi)=\min\Big\{\lVert\bf{\tilde y}\rVert_{1}: (\bf y,\bf{\tilde y})\in \cal{\tilde Y}(\bf x,\xi)\Big\}.
\end{align}
Since \eqref{eq_support_feasibility} is to minimize the $L_1$ norm of $\bf{\tilde y}$, the next result follows easily. 
\begin{lem}
\label{lem_feasibility_WCEV}
$\tilde Q_f(\bf x,\xi)>0$ iff $Q(\bf x,\xi)=+\infty$. Moreover, $\bf x$ is almost surely feasible to  $\bf{2-Stg \ DRO}$ iff the optimal value of the following problem equals 0. 
\begin{align}
\label{eq_feasibility_subprob}
\bf{WCEV}(\it F): \ 	\displaystyle\max_{P\in \cal P}E_P[\tilde Q_f(\bf x, \xi)]\equiv \max\Big\{\int_{\Xi} \tilde Q_f(\bf x, \xi)P(d\xi): P\in \cal{P}\Big\}.
\end{align}
\end{lem}
As a result, the worst-case expected value of $\tilde Q_f(\bf x,\xi)$ can be used to verify the (almost surely) feasibility of $\bf x$, and hence is denoted by $\bf{WCEV}(\it F)$. To keep our notation consistent, the one presented in \eqref{eq_WCEV}, after replacing $Q(\xi)$ by $Q(\bf x, \xi)$, is interchangeably referred to as $\bf{WCEV}(\it O)$ as it computes the worst-case expected recourse cost. We mention that when a mixed integer ambiguity set, denoted by $\cal P^I$, is adopted as the ambiguity set, Lemma~\ref{lem_feasibility_WCEV} still holds. Actually, both  $\bf{WCEV}(\it O)$ and  $\bf{WCEV}(\it F)$, as well as their extensions on $\cal P^I$, can be computed by the primal oracles developed in the previous section.

\begin{remark}
When the optimal value of $\bf{WCEV}(\it F)$ equals 0, it has an intuitive interpretation. That is, the probability over the subset of $\Xi$ satisfying $\tilde Q_f(\bf x,\xi)=0$ equals 1 and the probability over its complement 0 for any distribution within $\cal P$. When the optimal value of $\bf{WCEV}(\it F)$ is larger than 0, it means that some scenarios become infeasible for the given $\bf x$ and they have positive probabilities under some distribution in $\cal P$. For the latter case,  the scenario set generated by a primal oracle should be used within C\&CG in a way such that the current $\bf x$ is cut off, i.e., being excluded from future considerations. Note that in the remainder of this paper, we say $\bf x$ is feasible means that it is almost surely feasible.  \qed
\end{remark}

\subsubsection{Primary Components of the C\&CG Method for $\bf{2-Stg \ DRO}$}
The C\&CG method for $\bf{2-Stg \ DRO}$, which is referred to as C\&CG-DRO, involves a master problem and two subproblems. The master problem is a relaxation to $\bf{2-Stg \ DRO}$ and yields a lower bound. Two subproblems, which are for feasibility and optimality respectively, help us strengthen the master problem and derive an upper bound.   

\noindent (I): \textbf{The Master Problem of C\&CG-DRO}\\
 The master problem is defined on $\bsm{\widehat \xi}^o$ and $\bsm{\widehat \xi}^f$, two sets of fixed scenarios in $\Xi$. To differentiate it from the master problem of \textit{Oracle-2}, we call it the \textit{main master problem} ($\bf{MMP}$). In the following we first present a form developed using big-M technique, which is rather intuitive and interpretable. Then, with deep insights, we derive a big-M-free one that is compact and more rigorous. Hence, they are referred to as $\bf{MMP}_1$ and $\bf{MMP}_2$, respectively.  
\begin{subequations}
\label{eq_MMP1_CCG}
	\begin{align}
		\textbf{MMP}_1 : \ubar{w} = \min\limits_{\bf x \in \cal{X}} &f_1(\bf x) + \eta\\
		& \eta \geq \max\Big\{\sum_{\xi\in \bsm{\widehat \xi}^o} \eta^o_{\xi}p^o_{\xi}: (p^o_{\xi})_{\xi\in \bsm{\widehat \xi}^o}\in \cal P\Big\} \label{eq_WCEV_MMP1}\\
		& \Big\{\eta^o_{\xi}= f_2(\bf x, \xi,\bf y_{\xi})+\bf M\lVert\bf{\tilde y}^o_{\xi}\rVert_{1}, \ (\bf y_{\xi},\tilde {\bf y}^o_{\xi})\in \tilde{\cal Y}(\bf x,\xi)\Big\} \ \ \forall \xi\in \bsm{\widehat \xi}^o \label{eq_MMP_Oscenarios}\\
		& 0\geq \max\Big\{\sum_{\xi\in \bsm{\widehat \xi}^f} \eta^f_{\xi}p^f_{\xi}: (p^f_{\xi})_{\xi\in \bsm{\widehat \xi}^f}\in \cal P\Big\} \label{eq_WCEV_MMP2}\\
		& \Big\{\eta^f_{\xi}= \lVert\bf{\tilde y}^f_{\xi}\rVert_{1}, \ (\bf y_{\xi},\tilde{\bf y}^f_{\xi})\in \tilde{\cal Y}(\bf x,\xi)\Big\} \ \ \forall \xi\in \bsm{\widehat\xi}^f \label{eq_MMP_Fscenarios}
	\end{align}
\end{subequations}
We mention  that in $\eqref{eq_MMP_Oscenarios}$ $\tilde{\cal Y}(\cdot,\cdot)$ is employed, instead of original $\cal Y(\cdot,\cdot)$, and the sum of $f_2$ and big-M penalized $\lVert\bf{\tilde y}^o_{\xi}\rVert_{1}$ is assigned to $\eta^o_{\xi}$. It is worth highlighting that these components are fundamental in solving $\bf{2-Stg \ DRO}$ when the the feasibility issue arises in the recourse problem.  Consider a scenario $\xi^o$ in $\bsm{\widehat\xi}^o$. On one hand,  if $p^o_{\xi}$ equals $0$, i.e., $\xi^o$ virtually does not occur, it is not necessary to have the corresponding $\cal Y(\cdot,\xi^o)$ being non-empty. Or, alternatively the choice of $\bf x$ should not be affected by $\cal Y(\cdot,\xi^o)$. On the other hand, if $p^o_{\xi}$ is positive, indicating the occurrence of $\xi^o$ cannot be ruled out,  $\bf x$ must render $\cal Y(\cdot,\xi^o)$ non-empty to ensure a feasible recourse. Or, alternatively those causing $\cal Y(\cdot,\xi^o)$ empty cannot be in any optimal solution of $\bf{MMP}_1$. Those two points and the logic behind, which are the nature of DRO, can be well achieved by making use of artificial variables $\bf{\tilde y}^o_{\xi}$ and big-M coefficient $\bf M$ in $\eqref{eq_MMP_Oscenarios}$. 

 Note that  when $p^o_{\xi}=0$, $\eta$ in $\eqref{eq_WCEV_MMP1}$ is not affected by $\eta^o_{\xi}$,   and $\cal{\tilde Y}(\bf x, \xi^o)\neq \emptyset$ for any $\bf x$. Hence, we can conclude that $\xi^o$ has no impact on $\bf{MMP}_1$. Also, when $p^o_{\xi}>0$, the  big-M penalty term in $\eta^o_{\xi}$ will drive $\bf{MMP}_1$ to select an $\bf x$ such that $\bf{\tilde y}^o_{\xi}$ can be $0$, i.e., $\cal Y(\bf x,\xi^o)\neq \emptyset$. As a result, $\eqref{eq_MMP_Oscenarios}$ will disqualify $\bf x$ in any optimal solution of $\bf{MMP}_1$ if it causes the recourse problem of $\xi^o$ infeasible. Certainly, if the recourse problem is assumed to be feasible for any $\bf x\in \cal X$ and $\xi\in \Xi$, it is not needed to introduce variable $\bf{\tilde y}$ and set $\cal{\tilde Y}$.

\begin{remark}
 As the two inequalities in \eqref{eq_WCEV_MMP1} and \eqref{eq_WCEV_MMP2} refine the feasible set of $\bf x$ from optimality and feasibility perspectives, we refer to them as \textit{optimality} and \textit{feasibility} cutting planes, and $\widehat{\bsm \xi}^f$ and $\widehat{\bsm\xi}^o$ feasibility and optimality sets, respectively.  Note that different from most of classical cutting plane methods that generate new and independent cutting planes over iterations, $\bf{MMP}_1$ strengthens either the \textit{optimality} or the \textit{feasibility} cutting plane based on the augmented $\widehat{\bsm{\xi}}^o$ or $\bsm{\widehat \xi}^f$ in every iteration. \qed
\end{remark}   

By further investigating the critical interaction between a scenario and the associated probability, we actually can leverage it to eliminate the reliance on big-M. Similar to the unification idea presented in \cite{ZENG2013457} and \cite{zeng2022two}, the two sets of scenarios, for optimality and feasibility respectively, can be merged to have a unified set.  By doing so, we let $\widehat{\bsm \xi}=\widehat{\bsm \xi^f}\cup\widehat{\bsm \xi}^o$, through which we can build the following alternative formulation of $\bf{MMP}_1$. 
\begin{subequations}
	\label{eq_MMP2-W}
	\begin{align}
		\bf{MMP}_2 : \ubar{w} = \min\limits_{\bf x \in \cal{X}} &f_1(\bf x) + \eta\\
		& \eta \geq \max\Big\{\sum_{\xi\in \bsm{\widehat \xi}} \eta^o_{\xi}p^o_{\xi}: (p^o_{\xi})_{\xi\in \bsm{\widehat \xi}}\in \cal P\Big\} \label{eq_WCEV_MMP12}\\
		& \eta^o_{\xi}= f_2(\bf x, \xi,\bf y_{\xi})\ \ \forall \xi\in \bsm{\widehat \xi} \label{eq_MMP_Oscenarios2}\\
		& 0\geq \max\Big\{\sum_{\xi\in \bsm{\widehat \xi}} \eta^f_{\xi}p^f_{\xi}: (p^f_{\xi})_{\xi\in \bsm{\widehat \xi}}\in \cal P\Big\} \label{eq_WCEV_MMP22}\\
		& \eta^f_{\xi}= \lVert\bf{\tilde y}_{\xi}\rVert_{1} \ \ \forall \xi\in \bsm{\widehat\xi} \label{eq_MMP_Fscenarios2} \\
		& (\bf y_{\xi},\tilde{\bf y}_{\xi})\in \tilde{\cal Y}(\bf x,\xi) \ \ \forall \xi\in \bsm{\widehat\xi} \label{eq_MMP_yrecourse_relax}
	\end{align}
\end{subequations}
\begin{prop}
\label{prop_MMP2}
$\bf{MMP}_2$ is a valid relaxation and provides a lower bound to $\bf{2-Stg \ DRO}$. 
\end{prop}
\begin{proof}
	We first extend $\bf{MMP}_1$ by considering $\bm{\widehat\xi}$ that leads to the following auxiliary formulation.
	\begin{subequations}
		\begin{align}
			\bf{MMP}^a_1 : \ubar{w} = \min\limits_{\bf x \in \cal{X}} &f_1(\bf x) + \eta\\
			& \eta \geq \max\Big\{\sum_{\xi\in \bsm{\widehat \xi}} \eta^o_{\xi}p^o_{\xi}: (p^o_{\xi})_{\xi\in \bsm{\widehat \xi}}\in \cal P\Big\} \label{eq_MMP3_WCEV_o}\\
			& \Big\{\eta^o_{\xi}= f_2(\bf x, \xi,\bf y^o_{\xi}) + \bf{M}\lVert\bf{\tilde y}^o_{\xi}\rVert_{1}, \ (\bf y^o_{\xi},\tilde{\bf y}^o_{\xi})\in \tilde{\cal Y}(\bf x,\xi)\Big\} \ \ \forall \xi\in \bsm{\widehat\xi}\label{eq_MMP3_WCEV_c}\\
			& 0\geq \max\Big\{\sum_{\xi\in \bsm{\widehat \xi}} \eta^f_{\xi}p^f_{\xi}: (p^f_{\xi})_{\xi\in \bsm{\widehat \xi}}\in \cal P\Big\} \label{eq_MMP3_yrecourse_o}\\
			& \Big\{\eta^f_{\xi}= \lVert\bf{\tilde y}^f_{\xi}\rVert_{1},\ (\bf y^f_{\xi},\tilde{\bf y}^f_{\xi})\in \tilde{\cal Y}(\bf x,\xi)\Big\} \ \ \forall \xi\in \bsm{\widehat\xi}\label{eq_MMP3_yrecourse_c}
		\end{align}
	\end{subequations}
	Since $\widehat{\bsm \xi}=\widehat{\bsm \xi^f}\cup\widehat{\bsm \xi}^o \subseteq \Xi$,  $\bf{MMP}^a_1$ is a valid  relaxation that actually is stronger than $\bf{MMP}_1$.
	
	For $\bf{MMP}^a_1$, let $\{p_{\xi}^{o*}\}_{\xi\in \bsm{\widehat \xi}}$ denote an optimal solution to the LP in \eqref{eq_MMP3_WCEV_o}. Considering some $\xi'\in \bsm{\widehat \xi}$, we note that \eqref{eq_MMP3_yrecourse_o} and \eqref{eq_MMP3_yrecourse_c} will drive $\bf{\tilde y}^o_{\xi'}$ to be 0 if $p_{\xi'}^{o*}>0$. If this is not the case, we let $\big(p_{\xi'}^f, \bf{\tilde y}^f_{\xi'}\big) = \big(p_{\xi'}^{o*}, \bf{\tilde y}^o_{\xi'}\big)$. Then,  we have 
	$$\max\Big\{\sum_{\xi\in \bsm{\widehat \xi}} \eta^f_{\xi}p^f_{\xi}: (p^f_{\xi})_{\xi\in \bsm{\widehat \xi}}\in \cal P\Big\} \geq \eta_{\xi'}^fp_{\xi'}^f = \lVert\bf{\tilde y}^f_{\xi'}\rVert_{1}p_{\xi'}^f> 0,$$ 
	which is contradictory to inequality \eqref{eq_MMP3_yrecourse_o}. Given the arbitrarity of $\xi'$,  we can simply remove $\bf{M}$ in \eqref{eq_MMP3_WCEV_c} without affecting the optimality of $\bf{MMP}^a_1$. The updated \eqref{eq_MMP3_WCEV_c} is
	\begin{align}
		\label{eq_updated_no_MMP_bigM}
		\Big\{\eta^o_{\xi}= f_2(\bf x, \xi,\bf y^o_{\xi}), \ (\bf y^o_{\xi},\tilde{\bf y}^o_{\xi})\in \tilde{\cal Y}(\bf x,\xi)\Big\} \ \ \forall \xi\in \bsm{\widehat\xi}.
	\end{align}

	Next, we prove that two sets of $(\bf y_{\xi},\tilde{\bf y}_{\xi})$ variables are not necessary. Let $\Big(\cal P^a,\cal {\widehat{Y}}^a(\bf x)\Big)$ denote the feasible set defined by \eqref{eq_MMP3_yrecourse_o} and \eqref{eq_MMP3_yrecourse_c}, i.e., 
\begin{align*}\Big(\cal P^a, \cal {\widehat{Y}}^a(\bf x)\Big)=&\Big\{0\geq \max\big\{\sum_{\xi\in \bsm{\widehat \xi}} \eta^f_{\xi}p^f_{\xi}: (p^f_{\xi})_{\xi\in \bsm{\widehat \xi}}\in \cal P\big\} \\
	& \big\{\eta^f_{\xi}= \lVert\bf{\tilde y}^f_{\xi}\rVert_{1},\  (\bf y^f_{\xi},\tilde{\bf y}^f_{\xi})\in \tilde{\cal Y}(\bf x,\xi)\big\} \ \ \forall \xi\in \bsm{\widehat\xi}\Big\}.
\end{align*} 
\noindent Claim: There exists an optimal solution to  $\bf{MMP}^a_1$ such that its component, 
$(p^{o*}_{\xi},\bf y^{o*}_{\xi},\tilde{\bf y}^{o*}_{\xi})_{\xi\in \widehat{\bsm \xi}}$, belongs to $\Big(\cal P^a,\cal {\widehat{Y}}^a(\bf x)\Big)$. \\
\textit{Proof of Claim:} Suppose it is not true. So, there exists at least one scenario $\xi'\in\widehat{\bsm\xi}$ such that $\eta_{\xi'}^{f*} p_{\xi'}^{f*} > 0$. Again, we let $\big(p_{\xi'}^{o*}, \bf{ y}^{o*}_{\xi'}, \bf{\tilde y}^{o*}_{\xi'}\big) = \big(p_{\xi'}^{f*}, \bf{y}^{f*}_{\xi'}, \bf{\tilde y}^{f*}_{\xi'}\big)$. Then, for \eqref{eq_MMP3_WCEV_c} we have
\[\eta_{\xi'}^{o*}p_{\xi'}^{o*} = p_{\xi'}^{o*} \Big(f_2(\bf x, \xi',\bf y^{o*}_{\xi}) + \bf{M} \bf{\tilde y}^{o*}_{\xi'}\Big) > p_{\xi'}^{o*} f_2(\bf x, \xi',\bf y^{o*}_{\xi}),\]
which is contradictory to the updated \eqref{eq_MMP3_WCEV_c}, i.e., \eqref{eq_updated_no_MMP_bigM}. \qed

Hence, $(\bf y^o_{\xi},\tilde{\bf y}^o_{\xi})$ and $(\bf y^f_{\xi},\tilde{\bf y}^f_{\xi})$ can be unified into $(\bf y_{\xi},\tilde{\bf y}_{\xi})$ without affecting the optimality of $\bf{MMP}_1^a$. By doing so, $\bf{MMP}_1^a$ is equivalent to $\bf{MMP}_2$.
\end{proof}

\begin{remark}
\label{rmk:MMP_CCG_DRO}
$(i)$ Note that  ``$=$'' in \eqref{eq_MMP_Oscenarios} or \eqref{eq_MMP_Fscenarios} can be changed to ``$\leq$'' without affecting the optimality of $\bf{MMP}$, which  may lead to some computational improvement when $|\widehat{\bsm \xi}^o|$ or $|\widehat{\bsm\xi}^f|$ is large. More importantly, both $\bf{MMP}_1$ and $\bf{MMP}_2$ are  bilevel optimization formulations with two lower-level problems for feasibility and optimality respectively. 
For $\cal P$ defined in \eqref{eq_ambiguity}, the lower-level problems in \eqref{eq_WCEV_MMP1} and \eqref{eq_WCEV_MMP2} (or in \eqref{eq_WCEV_MMP12} and \eqref{eq_WCEV_MMP22}) are LPs and can be directly replaced by their optimality conditions or dual problems as in Appendix A.2. Actually, by making use of their specific structures and their dual problems, simpler single-level equivalent reformulations can be obtained. Consider $\bf{MMP}_2$ for demonstration, with $(\alpha^f,\bsm\beta^f)$ and $(\alpha^o,\bsm\beta^o)$ denoting  dual variables of  constraints of $\cal P$ in \eqref{eq_WCEV_MMP12} and \eqref{eq_WCEV_MMP22}. 
\begin{align}
	\begin{split}
		\label{eq_MMP_duality}
		\ubar{w} = \min\Big\{f_1(\bf x) + \eta\!: \ &\bf x\in \cal{X}, \ \eta\geq \alpha^o+\sum_{i=1}^m\gamma_i\beta^o_i, \ \  \alpha^o+\sum_{i=1}^m\psi_i(\xi)\beta^o_i\geq \eta^o_{\xi} \  \ \forall \xi\in \bsm{\widehat \xi},\\
		& 0\geq  \alpha^f+\sum_{i=1}^m\gamma_i\beta^f_i, \ \  \alpha^f+\sum_{i=1}^m\psi_i(\xi)\beta^f_i\geq \eta^f_{\xi}\ \  \forall \xi\in \bsm{\widehat \xi}, \\ 
		& \eqref{eq_MMP_Oscenarios2}, \ \eqref{eq_MMP_Fscenarios2},  \ \eqref{eq_MMP_yrecourse_relax}, \  \beta^o_i\geq 0, \ \ \forall i\in [m],\ \beta^f_i\geq 0 \ \ \forall i\in [m]
		\Big\}.
	\end{split}
\end{align}
We mention that \eqref{eq_MMP_duality} is a big-M-free single-level formulation. If $f_1$ ad $f_2$ are linear and $\cal X$ and $\cal Y(\cdot,\cdot)$ are linearly representable, it is a mixed integer linear program that is rather tractable by state-of-the-art professional solvers. \\
$(ii)$ For $\bf{MMP}_2$, it is not necessary to define the optimality cutting plane with respect to the whole set  $\bsm{\widehat \xi}$. It is worth noting from our numerical studies that $|\widehat{\bsm \xi^f}|$ could be more than an order of magnitude greater than $|\widehat{\bsm \xi}^o|$. Since the feasibility issue has not be analyzed in the existing literature, such a huge difference is very new and  unexpected. Hence, rather than employing $\bsm{\widehat \xi}$ in both optimality and feasibility cutting planes for consistency, it would be computationally more effective to build the first one with respect to  $\widehat{\bsm \xi}^o$ only. When computing large-scale instances, the benefit of this modification is often obvious.   \\
 $(iii)$ If the recourse problem is assumed to be feasible all the time, \eqref{eq_MMP_duality} reduces to the duality based master problem appeared in the literature  (e.g., \citep{zhao2018data,saif2021data,long2024supermodularity}). Yet, its derivation from $\bf{MMP}$ is intuitive, involves the knowledge of LP only, and more importantly, allows us to consider sophisticated $\cal P$ using mature optimization techniques. Unless specified otherwise,  \eqref{eq_MMP_duality} is utilized as $\bf{MMP}$ in our algorithm implementation. \qed
\end{remark}

\noindent (II): \textbf{Subproblems of C\&CG-DRO}\\
The two subproblems are referred to as main subproblems for consistency. As mentioned after Lemma \ref{lem_feasibility_WCEV}, the feasibility one is $\bf{WCEV}(\it F)$ in \eqref{eq_feasibility_subprob}, and the optimality one $\bf{WCEV}(\it O)$ in \eqref{eq_WCEV} (with $Q(\xi)$ replaced by $Q(\bf x, \xi)$). Because they can be solved by the oracles developed in Section \ref{sect_EVWCEV} in an almost identical  fashion, we next present the customization of those oracles on $\bf{WCEV}(\it O)$ only.  Specifically, for given $\bf x^*$, the finite mathematical program for  $\bf{WCEV}(\it O)$ can be easily obtained by modifying \eqref{eq_WCEV_opt} as in the following.
\begin{align}
	\label{eq_WCEV-OFMP}
\begin{split}
	\bf{WCEV(\it O)-FMP}: \  \eta^o(\bf x^*)=\max\Bigg\{\sum_{j=1}^{m+1} \eta_jp_j: 
	& (p_1,\dots, p_{m+1})\in \cal{P}, \\ 
	\  \xi_j\in \Xi \ \ \forall j\in [m+1], \  \Big\{\eta_j= \min \big\{f_2 (\bf x^*,\xi_j,\bf y_j), & \bf y_j\in \cal Y(\bf x^*,\xi_j)\big\}\Big\} \ \ \forall j\in[m+1]\Bigg\}.
\end{split}
\end{align}
Comparing \eqref{eq_WCEV_opt} and \eqref{eq_WCEV-OFMP}, the basic difference is that variable $\eta_j$ replaces $Q(\xi_j)$, and $\eta^j$ is set to the recourse cost for $\xi_j$ through a replica of the recourse problem parameterized by $\xi_j$ for all $j$. Note that the equality sign associated with $\eta_j$  can be changed to $\leq$ without affecting the whole formulation's optimality. 

If  $\bf{WCEV}(\it O)$ is solved by the CG-based \textit{Oracle-2}, the customized $\bf{PMP}$, denoted by $\bf{PMP}(\it O)$, is a simple LP with $Q(\xi^0_j)$ replaced by $Q(\bf x^*, \xi^0_j)$ in \eqref{eq_PMP_basic}. The customized $\bf{PSP}$, denoted by $\bf{PSP(\it O)}$, is in the following format. 
\begin{align}
	\label{eq_PSPQ}
	\bf{PSP(\it O)}:  \ v^*(\bf x^*,\bsm\xi^0_n) = \max_{\xi \in \Xi} \min_{\bf y\in \cal Y({\bf x^*,\xi})} f_2(\bf x^*,\xi,\bf y) - \alpha^* -\sum_{i=1}^m\psi_i(\xi)\beta^*_i.
\end{align}  
	Both \eqref{eq_WCEV-OFMP} and \eqref{eq_PSPQ} are bilevel optimization formulations that can be solved by the methods presented in Appendix A.2.	To unify our exposition, regardless of using \textit{Oracle-1} or \textit{Oracle-2} to solve $\bf{WCEV}(\it O)$ and $\bf{WCEV}(\it F)$, we let $\eta^o(\bf x^*)$ denote the optimal value of $\bf{WCEV}(\it O)$, $\hat{\bsm \xi}^o(\bf x^*)\equiv\{\xi^{o}_1,\dots, \xi^{o}_{|\hat{\bsm \xi}^o(\bf x^*)|}\}$ and $P^o\big(\hat{\bsm \xi}^o(\bf x^*)\big)$ the set of resulting scenarios and their probabilities; and $\eta^f(\bf x^*)$,  $\hat{\bsm \xi}^f(\bf x^*)$ and $P^f\big(\hat{\bsm \xi}^f(\bf x^*)\big)$ to their counterparts for $\bf{WCEV}(\it F)$.
\begin{remark}
	\label{Remark_9}
 As noted previously, when the recourse problem is an LP and there is no feasibility issue $\bf x\in \cal X$, Benders type of algorithms have been developed to solve $\mathbf{2-Stg \ DRO}$ \cite{bansal2018decomposition,duque2022distributionally,gamboa2021decomposition,gangammanavar2022stochastic,luo2022decomposition}. With the primal oracles for $\mathbf{WCEV}(\it O)$ and resulting $\hat {\bsm\xi}^o$, they naturally can be extended for possible enhancements. Specifically, let
the recourse problem be 
\begin{align}
\label{eq_recourse_LP} 
	\min \big\{ \mathbf{c}_2\mathbf y: \ \bf B\bf y \geq \bf b_2 - \bf A_2\bf x - \bf H \xi\big\}.
\end{align}
Note that for $\xi\in \hat{\bsm \xi}^o(\bf x^*)$ we can derive an optimal solution for the dual problem of \eqref{eq_recourse_LP}. Denoting them by 
$\{\pi^*_\xi\}_{\xi\in \hat{\bsm \xi}^o(\bf x^*)}$, Benders cutting planes, generated after solving the current $\mathbf{WCEV}(\it O)$, are presented in the following. 
\begin{align*}
	& \alpha^o+\sum_{i}\psi_i(\xi)\beta^o_i\geq \eta^{o}_{\xi} \ \ \forall \xi\in \hat{\bsm\xi}^o(\bf x^*)\\
	&\eta^{o}_{\xi} \geq (\bf b_2 - \bf A_2\bf x - \bf H\xi)^\intercal \pi^*_{\xi} \ \ \forall \xi\in \hat{\bsm \xi}^o(\bf x^*)\tag*{\qed}
\end{align*}

\end{remark}
Indeed, although it has not been investigated in the current literature,  Benders cutting planes to address the feasibility issue can be produced in the same fashion after solving $\mathbf{WCEV}(\it F)$. To differentiate from existing implementations and to be consistent, we refer to our new approach as Benders-DRO and the previous ones as basic Benders. Actually, as shown in the numerical study, Benders-DRO performs drastically better than basic Benders.

\subsection{Complete C\&CG Procedure for $\bf{2-Stg \ DRO}$}
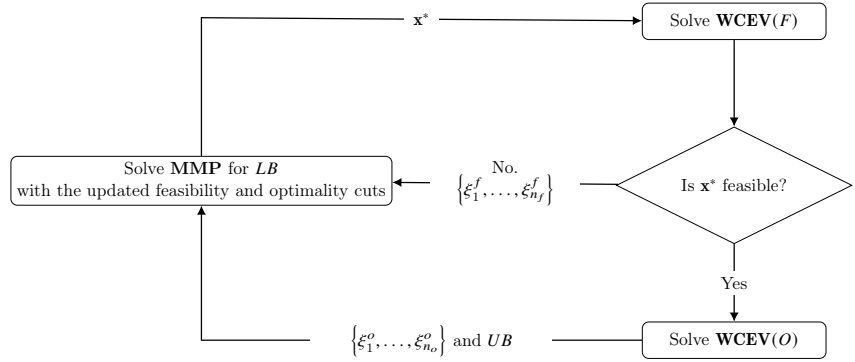
\begin{figure}[H]
	\centering
	\bigskip
	\tikzset{%
		>={Latex[width=2mm,length=2mm]},
		CCGCG/.style = {rectangle, rounded corners, draw=black,
			minimum width=4cm, minimum height=0.8cm,
			text centered},
		CCG/.style = {rectangle, rounded corners, draw=black,
			minimum width=4cm, minimum height=0.8cm,
			text centered},
		if/.style = {diamond, draw=black, aspect=2, text width=10em},
	}
	\resizebox{0.7 \textwidth}{!}{
		\begin{tikzpicture}[node distance=1.6cm,
			every node/.style={fill=white}, align=center]
			\node (MMP)             [CCG]              {Solve \textbf{MMP} for $LB$\\ with the updated feasibility and optimality cuts};
			\node (PP-F)     [CCGCG, right of = MMP, xshift = 10cm, yshift = 3.5cm]          {Solve $\bf{WCEV}(F)$};
			\node (feasible)     [if, below of=PP-F, yshift=-2cm]          {Is $\bf{x}^*$ feasible?};
			\node (PP-O)     [CCG, right of = MMP, xshift = 10cm, yshift = -3.5cm]          {Solve $\bf{WCEV}(O)$};
			\draw[->, thick]             (MMP) --++ (0,3.58) --node[text width=0.5cm]{$\bf{x}^*$}(PP-F);
			\draw[->, thick]             (PP-F) --(feasible);
			\draw[->, thick]             (feasible) -- node[text width=3cm]
			{No. \\
				$\Big\{\xi^f_1,\dots, \xi^f_{n_f}\Big\}$}(MMP);
			\draw[->, thick]             (feasible) -- node[text width=1.5cm]
			{Yes}(PP-O);
			\draw[->][thick, draw = black]             (PP-O)--++ (-11.6,0)-- node[xshift=5cm, yshift=-1.5cm,text width=5cm] 
			{$\Big\{\xi^o_1,\dots, \xi^o_{n_o}\Big\}$ and $UB$}(MMP);
	\end{tikzpicture}}
	\caption{The Schematic Flow of C\&CG-DRO\label{fig:CCG-DRO_sketch}}
\end{figure}

With the aforementioned problems defined, we are ready to  present the overall procedure of C\&CG method customized to compute $\bf{2-Stg \ DRO}$. Note that $LB$ and $UB$ denote lower and upper bounds respectively, $TOL$ is the optimality tolerance, and $t$ is the counter for iterations. Also, to facilitate our understanding, we sketch the logic and main steps in  Figure \ref{fig:CCG-DRO_sketch}.  \\

\refstepcounter{alg} 
\noindent\begin{minipage}{\textwidth}
    \noindent\textbf{Algorithm 1:} C\&CG-DRO\label{alg:2}\\[-0.8em]
    \makebox[0.92\textwidth]{\hrulefill}\\[-0.8em]
\end{minipage}
\begin{algorithmic}
		\State \textbf{Step 1}  Set $LB = -\infty$, $UB = +\infty$, $t=1$, and $\bsm{\widehat\xi}^o = \bsm{\widehat\xi}^f =\bsm{\widehat\xi}=\emptyset$.
		\State \textbf{Step 2}  \parbox[t]{\dimexpr\linewidth-\algorithmicindent}{Solve master problem \textbf{MMP}. If it is infeasible, report the infeasibility of $\bf{2-Stg \ DRO}$ and terminate. If it is unbounded, select an arbitrary new feasible solution $\bf x^*$. Otherwise, derive  optimal  value $\ubar w$ and solution $\bf x^*$, and update $LB = \ubar w$.}\\
		\State \textbf{Step 3}  \parbox[t]{\dimexpr\linewidth-\algorithmicindent}{Solve feasibility subproblem $\bf{WCEV(\it F)}$, derive its optimal value $\eta^f(\bf x^*)$, optimal set of scenarios $\hat{\bsm \xi}^f(\bf x^*)$ and associated probabilities (by one of oracles from Section \ref{sect_primal_DRO}).}\\
		\State \textbf{Step 4} Cases based on $\eta^f(\bf x^*)$
		 \begin{description}
	      \item[$\ \ \ $\textbf{Case A \ $\eta^f(\bf x^*) = 0$}]\\
		   $\ \ \ \ \ \ \ \ \ \ $$(i)$ \parbox[t]{\dimexpr\linewidth-\algorithmicindent}{Solve optimality subproblem $\bf{WCEV}(\it O)$, derive its optimal value $\eta^o(\bf x^*)$, optimal set of scenarios $\hat{\bsm \xi}^o(\bf x^*)$ and their probabilities (by one of oracles from Section \ref{sect_primal_DRO});}\\
		    $\ \ \ \ \ \  \ \ \ \ $$(ii)$ \parbox[t]{\dimexpr\linewidth-\algorithmicindent}{Update $\widehat{\bsm \xi}^o=\widehat{\bsm \xi}^o\cup \hat{\bsm \xi}^o(\bf x^*)$ (and accordingly $\widehat{\bsm \xi}$), and augment $\bf{MMP}$ with new variables and constraints accordingly. }\\
		  \item[$\ \ \ $\textbf{Case B \ $\eta^f(\bf x^*) > 0$}]\\
		  $\ \ \ \ \ \  \ \ \ \ $$(i)$ \parbox[t]{\dimexpr\linewidth-\algorithmicindent}{Update $\widehat{\bsm\xi}^f=\widehat{\bsm\xi}^f\cup\hat{\bsm\xi}^f(\bf x^*)$ (and accordingly $\widehat{\bsm \xi}$), and augment $\bf{MMP}$ with new variables and constraints accordingly;}\\
		  $\ \ \ \ \ \  \ \ \ \ $ $(ii)$ \parbox[t]{\dimexpr\linewidth-\algorithmicindent}{set $\eta^o(\bf x^*) = +\infty$. }
		  \end{description}		
		\State \textbf{Step 5} Update $UB = \min\{UB, f_1(\bf x^*)+\eta^o(\bf x^*)\}$. 
		\State \textbf{Step 6} \parbox[t]{\dimexpr\linewidth-\algorithmicindent}{If $UB - LB \leq TOL$, return $\bf x^*$ and terminate. Otherwise, set $t=t+1$ and go to Step 2.}
	\end{algorithmic}
\vspace{-1.4em}
\noindent\makebox[0.92\textwidth]{\hrulefill}\\

We note that Algorithm 1 is a rather vanilla version of  C\&CG-DRO. Some simple changes, as described in the following, often can yield significantly computational improvements, especially when subproblems are computed by CG-based \textit{Oracle-2}.   

\begin{description}
\item[\textbf{(C1)}] Rather than letting \textit{Oracle-2} start from scratch, $\bsm{\widehat \xi}^f$ and $\bsm{\widehat \xi}^o$ obtained up to date can be used as the initial scenario sets for it to solve $\bf{WCEV}(\it F)$ and $\bf{WCEV}(\it O)$, respectively. If this is the case, the following operation should be included between Step 2 and Step 3. \\
		 \textbf{ -- Step 2.a}  For $\xi^f\in \widehat{\bsm\xi}^f$, compute $\tilde Q(\bf x^*,\xi^f)$; and for $\xi^o\in \widehat{\bsm\xi}^o$, compute $Q(\bf x^*, \xi^o)$.\\
		 With this change, we can simplify \textbf{Step 4} as the following. \\ 
		\textbf{--} In \textbf{Step 4-Case A} and \textbf{Step 4-Case B}, update $\widehat{\bsm \xi}^f=\hat{\bsm \xi}^f(\bf x^*)$ and  $\widehat{\bsm \xi}^o=\hat{\bsm \xi}^o(\bf x^*)$, respectively. \\
 	We note that those changes are included in our default implementation of C\&CG-DRO with \textit{Oracle-2}.  
 
 \item[\textbf{(C2)}] Solving subproblems often leads to many scenarios with zero probability in $\hat{\bsm\xi}^f(\bf x^*)$ and $\hat{\bsm\xi}^o(\bf x^*)$, which could be more prominent when \textit{Oracle-2} is applied. To keep $\bf{MMP}$ more tractable, we only set $\widehat{\bsm \xi}^f$ and $\widehat{\bsm \xi}^o$ to include scenarios of non-zero probabilities from $\hat{\bsm \xi}^f(\bf x^*)$ and $\hat{\bsm \xi}^o(\bf x^*)$ in our default implementation of C\&CG-DRO. Yet, to keep the following theoretical analyses simple, C\&CG-DRO does not screen out scenarios of zero probability.

\item[\textbf{(C3)}] If we are concerned with the feasibility of $\bf{PMP}(\it O)$ (i.e., the pricing master problem of \textit{Oracle-2} for $\bf{WCEV}(\it O)$), instead of using Farka's Lemma modified implementation, we can simply use $\hat{\bsm\xi}^f(\bf x^*)$ to initialize   $\bf{PMP}(\it O)$. Since  $\bf{WCEV}(\it O)$ is to be solved in $\bf{Step \ 4}-\bf{Case \ A}$,  employing $\hat{\bsm\xi}^f(\bf x^*)$ guarantees that $\bf{PMP}(\it O)$ is feasible. Yet, as noted earlier,  set $\hat{\bsm\xi}^f(\bf x^*)$ might not be small and could slow down the computation of  $\bf{PMP}(\it O)$.
\end{description}

\begin{figure}[H]
	\centering
	\tikzset{%
		>={Latex[width=2mm,length=2mm]},
		beginend/.style = {rectangle, rounded corners, draw=black,
			minimum width=4cm, minimum height=1cm,
			text centered},
		process/.style = {rectangle, draw=black,
			minimum width=2cm, minimum height=0.5cm,
			text centered},
		if/.style = {diamond, draw=black, aspect=3, text width=10em},
		basic box/.style={shape=rectangle, align=center, draw=red, rounded corners},
		largebasic box/.style={shape=rectangle, align=center, draw=blue, fill = none, rounded corners},
	}
	\resizebox{1 \textwidth}{!}{
		\begin{tikzpicture}[node distance=1.6cm,
			every node/.style={fill=white}, align=center]
			\node (Setup)             [beginend]              {$LB = -\infty$, $UB = +\infty$, $t=1$, \\ and $\widehat{\bsm \xi} = \widehat{\bsm \xi}^o = \widehat{\bsm \xi}^f = \emptyset$};
			\node (SolveMMP)     [process, below of=Setup]          {Solve \textbf{MMP}};
			\node (Feasible)     [if, below of=SolveMMP, yshift=-0.5cm]          {Is \textbf{MMP} feasible?};
			\node (bounded)     [if, below of=Feasible, yshift=-1.5cm]          {Is \textbf{MMP} bounded?};
			\node (DROInfeasible)     [beginend, right of=Feasible, xshift=5cm]          {\textbf{2-Stg DRO} is infeasible};
			\node (Update)     [process, below of=bounded, yshift=-2cm]          {Compute $\{Q(\mathbf x^*, \xi^o)\}_{\xi^o \in \widehat{\bsm \xi}^o}$ and $\{\tilde Q_f(\mathbf x^*, \xi^f)\}_{\xi^f \in \widehat{\bsm \xi}^f}$};
			\node (SolvePMPF)     [process, below of=Update, xshift=-5cm, yshift=-0.5cm]          {Solve $\mathbf{PMP}(F)$};
			\node (SolvePSPF)     [process, below of=SolvePMPF, yshift=-1cm]          {Solve $\mathbf{PSP}(F)$};
			\node (vf)     [if, below of=SolvePSPF, yshift=-0.5cm]          {Is $v^{f*}(\mathbf x^*, \widehat{\bsm \xi}^f)\leq 0$?};
			\node (hatetaf)     [if, below of=vf, yshift=-1.5cm]          {Is $\ubar{\eta}^{f*}(\mathbf x^*) = 0$?};
			\node (SolvePMPO)     [process, below of=Update, xshift=4cm, yshift=-0.5cm]          {Solve $\mathbf{PMP}(O)$};
			\node (SolvePSPO)     [process, below of=SolvePMPO, yshift=-1cm]          {Solve $\mathbf{PSP}(O)$};
			\node (vo)     [if, below of=SolvePSPO, yshift=-0.5cm]          {Is $v^{o*}(\mathbf x^*, \widehat{\bsm \xi}^o)\leq 0$?};
			\node (TOL)     [if, below of=Update, yshift=-10cm]          {Is $UB-LB \leq TOL$?};
			\node (termination)     [beginend, below of=TOL, yshift=-1.5cm]          {Terminate};
			\draw[->]             (Setup) -- (SolveMMP);
			\draw[->]             (SolveMMP) -- (Feasible);
			\draw[->]             (Feasible) -- node[text width=4cm]
			{Yes}(bounded);
			\draw[->]             (Feasible) -- node[text width=1em]
			{No}(DROInfeasible);
			\draw[->]             (bounded)  -- ++(5,0) -- ++(0,-3.5) -- node[xshift=1.5cm, yshift=1.5cm, text width=4cm]
			{No. Find an arbitrary new feasible solution $\mathbf x^*$}(Update);
			\draw[->]             (bounded) -- node[text width=4cm]
			{Yes. Get $\mathbf x^*$ and set $LB=\ubar{w}$}(Update);
			\draw[->]             (Update)  -- ++(0,-1) -- ++(-5,0) -- ++(0,-0.5) --  (SolvePMPF);
			\draw[->]             (SolvePMPF) -- node[text width=3cm]
			{$\alpha^{f*}, \{\beta^{f*}_i\}$}(SolvePSPF);
			\draw[->]             (SolvePSPF) -- (vf);
			\draw[->]             (vf) -- node[text width=2cm]
			{Yes}(hatetaf);
			\draw[->]             (vf)  -- ++(-3.5,0) -- ++(0,4.7) -- node(renew-F)[xshift=-1.5cm, yshift=-2cm, text width=4cm]
			{No. Update $\widehat{\bsm \xi}^f$}(SolvePMPF);
			\draw[->]             (hatetaf) -- ++(4.5,0) -- ++(0,4.7) -- ++(0,3) -- node[xshift=-1.6cm, yshift=-3.5cm, text width=2cm]
			{Yes.}(SolvePMPO);
			\draw[->]             (hatetaf) -- ++(-6.5,0) -- ++(0,18.6) -- node(LeftNode)[xshift=-4.5cm, yshift=-8cm, text width=4cm]
			{No. Update $\widehat{\bsm \xi}$ and $UB$}(SolveMMP);
			\draw[->]             (SolvePMPO) -- node[text width=3cm]
			{$\alpha^{o*}, \{\beta^{o*}_i\}$}(SolvePSPO);
			\draw[->]             (SolvePSPO) -- (vo);
			\draw[->]             (vo)  -- ++(0,-2.4) -- ++(-4,0) -- ++(0,-0.5) --  node[xshift=2.1cm, yshift=0.9cm, text width=3.2cm]
			{Yes. Update $UB$}(TOL);
			\draw[->]             (vo)  -- ++(3.5,0) -- ++(0,4.7) -- node(renew-O)[xshift=1.5cm, yshift=-2cm, text width=4cm]
			{No. Update $\widehat{\bsm \xi}^o$}(SolvePMPO);
			\draw[->]             (TOL) -- node[text width=2cm]
			{Yes. Return $\mathbf x^*$}(termination);
			\draw[->]             (TOL) -- ++(10.5,0) -- ++(0,20.4) -- node(RightNode)[xshift= 4.5cm, yshift=-8cm, text width=4cm]
			{No. Update $\widehat{\bsm \xi}$}(SolveMMP);
			\begin{scope}[on background layer]
				\node[fit = (renew-F)(SolvePMPF)(SolvePSPF)(vf),basic box=white,
				inner ysep=1.5em,
				label={[xshift=4.5cm, yshift=-3cm, anchor=center, fill=white, draw, text width=2cm]{\textit{Oracle-2 }}}] {};
			\end{scope}
			\begin{scope}[on background layer]
				\node[fit = (renew-O)(SolvePMPO)(SolvePSPO)(vo),basic box=white,
				inner ysep=1.5em,
				label={[xshift=-4.5cm, yshift=-3cm, anchor=center, fill=white, draw, text width=2cm]{\textit{Oracle-2 }}}] {};
			\end{scope}
	\end{tikzpicture}}
	\caption{Flowchart for C\&CG-DRO algorithm with \textit{Oracle-2}}
	\label{CCG-DRO-flow-chart}
\end{figure}
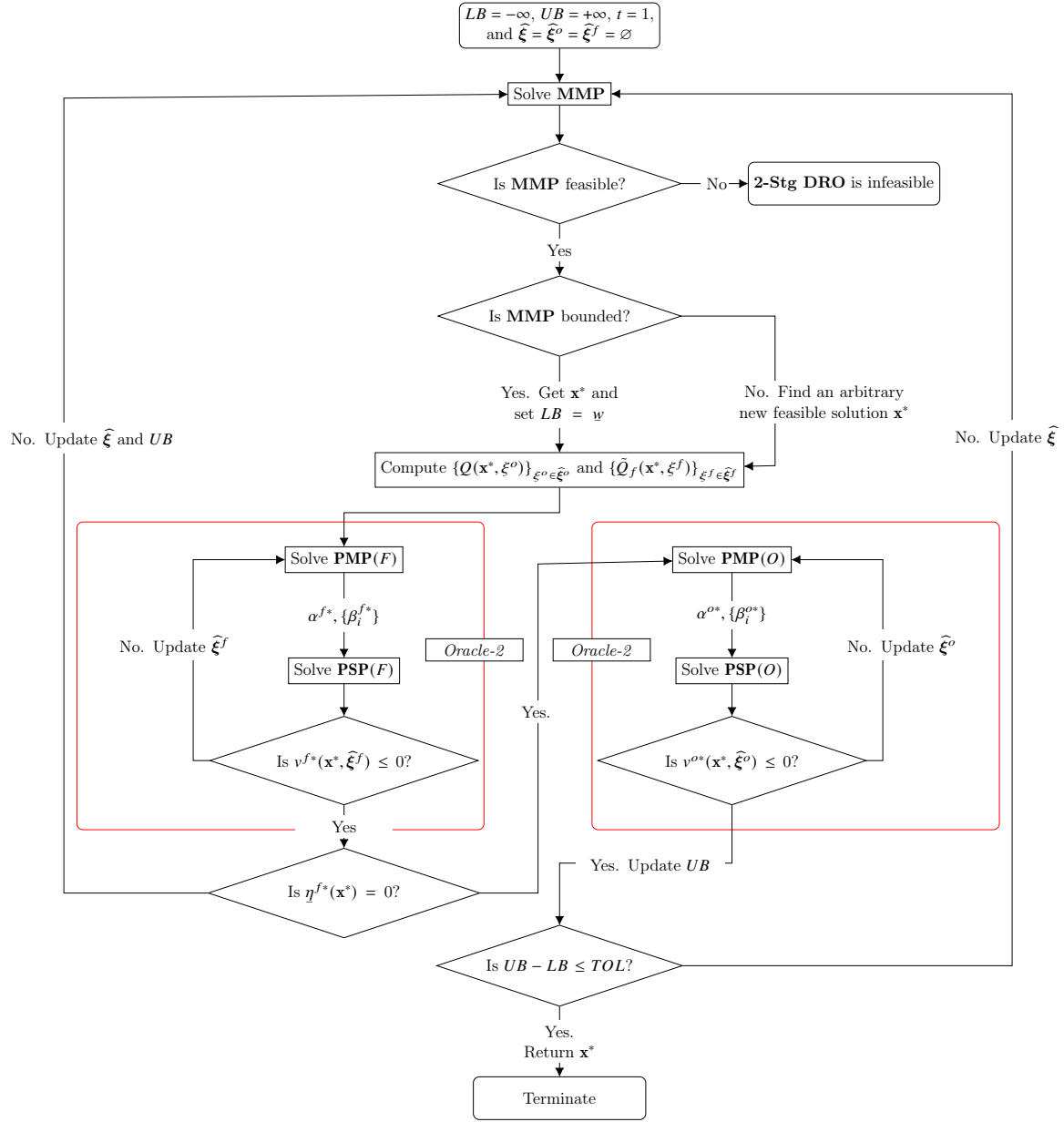

A detailed flow chart describing C\&CG-DRO with \textit{Oracle-2} is presented in Figure~\ref{CCG-DRO-flow-chart} at the end of this paper.
In the next subsection, we analyze C\&CG-DRO's  strength with respect to existing approaches, and its convergence and iteration complexity.

\subsection{Strength, Convergence and Iteration Complexity}
\label{subsect_conv_comp}
As previously mentioned,  the classical C\&CG implementation, referred to as basic C\&CG, is one of the most popular methods  for solving $\bf{2-Stg \ DRO}$. As a new and more intricate algorithm, it is significant to analytically demonstrate C\&CG-DRO's strength compared to that of basic C\&CG. Next, we show that this actually is the case under some mild conditions. Since the infeasibility issue of the recourse problem has not been studied in the literature, we assume  that there is no feasibility issue associated with $\bf x \in \cal X$. We also note that   $\textit{Oracle-2}$ is adopted in C\&CG-DRO  in the following analysis, because of its available infrastructure.

\begin{prop}  
	$(i)$ Assume that both  basic C\&CG and C\&CG-DRO have the same scenario set $\widehat{\bm \xi}^o$ before the start of iteration $t^0$. We further assume that the solution algorithms for their master and subproblems  are deterministic and $\bf{PMP}(\it O)$ has a unique shadow price with respect to $\widehat{\bm \xi}^o$. Then the lower bound derived from  basic C\&CG in iteration $t^0+1$ is dominated by that of C\&CG-DRO. \\
	$(ii)$ Assume that the (main)  master problems of basic C\&CG and C\&CG-DRO generate the same first-stage solution in iteration $t^0$. If the upper bound is updated in basic C\&CG in iteration $t^0$, this upper bound is dominated by that of C\&CG-DRO.    
\label{proposition 2}
\end{prop}
\begin{proof}
	$(i)$ As noted in Remark \ref{rmk:MMP_CCG_DRO}, \eqref{eq_MMP_duality}, the duality based reformulation of $\bf{MMP}_2$, serves as the (main) master problem for both basic C\&CG and C\&CG-DRO. We let $(\bf x^*, \alpha^*,\bm\beta^*, \dots)$ denote the optimal solution obtained from computing the (main) master problem for both algorithms, noting that $\widehat{\bsm \xi}=\widehat{\bsm \xi}^o$. With \textbf{(C1)} implemented and the stated assumptions, it can be seen that $(\alpha^*,\bm\beta^*)$ is the shadow price of $\bf{PMP}(\it O)$ for $\widehat{\bm \xi}^o$ and $\bf x^*$. Hence, the associated $\bf{PSP}(\it O)$ is  
	\begin{align*}
		\max_{\xi \in \Xi} Q(\bf x^*,\xi) - \alpha^* -\sum_{i=1}^m\psi_i(\xi)\beta^*_i.
	\end{align*}
	Also, the subproblem of basic C\&CG is
    \begin{align*}
		\max_{\xi \in \Xi} Q(\bf x^*,\xi) -\sum_{i=1}^m\psi_i(\xi)\beta^*_i.
	\end{align*}
   Clearly, they are essentially same and have the identical optimal solution. Denoting it by $\xi'$, we therefore have $\xi'\in \hat{\bsm \xi}^o(\bf x^*)\backslash\widehat{\bm \xi}^o$, the set of new scenarios identified after executing \textit{Oracle-2}. As a result, after computing the updated (main) master problems for C\&CG and C\&CG-DRO respectively in iteration $(t^0+1)$, the lower bound derived by the former method is weaker than that generated by the latter one.

   \noindent $(ii)$    Again, in iteration $t^0$, let $\bf x^*$ be the optimal first-stage solution output from computing the (main) master problems,  $\xi'$ the optimal solution from solving the subproblem of basic C\&CG, and its optimal value $\alpha' =  Q(\bf x^*, \xi')-\sum_{i=1}^m\psi_i(\xi')\beta_i^*$. Hence, we have $$\alpha' \geq Q(\bf x^*, \xi)-\sum_{i=1}^m\psi_i(\xi)\beta_i^* \ \ \forall \xi\in \Xi.$$
Upon termination of \textit{Oracle-2} with $\epsilon=0$,  the dual problem of $\bf{PMP}(O)$ is
\begin{align*}
	\ubar{\eta}^{o*}\big(\bf x^*, \hat{\bm \xi}^o(\bf x^*)\big) = 
	\min\Big\{\alpha + \sum_{i=1}^m \beta_i \gamma_i: \alpha + \sum_{i=1}^m \psi_i(\xi)\beta_i \geq Q(\bf x^*, \xi) \ \ \forall \xi\in \hat{\bm \xi}^o(\bf x^*), \ \beta_i\geq 0\ \ \forall i\in [m] \Big\}.
\end{align*}

Given that $\hat{\bm \xi}^o(\bf x^*)\subseteq \Xi$, it is clear that $(\alpha', \bm\beta^*)$ is a feasible solution to this dual problem. So, we have $\ubar{\eta}^{o*}\big(\bf x^*, \hat{\bm \xi}^o(\bf x^*)\big)\leq \alpha'+\sum_{i=1}^m\beta^*_i\gamma_i$, leading to 
$$f_1(\bf x^*)+\ubar{\eta}^{o*}\big(\bf x^*, \hat{\bm \xi}^o(\bf x^*)\big)\leq f_1(\bf x^*)+\alpha'+\sum_{i=1}^m\beta^*_i\gamma_i.$$
Note that the left-hand-side and right-hand-side of this inequality are benchmarked against the current upper bounds  and employed to update them if applicable by C\&CG-DRO and basic C\&CG respectively (e.g., Step 5 in C\&CG-DRO). Hence, if an update occurs to the upper bound in basic C\&CG, the new upper bound is dominated by the upper bound provided by C\&CG-DRO. 
\end{proof}

Next, we show that C\&CG-DRO eventually solves $\bf{2-Stg \ DRO}$.  We assume that \textit{Oracle-1} is adopted as it does not incur any numerical issue for $\bf{WCEV}$ problem.

\begin{theorem}
	\label{CCG-DRO-convergence}
	For a fixed $\bf x\in \cal X$, let $\Xi^o(\bf x) = \big\{\xi\in \Xi:\ Q(\bf x, \xi) <+\infty\big\}$ denote the feasible sample space under $\bf x$, and the assumption made in Theorem \ref{CG-epsilon} holds for $\bf x\in \cal X$, i.e.,  $$Q(\bf x, \xi)- \sum_{i}\psi_i(\xi)\beta_i$$
	is continuous, uniformly with respect to $\xi$ for $\bf x\in \cal X$ and $\bsm \beta \geq \bf 0$.\\
	$(i)$ If there is no feasibility issue for the recourse problem, i.e., all $\bf x\in \cal X$ are feasible, C\&CG-DRO returns an optimal (or $\varepsilon-$optimal) one in a finite number of iterations. \\
	$(ii)$ When the feasibility issue exists, we further assume $\Xi$ is a polytope, $\tilde Q_f (\bf x, \xi)$ is convex over $\Xi$, $Q (\bf x, \xi)$ is uniformly continuous over $\Xi^o(\bf x)$ for $\bf x \in \cal X$, and $\psi_i(\xi)$ are concave over $\Xi$ for $i\in [m]$. C\&CG-DRO either reports that $\bf{2-Stg \ DRO}$ is infeasible or  returns an optimal (or $\varepsilon-$optimal) one in a finite number of iterations.  
\end{theorem}
\begin{proof}
	We first present the proof for the first statement when $\bf{2-Stg\ DRO}$ has no feasibility issues.  Note that according to the definition of the uniform continuity, we have that for any given  $\varepsilon > 0$ there exists a $\delta > 0$ satisfying
	\[\Big\vert\big(Q(\bf x, \xi_1)- \sum_{i}\psi_i(\xi_1)\beta_i\big) - \big(Q(\bf x, \xi_2)- \sum_{i}\psi_i(\xi_2)\beta_i\big)\Big\vert \leq \varepsilon,\]
	if $\Vert \xi_1 - \xi_2 \Vert \leq \delta$ for  $\xi_1,\xi_2 \in \Xi$, $\bf x\in \cal X$, and $\bsm \beta \geq \bf 0$.

	\noindent Claim 1: Suppose that $\widehat{\bsm \xi}^o$  is currently available scenarios, $\bf x^*$ is an optimal solution to $\bf{MMP}$ defined with $\widehat{\bsm \xi}^o$, and $\hat{\bsm \xi}^o(\bf x^*)$ has been obtained by solving $\bf{WCEV(\it{O})}$ for $\bf x^*$ in iteration $t^0$. If the Hausdorff distance between sets  $\hat{\bsm \xi}^o(\bf x^*)$  and $\widehat{\bsm \xi}^o$, denoted by $h\big(\hat{\bsm \xi}^o(\bf x^*),\widehat{\bsm \xi}^o\big)$, is less than or equal to $\delta$, we have $UB \leq LB + \varepsilon$, i.e., $\bf x^*$ is an $\varepsilon-$optimal solution.\\
	\textit{Proof of Claim 1:} 
	Note that, for $\bf x^*$, its optimal value of $\bf{WCEV}(\it O)$ reduces to 
	\begin{align*}
		\eta^o(\bf x^*)=\max&\Big\{\sum_{\xi\in \hat{\bsm \xi}^o(\bf x^*)}p_{\xi} Q(\bf x^*, \xi): (p_{\xi})_{\xi\in \hat{\bsm \xi}^o(\bf x^*)}\in \cal{P}\Big\}\\
		= \min & \Big\{\alpha+\sum_{i=1}^m\gamma_i\beta_i\!: \ \alpha+\sum_{i}\psi_i(\xi)\beta_i\geq Q(\bf x^*, \xi) \  \ \forall \xi\in \hat{\bsm \xi}^o(\bf x^*) \Big\}. 
	\end{align*}
	
	For the dual form \eqref{eq_MMP_duality} of $\bf{MMP}$, let $\big(\bf x^*, \alpha^{o*},\ \bsm{\beta}^{o*}, (\eta^{o*}_{\xi})_{\xi\in \bsm{\widehat\xi}^o}, \dots \big)$ be  components associated with $\bf x^*$. Without loss of generality, we assume that $\eta^{o*}_{\xi}= Q(\bf x^*,\xi)$. We have 
	\begin{align*}
		LB = \ubar{w}^{t^0} &= f_1(\bf x^*) + \eta^*\\
		\eta^* = &\alpha^{o*}+\sum_{i=1}^m\gamma_i\beta^{o*}_i\\
		\alpha^{o*}+\sum_{i}\psi_i(\xi)\beta^{o*}_i&\geq  Q(\bf x^*, \xi) \  \ \forall \xi\in  \widehat{\bsm \xi}^o.
	\end{align*}

	By the definition of Hausdorff distance, we have
	$\max_{\xi_1\in \hat{\bsm \xi}^o(\bf x^*)
	}\Big\{\min_{\xi_2 \in \widehat{\bsm \xi}^o} \{\Vert \xi_1 - \xi_2\Vert \} \Big\} \leq \delta.$
	Following the uniform continuity assumption on $Q(\bf x, \xi)- \sum_{i}\psi_i(\xi)\beta_i$, we have
	\begin{align*}
		&\max_{\xi \in \hat{\bsm \xi}^o(\bf x^*)}\{Q(\bf x, \xi) - \sum_{i}\psi_i(\xi)\beta^{o*}_i\}  - \max_{\xi \in \widehat{\bsm \xi}^o}\{Q(\bf x, \xi) - \sum_{i}\psi_i(\xi)\beta^{o*}_i\}\\
		\leq &\max_{\xi \in \hat{\bsm \xi}^o(\bf x^*)}\Big\{\{Q(\bf x, \xi) - \sum_{i}\psi_i(\xi)\beta^{o*}_i\} - \max_{\xi' \in \widehat{\bsm \xi}^o: \Vert \xi' - \xi\Vert \leq \delta}\{Q(\bf x, \xi') - \sum_{i}\psi_i(\xi')\beta^{o*}_i\}\Big\}\\
		\leq & \max_{\xi \in \hat{\bsm \xi}^o(\bf x^*)}\varepsilon=\varepsilon.
	\end{align*}
	Then, it can be inferred that 
	\begin{align*}
		\alpha^{o*} &\geq \max_{\xi \in \widehat{\bsm \xi}^o}\{Q(\bf x, \xi) - \sum_{i}\psi_i(\xi)\beta^{o*}_i\}\\
		&=\max_{\xi \in \widehat{\bsm \xi}^o}\{Q(\bf x, \xi) - \sum_{i}\psi_i(\xi)\beta^{o*}_i\} - \max_{\xi \in \hat{\bsm \xi}^o(\bf x^*)} \{Q(\bf x, \xi) - \sum_{i}\psi_i(\xi)\beta^{o*}_i\} +\\
		& \max_{\xi \in \hat{\bsm \xi}^o(\bf x^*)}\{Q(\bf x, \xi) - \sum_{i}\psi_i(\xi)\beta^{o*}_i\}\\
		&\geq \max_{\xi \in \hat{\bsm \xi}^o(\bf x^*)} \{Q(\bf x, \xi) - \sum_{i}\psi_i(\xi)\beta^{o*}_i\} - \Big\vert \max_{\xi \in \widehat{\bsm \xi}^o}\{Q(\bf x, \xi) - \sum_{i}\psi_i(\xi)\beta^{o*}_i\}  -\\
		& \max_{\xi \in \hat{\bsm \xi}^o(\bf x^*)}\{Q(\bf x, \xi) - \sum_{i}\psi_i(\xi)\beta^{o*}_i\} \Big\vert\\
		&\geq  \max_{\xi \in \hat{\bsm \xi}^o(\bf x^*)} \{Q(\bf x, \xi) - \sum_{i}\psi_i(\xi)\beta^{o*}_i\} - \varepsilon.
	\end{align*}
	Therefore, $(\alpha^{o*}+\varepsilon,\ \bsm \beta^{o*})$ is feasible to the dual problem of $\bf{WCEV}(\it O)$, i.e., $UB \leq f_1(\bf x^*) + \eta^o(\bf x^*) \leq f_1(\bf x^*) + \alpha^{o*} + \varepsilon + \sum_{i=1}^m \gamma^i\beta_i^{o*} \leq \ubar{w}^{t_0}+\varepsilon = LB+\varepsilon$. \qed

	Given the definition of $\eta^o(\bf x^*)$, we can infer that $\bf x^*$ is an $\varepsilon-$optimal solution to $\bf{2-Stg \ DRO}$. 	Moreover, when $\hat{\bsm \xi}^o(\bf x^*) \subseteq \widehat{\bsm \xi}^o$, it is clear that we have $UB \leq LB$, rendering $\bf x^*$ an optimal solution to $\bf{2-Stg \ DRO}$.

	Similar to the proof of Theorem \ref{CG-epsilon} on compact sample space $\Xi$, we introduce balls with volume equal to $V(\frac{\delta}{2})$. According to Claim 1, we can conclude that after a finite number of iterations, C\&CG-DRO either returns an optimal solution or an $\varepsilon-$optimal solution when  no feasibility issue occurs.

	Next, we consider the second statement and assume that \textit{Oracle-1} returns an extreme point solution of $\Xi$ if such type of optimal solutions is available. 
	
	\noindent Claim 2: If $\bf{2-Stg \ DRO}$ is infeasible, C\&CG-DRO reports its infeasibility within a finite number of iterations. 
	
	\noindent \textit{Proof of Claim 2}: Since $\Xi$ is a polytope, $\tilde Q_f(\bf x, \xi)$ is convex and $\psi_i(\xi)$ are concave over $\Xi$ for $i\in [m]$, respectively, it follows from Corollary~\ref{CG-exact} that, if $\bf{MMP}$ is feasible, at least one new extreme point from $\mathsf{X\!V}(\Xi)$ (which does not belong to existing $\widehat{\bsm \xi}^f$) will be derived after solving $\bf{WCEV}(F)$. Given that $\mathsf{X\!V}(\Xi)$ is a finite set, it follows that $\bf{MMP}$ will become infeasible within a finite number of iterations.  It hence certifies that $\bf{2-Stg \ DRO}$ is infeasible. \qed
	
	Note that the feasibility check presented in $\bf{WCEV}(F)$ will only be performed  a finite number of times. By combining Claims 1 and 2,  C\&CG-DRO either reports that $\bf{2-Stg \ DRO}$ is infeasible or  returns an optimal (or $\varepsilon-$optimal) one in a finite number of iterations.  
\end{proof}

Next, we present a few results regarding the number of iterations of the algorithm. The following one can be easily obtained by Theorem~\ref{CG-epsilon}.

\begin{cor}
	For a given $\varepsilon$, assume  $\delta$ in the proof Theorem~\ref{CG-epsilon} is available. Then, the number of iterations of C\&CG-DRO before termination is bounded by $\ceil*{V(\hat d+\frac{\delta}{2})/V(\frac{\delta}{2})+|\mathsf{X\!V}(\Xi)|}$.
\end{cor}
When $\cal X$ or $\Xi$ is a finite set, a natural upper bound on the iteration complexity can be derived. Note that it is easy to show that repeated $\bf x$ or $\xi$ leads to $LB=UB$.  
\begin{prop}
	 $(i)$ If $\cal{X}$ is a finite set, the number of iterations is bounded by the cardinality of $\cal X$, i.e., $|\cal X|$; $(ii)$ If sample space $\Xi$ is finite, the number of iterations is bounded by the cardinality of $|\Xi|$.  
\end{prop}

Next, we discuss the case where both the recourse problem and the ambiguity set are linear, noting that it can be exactly solved  within a finite number of iterations.

\begin{prop}
	Assume that the recourse problem is an LP,  $\Xi$ is a polytope, and $\psi_i(\xi)$ are linear over $\Xi$ for $i\in [m]$.  C\&CG-DRO either reports that $\bf{2-Stg\ DRO}$ is infeasible, or converges to its exact solution within $|\mathsf{X\!V}(\Xi)|$ iterations.
\end{prop}
\begin{proof}
	According to the proof of Corollary~\ref{CG-exact} and Theorem~\ref{CCG-DRO-convergence}, since the set of extreme points of $\Xi$, i.e., $|\mathsf{X\!V}(\Xi)|$, is finite, C\&CG-DRO  either reports that $\bf{2-Stg\ DRO}$ is infeasible, or converges to an exact solution that belongs to $\mathsf{X\!V}(\Xi)$. Hence, the whole procedure completes in $|\mathsf{X\!V}(\Xi)|$ iterations.
\end{proof}

\section{Further Investigations on Solving $\bf{2-Stg \ DRO}$}
\label{sect_CCGvariants}
In this section, a couple of new variants of the primary C\&CG-DRO algorithm are presented to support us with a stronger solution capacity. We first develop a variant to compute $\bf{2-Stg \ DRO}$ with  Wasserstein metric-based ambiguity set. Then, we describe another variant to handle MIP ambiguity set, which has not been investigated in the literature yet and can be used to capture more sophisticated uncertainties. 

On one hand, the basic intuition behind the variant for $\bf{2-Stg \ DRO}$ with  Wasserstein metric-based ambiguity set is: given that we are considering a distribution close to an empirical one, the critical scenarios underlying that distribution should also be close to those of the empirical one. Hence, we should identify some non-trivial scenarios around each empirical sample (in an independent fashion) and then aggregate them to provide a basis for selection. We repeat this operation in the framework of C\&CG-DRO. On the other hand, the basic intuition to handle the MIP ambiguity set is: we generate scenarios and build a bilevel MIP (instead of a bilevel LP for convex ambiguity set) to obtain a lower bound approximation of the original $\bf{2-Stg \ DRO}$. As bilevel MIP is very hard to solve, we fix discrete variables in the lower-level problem to convert it into an LP and thus to derive a weaker lower bound. Such a lower bound can be improved through C\&CG-DRO iteratively.

\subsection{A Variant to Handle Wasserstein Metric-Based Ambiguity Set}
\label{subsect_Wasserstein}
Consider two distributions $P$ and $P^e$. The Wasserstein metric between them is defined  as 
\begin{align}\cal W(P, P^e) = \Big\{\inf_{K\in S(P, P^e)}\int_{\Xi}\int_{\Xi^e} \lVert  \xi - \zeta \rVert_{\bf p} K(d \zeta, d \xi)\Big\},
	\label{eq_Wasserstein}
\end{align}
where $S(P, P^e)$ denotes the collection  of joint probability distributions of $\xi$ and $\zeta$  with marginal distributions being $P$ and $P^e$, respectively, and $\lVert\cdot\rVert_{\bf p}$ denotes $L_{\bf p}$ norm. Accordingly, the Wasserstein metric-based ambiguity set,
which has received a lot of attention in the literature of data-driven decision making, is  
\begin{align*}
	\cal{P}^W = \Big\{P \in \cal{M}(\Xi,\cal{F}): \cal W(P, P^e) \leq \mathrm{r}\Big\},
\end{align*}
where $r$ is radius bound imposed on the Wasserstein metric.   
When $P^e$ is defined on a set of empirical samples, i.e.,  $\Xi^e\equiv\{\xi^e_1,\dots,\xi^e_N\}$ with $\{p^e_i\}_{i=1}^{N}$ being its probability mass function, \eqref{eq_Wasserstein} can be simplified to the following one by taking advantage of conditional probabilities~\cite{mohajerin2018data,zhao2018data},
\begin{align}
	\label{eq_Wass_conditional}
	\cal W(P, P^e) = \Big\{\min_{P_i \in \cal M(\Xi, \cal F), i\in [N]} \sum_{i=1}^{N}p_i^e\int_{\Xi} \lVert \xi -  \xi^e_i \rVert P_i(d \xi): \int_{\Xi} P_i(d\xi) = 1\ \ \forall i \in [N]\Big\},
\end{align}
noting that $P_i$ denotes the conditional probability distribution under the condition of $\xi^e_i$. Clearly, we have $P=\sum_{i=1}^Np_i^eP_i$.

\begin{remark}
	We mention that, with the help of the indicator function, Wasserstein metric-based ambiguity set can be described in the form of \eqref{eq_ambiguity} as in the following. 
	\begin{align*}
		\mathcal{P}^W = \Big\{P \in \mathcal{M}(\Xi,\mathcal{F}): &\int_{\Xi}  P(d\mathbf \xi) = 1, \int_{\Xi} \mathbbm{1}_{\{i=k\}} P(d\mathbf \xi) = p^e_k\ \ \forall k\!\in [N],\\
		& \int_{\Xi} \sum_{i=1}^N\sum_{k=1}^Np_i^e\mathbbm{1}_{\{i=k\}}  \Vert \xi - \xi^e_k \Vert_{\bf p} P_i(d\xi) \leq r \Big\},
	\end{align*}
	where the indicator function  $\mathbbm{1}_{\{i=k\}}$ equals 1 if $i=k$ and 0 otherwise. \qed
\end{remark}

With $\cal{P}^W$ defined on $\Xi^e$, it has been noted that the true data-generating probability distribution contained in $\cal{P}^W$ is of a high confidence, and, under some mild assumption, a finite-sample guarantee can be established for a solution to $\cal{P}^W$-based DRO \cite{mohajerin2018data}. 
We next present a new variant of $\bf{C\&CG-DRO}$ to compute such  $\cal{P}^W$-based  $\bf{2-Stg \ DRO}$,  which is referred to as $\bf{C\&CG-DRO}(\cal P^W)$ for simplicity. Also, in the context of Wasserstein metric-based ambiguity set, unless explicitly stated otherwise, we omit the subscript ${\bf p}$ for clarity in our exposition.      
  
\subsubsection{Computing $\bf{WCEV}$ Problems for Given $\bf x^*$}
To describe new variants of solution oracles for $\cal{P}^W$-based $\bf{WCEV}$'s,  we use $\bf{WCEV}(\it F)$ for illustration. Note that all results developed for $\bf{WCEV}(\it F)$ are applicable to its counterpart $\bf{WCEV}(\it O)$ by simply replacing $\tilde Q_f(\xi_j,\bf x^*)$ by $Q(\xi_j,\bf x^*)$.   Specifically, consider the following  $\bf{WCEV}(\it F)$ formulated with respect to $\cal P^W$.
	\begin{align}
	\begin{split}
		 \label{eq_WD_f}
		\bf{WCEV}(\it{F}) : \eta^f(\bf x^*) = \sup\limits_{P_i \in \cal M(\Xi, \cal F), i\in [N]} &\sum_{i = 1}^N p_i^e\int_{\Xi} \tilde Q_f(\bf x^*, \xi) P_i(d \xi)\\
		\text{s.t.}\, &\int_{\Xi} P_i(d\xi) = 1\ \ \forall i \in [N]\\
		&\sum_{i = 1}^Np_i^e \int_{\Xi} \lVert \xi - \xi^e_i \rVert P_i(d \xi) \leq r
	\end{split}
	\end{align}

In the following, we leverage a structural insight and the Pigeonhole principle, a classical result in combinatorics, to develop a novel and compact finite mathematical program that solves $\bf{WCEV}(\it F)$ problem.  
\begin{prop}
	An optimal solution to the following finite mathematical program solves $\bf{WCEV}(\it{F})$ problem in (\ref{eq_WD_f}) exactly.
\begin{align*}
	\bf{WCEV}(\it F)\bf{-FMP}: \max\Big\{&\sum_{i = 1}^N p^e_i[\tilde Q_f(\bf x^*,\xi_{i1}) p_{i1} + \tilde Q_f(\bf x^*,\xi_{i2}) p_{i2}]: p_{i1} + p_{i2} = 1\ \ \forall i \in [N],\\
	&\,\sum_{i = 1}^Np^e_i (\lVert \xi_{i1} - \xi^e_i \rVert p_{i1} + \lVert \xi_{i2}- \xi^e_i \rVert p_{i2}) \leq r,\  \xi_{i1}, \xi_{i2}\in \Xi \ \ \forall i\in [N]
	\Big\}
\end{align*}
\end{prop}
\begin{proof}
	According to Theorm \ref{thm_feasible_EV_sum}, the optimal value of (\ref{eq_WD_f}) can be obtained. Following similar steps to the proof of Proposition \ref{prop_WCEV_OPT} and applying the Richter-Rogosinski Theorem \cite{shapiro2001duality}, \eqref{eq_WD_f} has an optimal distribution with at most $(N + 1)$ scenarios of non-zero probabilities. Hence, given the definition of $P_i$, \eqref{eq_WD_f} is equivalent to 
	 \begin{align*}
	\bf{WCEV}(\it F)\bf{-FMP}^0:    \max\Big\{&\sum_{i = 1}^N \sum_{j=1}^{N+1} p^e_i\tilde Q_f(\bf x^*,\xi_{ij}) p_{ij}: \sum_{j=1}^{N+1}p_{ij} = 1\ \ \forall i \in [N],\\
	 	&\,\sum_{i = 1}^N\sum_{j=1}^{N+1}p^e_i \lVert \xi_{ij} - \xi^e_i \rVert p_{ij} \leq r, \xi_{ij} \in \Xi \ \ \forall i\in [N],\ \forall j\in [N+1] \Big\}.
	 \end{align*} 
  Assume an optimal solution to $\bf{WCEV}(\it F)\bf{-FMP}^0$, denoted by $\big(\xi^*_{ij}, p^*_{ij}\big)_{i\in[N],j\in[N+1]}$, is available. By fixing $\xi_{ij}=\xi^*_{ij}$, $\bf{WCEV}(\it F)\bf{-FMP}^0$ reduces to an LP with $(N+1)$ constraints. Hence, according to the theory of LP, there exists an optimal solution with at most $(N+1)$ probability variables $p_{ij}$ being non-zeros. 

 Given that $\sum_{j=1}^{N+1}p_{ij} = 1$ must be satisfied, it follows that we have at least one $p_{ij}>0$ for every $i$. Then, by the Pigeonhole principle, we further have that no more than two $p_{ij}>0$ for every $i$. As a result, for $i$, we can reduce $(N+1)$ variables  to just two variables, i.e., $p_{i1}$ and $p_{i2}$, through which $\bf{WCEV}(\it F)\bf{-FMP}^0$ becomes $\bf{WCEV}(\it F)\bf{-FMP}$. 
\end{proof}

It is worth pointing out that \eqref{eq_WD_f} (as well as its counterpart for $\bf{WCEV}(\it O)$) actually has a decomposable structure. Note in the objective function of \eqref{eq_WD_f} that 
\begin{align*}
	\sum_{i = 1}^N p_i^e\int_{\Xi} \tilde Q_f(\bf x^*, \xi) P_i(d \xi)
	&=p_1^e\int_{\Xi} \tilde Q_f(\bf x^*, \xi) P_1(d \xi)+\dots+p_N^e\int_{\Xi} \tilde Q_f(\bf x^*, \xi) P_N(d \xi)
	\\
	&=p_1^e\int_{\Xi} \tilde Q_f(\bf x^*, \xi^1) P_1(d \xi^1)+\dots+p_N^e\int_{\Xi} \tilde Q_f(\bf x^*, \xi^N) P_N(d \xi^N),
\end{align*}
which allows us to modify \textit{Oracle-2} to tackle $P_i$ parallelly. Let $\hat{\bsm \xi}^{f}_i=\{\xi^f_{i,1},\dots,\xi^f_{i,n_i}\}$ denote a set of scenarios introduced for $P_i$, $i=1,\dots, N$, and $\widehat{\bsm\xi}^f=(\hat{\bsm \xi}^{f}_1,\dots,\hat{\bsm \xi}^{f}_N)$. We construct the following $\bf{PMP}(\it F)$ defined on top of $\widehat{\bsm\xi}^f$. Note that dual variables are indicated after the colon for each constraint. 
\begin{subequations}
	\begin{align}
		\bf{PMP}(\it F) : \ubar{\eta}^{f*}(\bf x^*, \widehat{\bsm\xi}^f) = \max &\sum_{i = 1}^N p^e_i \sum_{j=1}^{n_i} \tilde Q_f(\bf x^*,\xi^{f}_{i,j}) p_i( \xi^{f}_{i,j})\\
		\text{s.t.}\, &\sum_{j = 1}^{n_i} p_i( \xi^{f}_{i,j}) = 1\ \  \forall i \in [N]\quad : \alpha^{f}_i\\
		&\sum_{i = 1}^N p^e_i \sum_{j=1}^{n_{i}} \lVert \xi^{f}_{ij} - \xi^e_i \rVert p_i( \xi^{f}_{i,j}) \leq r\quad: \beta^f\\
		& p_i(\xi^{f}_{i,j})\geq 0\ \ \forall j\in[n_i],\ \forall i\in[N]	
\end{align}
\end{subequations}
Then, with $\bf{PMP}(\it F)$'s shadow price $(\alpha^{f*}_1, \cdots, \alpha^{f*}_N, \beta^{f*})$ and the stricture of $\bf{PMP}(\it F)$,  we define $\bf{PSP}^{\it i}(\it{F})$ for $i\in [N]$ as follows.  
\begin{equation}
		\bf{PSP}^\it{i}({\it F}) : v^{i*}(\bf x^*, \widehat{\bsm\xi}^f) = \max_{\xi\in \Xi} \ p^e_i\big(\tilde Q_f(\bf x^*, \xi) - \beta^{f*}\lVert \xi^e_{i} - \xi \rVert\big) - \alpha^{f*}_i.
\end{equation}

\begin{remark}
$(i)$ Note that, except $\tilde Q_f(\bf x^*, \xi)$, the objective function of $\bf{PSP}^i({\it F})$ is concave for $L_{\bf p}$-norm if $\bf p\geq 1$ or linearizable (by introducing binary variables) if $\bf p=0$. So, it does not impose additional computational challenge, compared to standard $\bf{PSP}({\it F})$ constructed in Section~\ref{subsect_CG}. \\
$(ii)$ Regarding customization of  \textit{Oracle-2} to handle $\cal P^w$, for shadow price $(\alpha^{f*}_1, \cdots, \alpha^{f*}_N, \beta^{f*})$, all $\bf{PMP}^{\it i}(\it{F})$'s will be solved in an independent or parallel fashion for $i\in [N]$. Once a scenario with  positive reduced cost, denoted by $\xi^{f*}_i$, is  identified by $\bf{PSP}^{\it i}({\it F})$, $\hat{\bsm \xi}_i$ will be augmented as $\hat{\bsm \xi}^f_i=\hat{\bsm \xi}^f_i\cup\{\xi^{f*}_i\}$. Then, $\bf{PMP}(\it F)$, with the updated $\hat{\bsm\xi}^f_i$ for $i\in [N]$, will be used to generate a new shadow price. Hence, the CG procedure terminates if none of $\bf{PMP}^{\it i}(\it{F})$ produces a scenario of positive reduced cost, which is the primary difference to   \textit{Oracle-2} presented in Section~\ref{subsect_CG}. Upon termination, we note that $\hat{\bsm \xi}^f_i\cap\hat{\bsm \xi}^f_j$ might not be empty for $i\neq j$. It suggests that scenarios in the overlap are critical to both $\xi^e_i$ and $\xi^e_j$. \qed
\end{remark}

\subsubsection{Customizing C\&CG-DRO to Handle $\cal P^W$}
Similar to $\widehat{\bsm\xi}^f$, we let $\widehat{\bsm\xi}^o=(\hat{\bsm \xi}^{o}_1,\dots,\hat{\bsm \xi}^{o}_N)$, i.e., optimality sets for  $i\in[N]$, and $\widehat{\bsm \xi} = \widehat{\bsm\xi}^f \cup \widehat{\bsm\xi}^o$. Note that they have been updated up to now after solving $\bf{WCEV}(\it F)$ and $\bf{WCEV}(\it{O})$ in previous iterations, respectively. Then, the main master problem of $\bf{2-Stg \ DRO}$ defined on $\cal P^W$ can be reformulated as:
\begin{subequations}
	\begin{align}
		\bf{MMP}^W : \ubar{w} = \min\limits_{\bf x \in \cal{X}}\,\, &f_1(\bf x) + \eta\\
		&\eta \geq \max\Big\{\sum_{i = 1}^N p^e_i \sum_{\xi \in \hat{\bsm \xi}_i} \eta_{\xi}^o p^o_{\xi}: (p^o_{\xi})_{\xi \in \widehat{\bsm \xi}} \in \cal P^W\Big\}\label{eq_WD_WCEV_o}\\
		&\eta^o_{\xi} = f_2(\bf x, \xi,\bf y_{\xi})\ \ \forall \xi \in \widehat{\bsm \xi}\label{eq_WD_WCEV_c}\\
		&0 \geq \max\Big\{\sum_{i = 1}^N p^e_i \sum_{\xi \in \widehat{\bsm \xi}_i} \eta^f_{\xi} p^f_{\xi}: (p^f_{\xi})_{\xi \in \widehat{\bsm \xi}} \in \cal P^W\Big\}\label{eq_WD_yrecourse_o}\\
		&\eta^f_{\xi} = \lVert\bf{\tilde y}_{\xi}\rVert_{1}\ \  \forall \xi \in \widehat{\bsm \xi}\label{eq_WD_yrecourse_c}\\
		&(\bf y_{\xi}, \bf{\tilde y}_{\xi}) \in \tilde{\cal Y}(\bf x, \xi)\ \ \forall \xi \in \widehat{\bsm \xi}\label{eq_WD_ydomain}
	\end{align}
\end{subequations}

 Similar to $\bf{MMP}_2$ in  \eqref{eq_MMP2-W}, $\bf{MMP}^W$ is a bilevel optimization formulation and can be easily converted into a single-level. We next present the single-level one derived by using the strong duality-based reformulation technique.  
 \begin{align}
 	\begin{split}
 		\bf{MMP}^W : \ubar{w} = \min\Big\{f_1(\bf x) + \eta: & \bf x\in \cal X,\ \eta \geq \sum_{i = 1}^N\alpha^{o}_i + r\beta^o,\\
 		& \Big\{\alpha^o_i +p^e_i\lVert \xi  - \xi^e_i \rVert \beta^o \geq p^e_i\eta^o_{\xi}\ \ \forall \xi \in \hat{\bsm \xi}_i\Big\} \ \ \forall i \in [N]\\
 		&0 \geq \sum_{i = 1}^N\alpha^f_i + r\beta^f, \\
 		& \Big\{\alpha^f_i +p^e_i\lVert \xi - \xi^e_i \rVert \beta^f \geq p^e_i \eta^f_{\xi} \ \ \forall \xi \in \hat{\bsm \xi}_i\Big\}\ \  \forall i \in [N]\\
 		&\eqref{eq_WD_WCEV_c},\ \eqref{eq_WD_yrecourse_c},\ \eqref{eq_WD_ydomain},\ \beta^o \geq 0,\ \beta^f\geq 0 \Big\}.
 	\end{split}
 \end{align}
With the aforementioned $\bf{MMP}^W$ and those defined to solve $\bf{WCEV}$ problems, the customization of  C\&CG-DRO to solve $\mathbf{2-Stg \ DRO}$ with $\cal P^W$ can be easily obtained. As the necessary changes are rather straightforward, we do not provide detailed descriptions.

\begin{remark} 
$(i)$ Actually, as Wasserstein metric-based ambiguity set is defined on top of $\Xi^e$, it is rather straightforward to initialize $\hat{\bsm \xi}^{o}_i$ by $\xi^e_i$. Moreover, given that $p^e_i>0$, the associated recourse problem for $\xi^e_i$ must be feasible for any choice of $\bf x$. In our numerical study, when the radius $r$ is small, this initialization strategy is computationally very effective, while it is less effective when $r$ is large. Also, when $N$ is large, using the whole $\Xi^e$ for initialization might not be computationally effective. If this is the case, we can consider to employ a more adversarial subset of $\Xi^e$ for  initialization.  \\
$(ii)$ We can also modify the procedure of generating Benders cutting planes according described in Remark  \ref{Remark_9} to handle $\cal P^W$.
We consider the recourse problem in \eqref{eq_recourse_LP} and  Benders cuts for optimality for demonstration. Those cuts, generated after solving the current $\mathbf{WCEV}(\it O)$ with $\bf x=\bf x^*$, are in the following forms. 
\begin{align*}
	& \Big\{\alpha^o_i +p^e_i\lVert \xi  - \xi^e_i \rVert \beta^o \geq p_e^i\eta^o_{\xi}\ \ \forall \xi \in \hat{\bsm \xi}^o_i(\bf x^*)\Big\} \ \ \forall i \in [N]\\
	&\Big\{\eta^{o}_{\xi} \geq (\bf b_2 - \bf A_2\bf x - \bf H\xi)^\intercal \pi^*_{\xi} \ \ \forall \xi\in \hat{\bsm \xi}_i^o(\bf x^*)\Big\} \ \ \forall i \in [N]
\end{align*} \qed
\end{remark}

\subsection{An Extension to Handle Mixed Integer Ambiguity Set}\label{subsect_MIPAmbiguity}
As mentioned, so far all ambiguity sets in the DRO literature are assumed to be convex in $P$, although the underlying sample spaces can be either continuous or discrete. Such an assumption is crucial to duality based reformulations, and actually is the enabling structure for all known solution methods. Nevertheless, it seriously restricts our modeling capacity to describe and analyze real world problems. Even some simple situations, as shown next, cannot be represented by any convex set, while they can be captured by mixed integer sets.

\begin{example}
	In addition to constraints in \eqref{eq_ambiguity}, we consider a situation  where the first moment only belongs to one of two intervals, i.e., $[l_1,u_1]$ and $[l_2,u_2]$. It can be seen that it is not convex in $P$, as a convex combination of two legitimate distributions may not yield a legitimate one.  Yet, it can be represented by introducing a new binary variable, $z$, and incorporating the following constraints
	\begin{equation}
		\label{eq_MIP_example}
		l_1z+l_2(1-z)\leq \int_{\Xi} \xi P(d\xi)\leq u_1z+ u_2(1-z),
	\end{equation}
	resulting in a mixed integer representation for the updated ambiguity set $\cal{P}$. $\square$
\end{example}

It is easy to see that for both $z=0$ and $1$, the corresponding $\bf{WCEV}$ problems have optimal discrete probability distributions. Hence the original one does, too. Nevertheless, such mixed integer ambiguity sets are infeasible to the popular duality based approaches. We highlight that they actually can be addressed  by solution methods developed from the primal perspective, through which our modeling and solution capacity on DRO can be greatly improved. Specifically, we mainly consider the following mixed 0-1 ambiguity set to present our algorithm development for the associated $\bf{2-Stg \ DRO}$,  where $\cal{Z}\subseteq \{0,1\}^{n_{\xi^I}}$ denotes the feasible set for binary vector $\bf z$. Note that it can be further extended to handle mixed integer Wasserstein metric-based ambiguity sets.
\begin{align}
	\begin{split}
		\label{eq_IP_ambuiguity}
		\cal{P}^{I} = \Big\{(P, \bf z) \in \cal{M}(\Xi,\cal{F})\times\cal{Z}: &\int_{\Xi} P(d\xi) = 1, \, E_{P}[\psi_i(\xi, \bf z)] \leq \gamma_i(\bf z)  \ \ \forall i \in [m] \Big\}.
	\end{split}
\end{align}

\subsubsection{Computing $\bf{WCEV}$ Problems for A Given $\bf x^*$}

It can be seen that for any $\bf z\in \cal Z $, there are still $m+1$ constraints defining $P$. Hence, the continuous finite mathematical program $\bf{WCEV-FMP}$ presented in~Section \ref{sect_EVWCEV} can be easily extended to the following mixed integer one to compute $\bf{WCEV}(\it F)$ problem. We note again that all results developed for $\bf{WCEV}(\it F)$ problem are applicable to its counterpart $\bf{WCEV}(\it O)$.  

\begin{cor}
	Assume that the sufficient conditions presented in Proposition \ref{prop_WCEV_OPT} are satisfied for every $\bf z\in \cal Z$. Then, $\bf{WCEV}$ problem is equivalent to 
	\begin{align*}
		\begin{split}
			\bf{WCEV}^I(\it F)-\bf{FMP}: \max\Big\{ &\sum_{j=1}^{m+1} p_j\tilde Q_f(\xi_j,\bf x^*):  \sum_{j=1}^{m+1}p_j=1, \ \sum_{j=1}^{m+1}p_j\psi_i(\xi_j,\bf z)\leq \gamma_i(\bf z)\ \ \forall i\in[m], \\ 
			&\xi_j\in \Xi \ \ \forall j\in [m+1],  \ p_j\geq 0 \ \ \forall j\in [m+1], \ \bf z\in \cal Z\Big\}. 
		\end{split}
	\end{align*}
\end{cor}

\setcounter{example}{0}
\begin{example}[Continued]
	As constraints in \eqref{eq_MIP_example} are appended to those in \eqref{eq_ambiguity}, the number of constraints is $m+3$. Then, the corresponding $\bf{WCEV}^I-\bf{FMP}$ is obtained by augmenting the original $\bf{WCEV-FMP}$ with binary variable  $z$, continuous variables $p_{m+2}$ and $p_{m+3}$, and $\xi_{m+2}$ and $\xi_{m+3}$, and with constraint
	\begin{equation}
		l_1z+l_2(1-z)\leq \sum_{j=1}^{m+3} p_j\xi_j\leq u_1z+u_2(1-z). \tag*{\qed}
	\end{equation}
\end{example}

Compared to original $\bf{WCEV(\it F)-FMP}$,  this mixed integer $\bf{WCEV}^I(\it F)-\bf{FMP}$ is even more computationally challenging. With the strong performance of \textit{Oracle-2},  it would be desired to extend it to compute  $\bf{WCEV}$ problems with respect to $\cal P^I$.  As noted earlier, the challenge of a mixed integer master problem (e.g., $\bf{PMP}$ with a mixed integer ambiguity set) can be addressed by customizing the classical B\&P procedures on top of CG~\cite{barnhart1998branch}. Rather than designing and implementing traditional branch-and-bound subroutines, we extend \textit{Oracle-2} in a novel fashion that leverages strong features of professional solvers to minimize the extra development burden. In this subsection, we assume that $\cal P^I$ is not empty for any $\bf z\in \cal Z$.  

First, we  present the new $\bf{PMP}$ defined for $\cal{P}^{I}$, denoted by $\bf{PMP}^I$.
\begin{align}
	\label{eq_PMP-MI}
	\begin{split}
		\bf{PMP}^I(\it F): \underline \eta^*(\bsm\xi^0_n)= \max & \Big\{\sum_{j=1}^{n} \tilde Q_f(\xi^0_j,\bf x) p(\xi^0_j), \ \sum_{j=1}^{n} p(\xi^0_j) = 1, \ p(\xi^0_j)\geq 0 \ \ \forall j\in [n], \\
		& \sum_{j=1}^{n} \psi_i(\xi^0_j,\bf z)p(\xi^0_j)\leq \gamma_i (\bf z)\ \ \forall i \in [m], \ \bf z\in \cal{Z} \Big\}.
	\end{split}
\end{align}

Note that $\bf{PMP}^I$ is a mixed integer program, not a linear program. In the following, we define the pricing subproblem $\bf{PSP}^I$. It has a bilevel optimization structure, which is substantially different from classical pricing subproblems.  
\begin{align}
	\begin{split}
		\bf{PSP}^I(\it F) : v^*(\bsm\xi^0_n) = \max\Bigg\{\tilde Q_f(\xi,\bf x) -   \hat\alpha-&\sum_{i=1}^m\psi_i(\xi,\bf z)\hat\beta_i: \xi\in \Xi,  \ \bf z\in \cal{Z} \label{eq_bilevel_u1}
	\end{split}\\
	\begin{split}
		\big(\hat P, (\hat\alpha, \hat{\bsm\beta})\big)\in \arg & \max \Big\{\sum_{j\in 1}^n \tilde Q_f(\xi^0_i,\bf x) p(\xi^0_i): \sum_{j=1}^{n} p(\xi^0_j) = 1 \\ 
		\sum_{j=1}^{n} \psi_i(\xi^0_j,\bf z) p(\xi^0_j) &\leq \gamma_i(\bf z)\ \ \forall i \in [m], \  
		p(\xi^0_j)\geq 0 \ \ \forall j\in [n]\Big\}\Bigg\}. \label{eq_bilevel_1} 
	\end{split}
\end{align}
Note that the lower-level problem in \eqref{eq_bilevel_1} is actually the continuous portion of $\bf{PMP}^I(\it F)$. With a slight abuse of notation, we let $\big(\hat P, (\hat\alpha, \hat{\bsm\beta})\big)$ denote a pair of  optimal primal and dual  solutions to \eqref{eq_bilevel_1} for $\bf z$ selected by the upper-level problem. So, the overall bilevel optimization problem seeks to choose  $(\xi,\bf z)$ that maximizes the reduced cost based on $(\hat\alpha, \hat{\bsm \beta})$ feedback from the lower-level problem. 
We can further furnish this bilevel optimization problem with the following inequality. It stipulates that, with new $\xi$, the upper bound on the WCEV should be larger than or equal to the optimal value of $\bf{PMP}^I(\it F)$, which helps us to reduce the generation of unnecessary columns in the execution of \textit{Oracle-2}. 

$$\sum_{j\in 1}^n \tilde Q_f(\xi^0_j,\bf x)\hat p(\xi^0_j)+\tilde Q_f(\xi,\bf x) - \hat\alpha-\sum_{i=1}^m\psi_i(\xi,\bf z)\hat\beta_i\geq \underline \eta^*(\bsm\xi^0_n). \label{eq_bilevel_3}$$

Indeed,  bilevel optimization formulation $\bf{PSP}^I$ can be treated as an integration of the master and subproblems of the conventional CG algorithm, where all $\bf z$'s in  $\cal Z$ are evaluated in \eqref{eq_bilevel_u1} when maximizing the reduced cost. Hence, the next result follows directly, which extends that presented in Proposition \ref{prop_CG_basic}.   

\begin{prop}
	\label{prop_CG_MIP_basic}
	Suppose that $\bf{PMP}^I(\it F)$ is feasible and both $\bf{PMP}^I(\it F)$ and $\bf{PSP}^I(\it F)$ are solved to optimality. We have 
	\begin{equation}
		\underline{\eta}^*(\bsm \xi^0_n)  \leq \max_{ P\in \cal P^I} E_P[\tilde Q_f(\bf \xi)]\leq \underline{\eta}^*(\bsm\xi^0_n) +v^*(\bsm\xi^0_n). \tag*{\qed}
	\end{equation}
\end{prop}

Regarding the modification of \textit{Oracle-2} to handle $\cal P^I$, it can be done by simply using $\bf{PMP}^I(\it F)$ and $\bf{PSP}^I(\it F)$ to replace their counterparts, and removing the shadow price output in \textbf{Step 2}. The resulting variant is referred to as \textit{Oracle-2}$^I$.  
\begin{remark}
	\label{rem_bilevel}
	$(i)$ When solving bilevel optimization problem $\bf{PSP}^I$ by a contemporary MIP solver through its singe-level reformulation (see Appendix A.2), it is common to have both optimal primal and dual solutions, i.e., $\big(\hat P, (\hat\alpha, \hat{\bsm\beta})\big)$ are available. Hence, \textit{Oracle-2}$^I$ is rather straightforward to implement, compared to the traditional B\&P method.  \\  
	$(ii)$
	In our computational study, a hybrid implementation including both \textit{Oracle-2} and variant  \textit{Oracle-2}$^I$ is often computationally more efficient. 	Specifically, let $\cal P(\bf z)$ denote the ambiguity set for given $\bf z$. For a fixed $\bf z$ and the associated $\cal P(\bf z)$,  we first run \textit{Oracle-2} to generate a set of $\xi$'s until $v^*(\bsm\xi^0_n)$ becomes 0. Then, \textit{Oracle-2}$^I$ is called to derive a new $\bf z$ with positive $v^*(\bsm\xi^0_n)$, which allows us to switch back to run \textit{Oracle-2}. We repeat those steps until no $\bf z$ with $v^*(\bsm\xi^0_n)>0$, which terminates the whole procedure. Note that this fashion of implementation reduces the number of calls to solve bilevel optimization problems.   	\qed
\end{remark}

\subsubsection{Customizing C\&CG-DRO to Handle $\cal P^I$}

With an augmented $\widehat{\bsm\xi}$, we can build a main master problem with respect to $\cal P^I$, which is a a bilevel optimization formulation with MIP lower-level problem(s). Its optimal solution  can be obtained by employing a C\&CG variant specialized  for such type of bilevel MIP optimization \cite{zeng2014solving}. Nevertheless, rather than searching deeply for strongest $\bf x$ by solving a bilevel MIP program, we can build and solve the following main master problem that is simpler and can be solved more efficiently.      

Note that, for a given $\bf x^*$, computing $\bf{WCEV}(\it F)$ (and $\bf{WCEV}(\it O)$, respectively) problem actually returns $\big(\bf z^{f*}(\bf x^*),\hat{\bsm \xi}^f(\bf x^*)\big)$ (and $\big(\bf z^{o*}(\bf x^*),\hat{\bsm \xi}^o(\bf x^*)\big)$, respectively). Hence, similar to $\widehat{\bsm \xi}^f$ and $\widehat{\bsm \xi}^o$, we 
assume that, before a new C\&CG iteration, two sets of $\bf z$'s have been obtained accumulatively from computing $\bf{WCEV}(\it F)$ and $\bf{WCEV}(\it O)$ problems in previous iterations, denoted by $\widehat{\cal Z}^f$ and $\widehat{\cal Z}^o$ respectively. Also, let $\widehat{\cal Z}=\widehat{\cal Z}^f\cup\widehat{\cal Z}^o$, and $\cal P^I(\bf z')$ denote $\cal P^{I}$ with $\bf z$ fixed to $\bf z'$. Then, we construct and consider the following main master problem. As it is an extension of $\bf{MMP}$ in \eqref{eq_MMP2-W}, we refer to it as  $\bf{MMP}^I$.    
\begin{subequations}
\label{eq_MMP_PI}
	\begin{align}
		\bf{MMP}^I_2 : \ubar{w}' = \min\limits_{\bf x \in \cal{X}} \  &f_1(\bf x) + \eta\\
		 \eta \geq & \max\Big\{\sum_{\xi\in \bsm{\widehat \xi}} \eta^o_{\xi}p^o_{\xi,\bf z}: (p^o_{\xi_1,\bf z},\dots, p^o_{\xi_{|\widehat{\bsm \xi}|},\bf z})\in \cal P(\bf z)\Big\}\ \ \forall \bf z\in \widehat{\cal Z}\label{eq_WCEV_MMP1_z}\\
		 0\geq &\max\Big\{\sum_{\xi\in \bsm{\widehat \xi}} \eta^f_{\xi}p^f_{\xi,\bf z}: (p^f_{\xi_1,\bf z},\dots, p^f_{\xi_{|\bsm{\widehat \xi}|},\bf z})\in \cal P(
		\bf z)\Big\}\ \ \forall \bf z\in \widehat{\cal Z}\label{eq_WCEV_MMP2_z}\\
		& \ \ \ \ \ \eqref{eq_MMP_Oscenarios2}, \ \eqref{eq_MMP_Fscenarios2}, \eqref{eq_MMP_yrecourse_relax} \notag
	\end{align}
\end{subequations}
Since the lower-level problems in \eqref{eq_WCEV_MMP1_z} and \eqref{eq_WCEV_MMP2_z} are LPs, $\bf{MMP}^I$ can be converted into a single-level formulation. Especially by applying the duality-based technique, a big-M-free one, similar to \eqref{eq_MMP_duality}, can be obtained to facilitate easy computation.  We also note that \eqref{eq_WCEV_MMP1_z} and \eqref{eq_WCEV_MMP2_z} are defined over $\widehat{\bsm \xi}$ and $\widehat{\cal Z}$ to simplify $\bf{MMP}^I$'s structure. Actually, given the scale of  $|\widehat{\bsm \xi}|\times|\widehat{\cal Z}|$, it could be computationally more friendly to define \eqref{eq_WCEV_MMP1_z} on $\widehat{\bsm\xi}^o$ and $\widehat{\cal Z}^o$ only, which might not be significantly weaker than the current one.

With all master and subproblems revised according to $\cap P^I$, we next list modifications on particular steps of $\bf{C\&CG-DRO}.$
\begin{description}
	\item[\textbf{Step 1}:] Additional initialization includes 
	$\widehat Z^f=\widehat Z^o=\emptyset$. 
	\item[\textbf{Step 2}:] $\bf{MMP}$ is replaced by (the single-level reformulation of)  $\bf{MMP}^I$.
	\item[\textbf{Step 3}:] $\bf{WCEV}(\it F)$ is replaced by $\bf{WCEV}^I(\it F)$, and the variant of the adopted oracle also outputs optimal $\bf z^f(\bf x^*)$.
	\item[Step 4-Case A:] $\bf{WCEV}(\it O)$ and $\bf{MMP}$ are replaced by $\bf{WCEV}^I(\it O)$ and $\bf{MMP}^I$, the adopted oracle also outputs optimal $\bf z^o(\bf x^*)$, and additional update operation $\widehat Z^o=\widehat Z^o\cup\{\bf z^o(\bf x^*)\}$ is included. 
	\item[Step 4-Case B:] $\bf{WCEV}(\it F)$ is replaced by  $\bf{WCEV}^I(\it F)$, and additional update operation $\widehat Z^f=\widehat Z^f\cup\{\bf z^f(\bf x^*)\}$ is included.
\end{description}
With the aforementioned changes, the updated algorithm is referred to as $\bf{C\&CG-DRO(\cal P}^I)$. We mention that it can be further extended to handle Wasserstein metric-based mixed integer ambiguity sets. As this extension is rather a straightforward integration with the results presented in Section \ref{subsect_Wasserstein}, we omit its description to avoid redundancy. 
\begin{remark}
$(i)$ Analyses on the convergence and iteration complexity for $\bf{C\&CG-DRO}(\cal P^W)$ follows directly from Sections \ref{subsect_CG_convergence} and \ref{subsect_conv_comp} as $\cal P^W$ belongs to the general ambiguity set presented in \eqref{eq_ambiguity}. Such analyses for  $\bf{C\&CG-DRO(\cal P}^I)$ can be developed in a way similar to those presented Sections \ref{subsect_CG_convergence} and \ref{subsect_conv_comp}, noting that $\cal P^I$ reduces to the form of $\cal P$ for any fixed $\bf z$ and set $\cal{Z}$ is finite.  \\
$(ii)$We mention that the presented C\&CG-DRO, together with oracles for $\mathbf{WCEV}$ problems, actually provides a strong and flexible platform to compute a broad class of $\mathbf{2-Stg \ DRO}$ problems. On one hand, the ambiguity set does not need to adhere to any standard form or be defined by any specific mathematical structure. For example, it can be captured by mixing moment inequalities and Wasserstein metric-based consideration. Also, it is comparable to solution methodologies developed for RO and SP.  Note that, as long as master and sub- problems can be solved, there are no strict restrictions imposed on the underlying decision-making problem. More importantly, the overall development requires only a basic understanding of probability, LP, and MIP. This simplicity makes it easy to understand, modify, and debug. \qed
\end{remark}
\section{Numerical Studies}
\label{sect_computation}

In this section, we present and discuss numerical results obtained from computational experiments. Our focus is on testing, evaluating and analyzing $\mathbf{C\&CG-DRO}$ and its variants in computing $\mathbf{2-Stg \ DRO}$ instances. The facility location model with different structures or considerations is employed as the testbed, which often arises from various applications in logistics, supply chain and healthcare systems.   Data regarding clients and facilities' locations, distances and basic demands are adopted from \cite{snyder2005reliability}.  All solution methods are implemented by Python 3.6, with professional MIP solver Gurobi 9.5.2 on a Windows PC with E5-1620 CPU and 32G RAM. Unless noted otherwise, the time limit is set to 7200s, and the relative optimality tolerance of any algorithm/solver is set to .5\%.

\subsection{Distributionally Robust Facility Location Models}\label{section 4.1}
Consider a facility location problem that builds $p\geq 1$ facilities. Let $I$ represent the set of client sites, and $J\subseteq I$ the set of potential facility sites. The parameter $c_{ij}$ captures the service cost of a unit demand from client $i$ served by facility $j$, with $c_{ij} = 0$ if $i=j$. Moreover, $d_i$ is the demand of client $i$, and $f_j$ is the fixed construction cost of a facility at $j$ with $F_j$ being its capacity once the facility is established.  The decision maker seeks a solution of the minimum cost consisting of the fixed construction cost and the weighted sum of the service cost for the normal situation (i.e., nominal demand $\bar{\bf d}$ with no disruptions) and the expected worst-case service cost within an ambiguity set. Under the DRO scheme, we investigate the impact of three major factors. They are the sample space, which could be either continuous or discrete, the ambiguity set, which could be either moment- or Wasserstein metric-based, and the recourse problem, which may or may not have the feasibility issue. Hence, there are 8 different combinations that will be used to  generate instances and to perform our computational study. In the following, we present mathematical formulations, with $\rho$ being the weight coefficient. 

The first basic formulation considers continuous random demands $\bf d=(d_1,\dots, d_{|I|})$.  
\begin{align}
    \mathbf{FL-DRO (\bf d)}: \min_{(\mathbf x,\mathbf y) \in \mathcal X} \sum_{j\in J}f_jy_j+ \rho \sum_{i\in I}\sum_{j\in J} c_{ij} x_{ij} + (1-\rho) \max_{P \in \mathcal P_{\bf{d}}} E_P [Q^{\bf d}(\bf y, \mathbf d)]
\end{align}
{with} $\mathcal X = \Big\{(\mathbf x, \mathbf y) \in \mathbb R^{|I|\times |J|} \times \{0, 1\}^{|J|}: \displaystyle\sum_{j\in J} x_{ij} \geq \bar d_i \ \ \forall i\in I, \ \displaystyle\sum_{j\in J} y_j = p,  \ \displaystyle\sum_{i\in I} x_{ij} \leq F_j y_j \ \ \forall j \in J\Big\}$, {and}
\begin{subequations}
	\begin{align*} Q^{\bf d}(\mathbf y, \bf d) = \min \Big\{\sum_{i\in I}\sum_{j\in J} c_{ij} w_{ij}: & \sum_{j\in J} w_{ij} \geq d_i \ \ \forall i \in I, \ \sum_{i\in I} w_{ij} \leq F_j y_j \ \ \forall j \in J, 
	\\ & w_{ij}\geq 0 \ \ \forall i\in I,\ \forall j\in J\Big\}.
	\end{align*}
\end{subequations}
Note that variables $(\mathbf x, \mathbf y)$ in  $\mathcal X$ are continuous and binary respectively, representing demand allocations and yes/no construction decisions. Constraints in $\mathcal X$ require that all nominal demands are satisfied, the demand allocation can only be made if the facility is constructed, and the total allocation to that facility is subject to its capacity. For $Q^{\bf d}(\bf y,\bf d)$, $\bf w$ represents the demand allocation after the randomness of demand is materialized.     

We highlight that when $F_j$ is sufficiently large, the whole formulation reduces to the uncapacitated model and $Q^{\bf d}(\bf y, \bf d)$ is always feasible regardless the choice of $\bf y$. Otherwise, $\bf y$ needs to be selected to ensure the feasibility of $Q^{\bf d}(\bf y,\bf d)$, which requires to eliminate infeasible $\bf y$'s in our computation. We consider both cases in our numerical study to understand and evaluate  the challenge and the impact of $Q$'s feasibility issue in $\mathbf{2-Stg \ DRO}$. As noted earlier, existing algorithms are not able to handle this challenge.

Regarding the underlying ambiguity set, the sample space of random demand is $\mathcal D = \Big\{\mathbf d \in \mathbb R^{|I|}: {d}^{-}_i \leq d_i \leq d^+_i,\, \forall i \in [I], \Big\}$. We mainly consider the following moment- and Wasserstein metric-based ones (using $L_1$ norm), denoted by $	\cal P^m_{\bf d}$ and $\cal P^w_{\bf d}$, respectively.
\begin{align}
	\label{eq_case_ambi_sets_C}
	\cal P^m_{\bf{d}} = \Big\{P \in \mathcal M(\mathcal D, \mathcal F) :  E_{P} [\mathbf d] \leq \tilde{\bm d}\Big\}; \  \ \ 
	\ \mathcal P^W_{\bf d} =  \Big\{P \in \mathcal M(\mathcal D, \mathcal F) : \mathcal W(P, P^e)\leq \rm r_{\bf d}\Big\}
\end{align}
As C\&CG-DRO (and its variants) is generally applicable, more sophisticated ambiguity sets are considered in Section \ref{subsect_compu_complexAMB}.

Empirical distribution $P^e$ in $\mathcal P^W_{\bf d}$ consists of a set of random samples drawn from $\mathcal{D}$, each with equal probability. Besides continuous demand, we also study discrete disruptions that cause a facility to be unavailable. The demands are then served by the survived facilities. The basic formulation for binary random disruptions, $\bf u=(u_1,\dots, u_{|J|})$, is  
\begin{align}
	\mathbf{FL-DRO (\bf u)}: \min_{(\mathbf x,\mathbf y) \in \mathcal X} \sum_{j\in J}f_jy_j+ \rho \sum_{i\in I}\sum_{j\in J} c_{ij} x_{ij} + (1-\rho) \max_{P \in \mathcal P_{\bf u}} E_{P} [Q^{\bf u}(\bf y, \mathbf u)].
\end{align}
For such random factor, we consider up to $k (\leq p-1)$ disruptions in set $J$, and hence the sample space is  $\mathcal U = \{\mathbf u \in \{0, 1\}^{|J|}: \sum_{j\in J}u_j \leq k\, \forall j\in J\}.$ Such type of sample space, although is finite, is generally exponential with respect to $|J|$, rendering the enumeration is practically infeasible.   The corresponding recourse problem is 
\begin{subequations}
	\begin{align}
		Q^{\bf u}(\mathbf y, \mathbf u) = \min\Big\{\sum_{i\in I}\sum_{j\in J} c_{ij} w_{ij}: & 
	\sum_{j\in J} w_{ij} \geq d_i\ \ \forall i \in I, \ 
		 \sum_{i\in I} w_{ij} \leq F_j y_j \ \ \forall j \in J, \\  \sum_{i\in I} & w_{ij} \leq F_j(1-u_j) \ \ \forall j \in J, \ w_{ij}\geq 0 \ \ \forall i\in I,\  \forall j \in J \Big\}.
	\end{align}
\end{subequations}
Similarly, following ambiguity sets are investigated.
\begin{align}
	\label{eq_case_ambi_sets_C}
	\cal P^m_{\bf u} = \Big\{P \in \mathcal M(\mathcal U, \mathcal F) :  E_{P} [\mathbf u] \leq \tilde{\bm k}\Big\}; \  \ \ 
	\ \mathcal P^W_{\bf u} =  \Big\{P \in \mathcal M(\mathcal U, \mathcal F) : \mathcal W(P, P^e)\leq \rm{r}_{\bf u}\Big\}.
\end{align}
Also, empirical distribution $P^e$ includes random samples of equal probability drawn from discrete sample space.

\subsection{Computational Results of Oracles for $\mathbf{WCEV}$ Problems}
We first perform a set of experiments to evaluate two primal oracles for $\mathbf{WCEV}$ problems, including $\mathbf{WCEV}(\it F)$ and $\mathbf{WCEV}(\it O)$. Instances belonging to ``Small'' group are with $|I|=|J|=5$, $p=3$. For the applicable combinations, the size of empirical set is $5$, $\rm r_{\bf d}=2$, $\rm{r}_{\bf u}=0.5$, and $k=2$ if $\mathcal U$ is adopted. For instances belonging to ``Large'' group, all parameters are same except that $|I|=|J|=8$ and the size of the empirical set is $10$. 
 
All results are presented in Table \ref{tab:WCEV_oracles}, with the time limit set to 10 minutes. where label ``T'' indicates the algorithm terminates before generating an optimal solution. Under this situation, the optimality gap is reported if available, or marked as ``-'' otherwise. On one hand, it can be easily seen that \textit{Oracle-1} is practically infeasible to apply. Only 3 out of 16 instances can be solved to optimality, with non-trivial amount of computational times.  On the other hand,  the computational time of \textit{Oracle-2} for all instances is negligible. Overall, we can estimate that \textit{Oracle-2} is drastically faster, potentially by many orders of magnitude. Hence, we just adopt $\textit{Oracle-2}$ in our remaining experiments.

\begin{table}[H]
	\caption{Computational Results of Two Oracles for $\mathbf{WCEV}$ Problems}
	\begin{center}
		\scalebox{0.7}{
			\begin{tabular}{ccc|cccc|cccc|}
				\hline
				\multicolumn{1}{|c}{\multirow{3}{*}{\begin{tabular}{@{}c@{}} Prob.\\ Size\end{tabular}}} & \multicolumn{1}{|c}{\multirow{3}{*}{\begin{tabular}{@{}c@{}} Ambiguity\\ Sets\end{tabular}}} & \multicolumn{1}{c|}{\multirow{3}{*}{\begin{tabular}{@{}c@{}} $\mathbf{WCEV}$\end{tabular}}} & \multicolumn{4}{c|}{Continuous Sample Space $\mathcal D$}& \multicolumn{4}{c|}{Discrete Sample Space $\mathcal U$}\\ \cline{4-11}
				\multicolumn{1}{|c}{} &\multicolumn{1}{|c}{} & & \multicolumn{2}{c|}{\textit{Oracle-1}} & \multicolumn{2}{c|}{\textit{Oracle-2}}& \multicolumn{2}{c|}{\textit{Oracle-1}} & \multicolumn{2}{c|}{\textit{Oracle-2}} \\ \cline{4-11}
				\multicolumn{1}{|c}{} &\multicolumn{1}{|c}{} &\multicolumn{1}{c|}{} & \multicolumn{1}{c}{Time(s)} &\multicolumn{1}{c|}{Gap(\%)} & \multicolumn{1}{c}{Time(s)} &\multicolumn{1}{c|}{Gap(\%)} & \multicolumn{1}{c}{Time(s)} &\multicolumn{1}{c|}{Gap(\%)} & \multicolumn{1}{c}{Time(s)} &\multicolumn{1}{c|}{Gap(\%)} \\ \hline \hline

				\multicolumn{1}{|c}{\multirow{4}{*}{Small}} & \multicolumn{1}{|c}{\multirow{2}{*}{Moment}} & \multicolumn{1}{c|}{$\bf{WCEV}(\it{F})$} &T &249  &\multicolumn{1}{|c}{0.08} &0  &T &140  &\multicolumn{1}{|c}{0.02} &0 \\
				\multicolumn{1}{|c}{ } & \multicolumn{1}{|c}{ } & \multicolumn{1}{c|}{$\bf{WCEV}(\it{O})$}  &T &11  &\multicolumn{1}{|c}{0.07} &0  &113.9 &0  &\multicolumn{1}{|c}{0.04} &0 \\\clineB{2-11}{2.5}
				\multicolumn{1}{|c}{\multirow{1}{*}{}} & \multicolumn{1}{|c}{\multirow{2}{*}{Wasserstein}} & \multicolumn{1}{c|}{$\bf{WCEV}(\it{F})$} &T &- &\multicolumn{1}{|c}{0.13} &0  &23.5 &0  &\multicolumn{1}{|c}{0.01} &0 \\
				\multicolumn{1}{|c}{ } & \multicolumn{1}{|c}{ } & \multicolumn{1}{c|}{$\bf{WCEV}(\it{O})$}  &T &89  &\multicolumn{1}{|c}{0.9} &0  &6.03 &0  &\multicolumn{1}{|c}{0.03} &0 \\\hline
				
				\multicolumn{1}{|c}{\multirow{4}{*}{Large}} & \multicolumn{1}{|c}{\multirow{2}{*}{Moment}} & \multicolumn{1}{c|}{$\bf{WCEV}(\it{F})$} &T &353  &\multicolumn{1}{|c}{0.00} &0  &T &1604  &\multicolumn{1}{|c}{0.03} &0 \\
				\multicolumn{1}{|c}{ } & \multicolumn{1}{|c}{ } & \multicolumn{1}{c|}{$\bf{WCEV}(\it{O})$}  &T &56.5  &\multicolumn{1}{|c}{0.44} &0  &T &385  &\multicolumn{1}{|c}{0.14} &0 \\\clineB{2-11}{2.5}
				\multicolumn{1}{|c}{\multirow{1}{*}{}} & \multicolumn{1}{|c}{\multirow{2}{*}{Wasserstein}} & \multicolumn{1}{c|}{$\bf{WCEV}(\it{F})$} &T &- &\multicolumn{1}{|c}{0.23} &0  &T &14.5  &\multicolumn{1}{|c}{0.05} &0 \\
				\multicolumn{1}{|c}{ } & \multicolumn{1}{|c}{ } & \multicolumn{1}{c|}{$\bf{WCEV}(\it{O})$}  &T &578  &\multicolumn{1}{|c}{3.99} &0  &T &- &\multicolumn{1}{|c}{0.07} &0 \\\hline
		\end{tabular}}
	\end{center}
	\label{tab:WCEV_oracles}
\end{table}

\subsection{Computational Results of the Uncapacitated Case}
Recall that the existing approaches of implementing C\&CG and Benders decomposition  to solve $\mathbf{2-Stg \ DRO}$ are referred to as basic C\&CG and basic Benders decomposition methods, respectively. In this subsection, we set $F_j$ large than all possible total demands for all $j$, i.e., the case without capacity restrictions. So, there is no feasibility challenge associated with the recourse problem, allowing us to benchmark our work with respect to the current literature. Next, we first investigate instances of $\mathbf{FL-DRO(\bf d)}$ and then those of $\mathbf{FL-DRO(\bf u)}$. Regarding moment-based ambiguity sets, we introduce parameter $\mathfrak{r}$ and a random vector  $\mathfrak{v}$ to generate random $\tilde{\bm d}$ or $\tilde{\bm k}$: $\tilde {\bm d}=\bar{\bf d}(1+\mathfrak{r}\times\mathfrak{v})$ with $\mathfrak{v}\in [-0.5*\bf 1, 0.5*\bf 1]$ and $\tilde {\bm k}=\mathfrak{r}\times\mathfrak{v}$ with $\mathfrak{v}\in [\bf 0, \bf 1]$. Also, we have $\rho=0.5$ in all computational studies.  

\subsubsection{Results for Instances of Continuous Sample Space}
For instances of $\mathbf{FL-DRO(\bf d)}$ with moment-based ambiguity set $\mathcal P^m_{\bf d}$, numerical results of basic C\&CG, C\&CG-DRO, basic Benders, and Benders-DRO algorithms are presented in Table~\ref{tbl:2-Stg-FLDDRO(M)}.
Columns ``LB'', ``UB'', ``Iter.'' and ``$|\widehat{\bsm \xi}|$''  report the lower and upper bounds, the number of main iterations and the number of scenarios generated  when an algorithm terminates, respectively.  Similar to Table \ref{tab:WCEV_oracles}, we mark the entry with ``T'' for ``Time(s)'' or ``-'' for ``Gap(\%)'' if the relevant data is not available.  
Obviously, Benders type of algorithms are practically infeasible, noting that they are extremely uncompetitive compared to their C\&CG counterparts. As for basic C\&CG and C\&CG-DRO, we also observe that the latter exhibits a clear and consistent advantage. For relatively easy instances that basic C\&CG can solve within a few minutes, C\&CG-DRO is faster by an order of magnitude. For more challenging instances, C\&CG-DRO can be up to several hundred times faster. One explanation for such drastically different performances is that basic C\&CG needs to perform many more iterations, given that it only generates a single scenario for every iteration. Another explanation is that C\&CG-DRO identifies critical scenarios in an effective fashion. Note that it produces much fewer scenarios for main master problems, compared to basic C\&CG,  that are still able to provide a strong support to define  worst-case distributions.
 Similar observations can be made between Basic Benders and Benders-DRO, i.e., Benders decomposition with \textit{Oracle-2} to generate Benders cuts, where the latter one is also much faster. Also, we observe that Benders type of algorithms always produce a huge number of Benders cuts, which, however, are rather ineffective in the derivation of optimal solutions.  
 
\begin{figure}[H]
	\bigskip
	\begin{tikzpicture}[scale=0.85]
		\begin{axis}[
			xlabel=$ln(\mbox{time})$,
			ylabel= Value of objective function,
			xmin=-5,xmax=9,
			ymin=0,ymax=3000,
			ytick={0,500,...,3000},
			xtick={-5,-3,...,9},
			grid,
			grid style={dotted},
			legend style={at={(0.5,1.15)},  
				anchor=north,legend columns=2},
			legend entries = { Basic C\&CG,  C\&CG-DRO}]
			
			\addplot [line width=0.5mm, red] table {Demand_FLDIUCCG_UB.dat};
			\addplot [line width=0.5mm, blue] table {Demand_FLDIUCCGCG_UB.dat};
			\addplot [line width=0.5mm, red] table {Demand_FLDIUCCG_LB.dat};
			\addplot [line width=0.5mm, blue] table {Demand_FLDIUCCGCG_LB.dat};
			\addplot[color=blue, solid, fill = brown, mark=10-pointed star, mark size=3pt] coordinates {(2.933, 122.76)};
			\addplot[color=red, solid, fill = orange, mark=*, mark size=3pt] coordinates {(8.63, 122.76)};
		\end{axis}
		\node at (rel axis cs:0.5,-0.25) {\small (a) Convergence over Time};
	\end{tikzpicture}
	\qquad
	\begin{tikzpicture}[scale=0.85]
		\begin{axis}[
			xlabel=Number of iterations,
			ylabel=Vaule of objective function,
			xmin=0,xmax=120,
			ymin=0,ymax=3000,
			ytick={0,500,...,3000},
			xtick={0, 20, ..., 120},
			grid,
			grid style={dotted},
			legend style={at={(0.5,1.15)},  
				anchor=north,legend columns=2},
			legend entries = {Basic C\&CG, C\&CG-DRO}]
			
			\addplot [line width=0.5mm, red] table {Demand_FLDIUCCG_UB_it.dat};
			\addplot [line width=0.5mm, blue] table {Demand_FLDIUCCGCG_UB_it.dat};
			\addplot [line width=0.5mm, red] table {Demand_FLDIUCCG_LB_it.dat};
			\addplot [line width=0.5mm, blue] table {Demand_FLDIUCCGCG_LB_it.dat};
			\addplot[color=blue, solid, fill = blue, mark=10-pointed star, mark size=3pt] coordinates {(3, 122.76)};
			\addplot[color=red, solid, fill = orange, mark=*, mark size=3pt] coordinates {(110, 122.76)};
		\end{axis}
		\node at (rel axis cs:0.5,-0.25) {\small (b) Convergence over Iterations};
	\end{tikzpicture}
	\caption{Convergence Plots of Computing $\mathbf{FL-DRO}(\bf u)$ with $\cal P^m_{\bf d}$}\label{FLD-case_moment}
\end{figure}
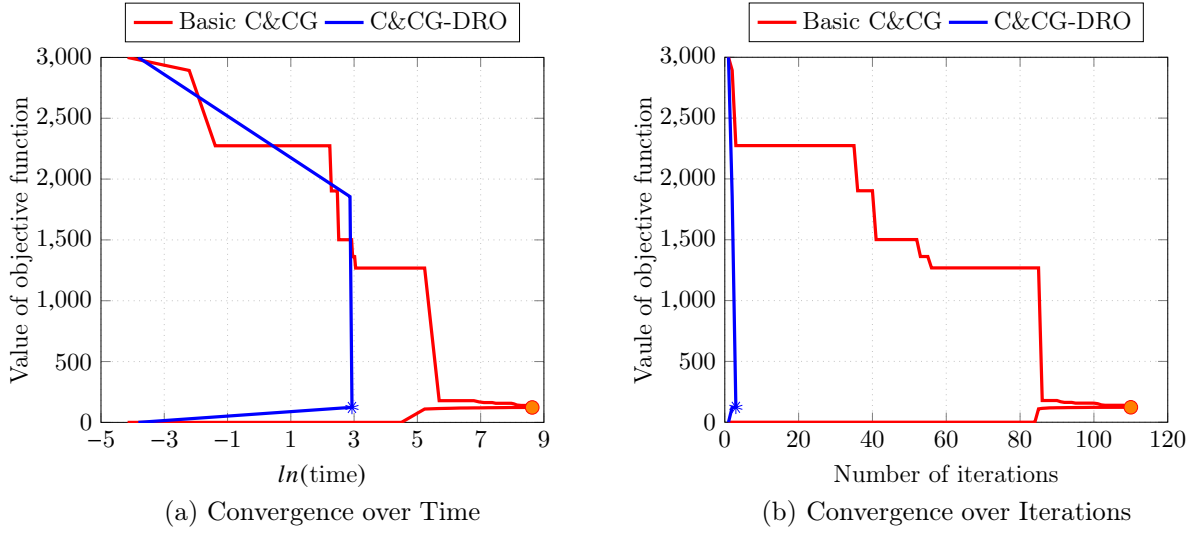

 The overall performance profiles of all these four algorithms are plotted in Figure \ref{fig:4_algorithms}, where the dominance of C\&CG-DRO is straightforward. To deeply understand algorithms' dynamic behaviors, 
Figure \ref{FLD-case_moment} presents the convergence trajectories of basic C\&CG and C\&CG-DRO for the case where $|I|=30$, $\mathfrak{r}=0.65$ and $p=6$. We do not include  those of Benders type of algorithms due to their very weak performances. It can be seen that C\&CG-DRO quickly converges to an optimal solution, especially with the optimal solution identified within a couple of iterations. This observation also confirms the importance of solving $\mathbf{WCEV}$ problem. It not only provides an exact evaluation for a first-stage solution, but also generates a set of scenarios that are critical to select a high quality first-stage solution.     

The aforementioned pattern among those algorithms holds for computational results of instances with Wasserstein metric-based ambiguity set $\cal P^W_{\bf d}$, which are presented in Table~\ref{tbl:2-Stg-FLDDRO(W)}. Typically, C\&CG-DRO is faster than basic C\&CG by 1 to 2 orders of magnitude, while Benders type of algorithms remain very  uncompetitive. It is also interesting to point out, regardless the solution algorithm we employed,  one difference between results for $\mathbf{FL-DRO(\bf d)}$ with  $\mathcal P^m_{\bf d}$ and with $\cal P^W_{\bf d}$: the number of (main) iterations for the former one is much larger than that of the latter one. We believe that the reason behind is that the empirical set underlying $\cal P^W_{\bf d}$ has a determinant impact on its structure. When more samples are available, it is easier to identify the worst-case distribution, which is observed in Table~\ref{tbl:2-Stg-FLDDRO(W)}. Moreover, it is worth highlighting that for some instances, $|\widehat{\bsm \xi}|$ is equal to the product between the number of iterations and $N+1$, e.g., instances (in the form of $I-\rm{r}_{\bf d}-N-p$) $15-10-50-6$ and $30-10-100-10$. Recall that  the theoretical analysis presented in Propositions \ref{prop_WCEV_OPT} and the structure of $\cal P^W_{\bf d}$ indicate that the number of scenarios with non-zero probability is not more than $N+1$ (for every iteration). Clearly, it is verified by the results of those instances, which actually are obtained by \textit{Oracle-2} that does not depend on Proposition \ref{prop_WCEV_OPT}.

 	\begin{table}[H]
 	\caption{Computational Results of Uncapacitated Models with $\cal P^m_{\bf d}$}
 		\begin{center}
 			\scalebox{0.72}{
 				\begin{tabular}{c|cc|cccccc|cccccc|}
 					\hline
 					\multicolumn{1}{|c|}{\multirow{2}{*}{$|I|$}} & \multicolumn{1}{c}{\multirow{2}{*}{$\mathfrak{r}$}} & \multicolumn{1}{c|}{\multirow{2}{*}{$p$}} & \multicolumn{6}{c|}{Basic C\&CG}& \multicolumn{6}{c|}{C\&CG-DRO}\\ \cline{4-15}
 					\multicolumn{1}{|c|}{} & & & {LB} & {UB} & {Gap(\%)} & {Iter.} & {$|\bsm{\widehat \xi}|$} & {Time(s)}& {LB} & {UB} & {Gap(\%)} & {Iter.} & {$|\bsm{\widehat \xi}|$} & {Time(s)}\\ \hline \hline
 					\multicolumn{1}{|c}{\multirow{7}{*}{15}} & \multicolumn{1}{|c}{\multirow{2}{*}{0.65}} & \multicolumn{1}{c|}{6} &51.3 &51.3 &0  &52 &52 &18.1 &51.3 &51.3 &0  &3 &13 &1.9 \\
 					\multicolumn{1}{|c}{ } & \multicolumn{1}{|c}{ } & \multicolumn{1}{c|}{10}  &16.5 &16.5 &0  &71 &71 &9.7 &16.5 &16.5 &0  &3 &12 &1.4 \\\clineB{2-15}{2.5}
 					\multicolumn{1}{|c}{ } & \multicolumn{1}{|c}{\multirow{2}{*}{0.80}} & \multicolumn{1}{c|}{6} &49.2 &49.2 &0  &83 &83 &16.9 &49.1 &49.2 &0  &3 &19 &1.7 \\
 					\multicolumn{1}{|c}{ } & \multicolumn{1}{|c}{ } & \multicolumn{1}{c|}{10} &15.8 &15.8 &0  &50 &50 &6.0 &15.8 &15.8 &0  &3 &14 &1.5  \\\clineB{2-15}{2.5}
 					\multicolumn{1}{|c}{ } & \multicolumn{1}{|c}{\multirow{2}{*}{0.95}} & \multicolumn{1}{c|}{6} &47.2 &47.2 &0  &58 &58 &15.1 &47.2 &47.2 &0  &4 &22 &2.7 \\
 					\multicolumn{1}{|c}{ } & \multicolumn{1}{|c}{ } & \multicolumn{1}{c|}{10} &15.1 &15.1 &0  &66 &66 &7.9 &15.1 &15.1 &0  &3 &12 &1.2 \\\clineB{2-15}{2.5}
 					\multicolumn{1}{|c}{ } & \multicolumn{2}{|c|}{\textbf{Average}}  & & &\textbf{0} &\textbf{63.3} &\textbf{63.3} &\textbf{12.3} & & &\textbf{0} &\textbf{3.2} &\textbf{15.3} &\textbf{1.7} \\\cline{1-1}\clineB{2-15}{2.5}

 					\multicolumn{1}{|c}{\multirow{7}{*}{30}} & \multicolumn{1}{|c}{\multirow{2}{*}{0.65}} & \multicolumn{1}{c|}{6} &122.8 &122.8 &0  &109 &109 &5529.5 &122.7 &122.8 &0  &3 &29 &20.3\\
 					\multicolumn{1}{|c}{ } & \multicolumn{1}{|c}{ } & \multicolumn{1}{c|}{10}  &70.8 &70.8 &0  &40 &40 &120.4 &70.8 &70.8 &0  &3 &29 &8.7 \\\clineB{2-15}{2.5}
 					\multicolumn{1}{|c}{ } & \multicolumn{1}{|c}{\multirow{2}{*}{0.80}} & \multicolumn{1}{c|}{6} &118.3 &118.7 &0  &46 &46 &770.0 &118.2 &118.3 &0  &3 &41 &17.7  \\
 					\multicolumn{1}{|c}{ } & \multicolumn{1}{|c}{ } & \multicolumn{1}{c|}{10} &68.2 &68.2 &0  &47 &47 &191.2 &68.2 &68.2 &0  &3 &25 &8.8 \\\clineB{2-15}{2.5}
 					\multicolumn{1}{|c}{ } & \multicolumn{1}{|c}{\multirow{2}{*}{0.95}} & \multicolumn{1}{c|}{6} &113.9 &114.1 &0  &76 &76 &3782.2 &113.9 &113.9 &0  &3 &43 &43.6 \\
 					\multicolumn{1}{|c}{ } & \multicolumn{1}{|c}{ } & \multicolumn{1}{c|}{10} &65.1 &65.1 &0  &58 &58 &337.3 &65.1 &65.1 &0  &3 &27 &7.7 \\\clineB{2-15}{2.5}
 					\multicolumn{1}{|c}{ } & \multicolumn{2}{|c|}{\textbf{Average}}  & & &\textbf{0} &\textbf{62.6} &\textbf{62.6} &\textbf{1788.43} & & &\textbf{0} &\textbf{3.0} &\textbf{32.3} &\textbf{17.8} \\\hline
 					
 					\multicolumn{1}{|c|}{\multirow{2}{*}{$|I|$}} & \multicolumn{1}{c}{\multirow{2}{*}{$\mathfrak{r}$}} & \multicolumn{1}{c|}{\multirow{2}{*}{$p$}} & \multicolumn{6}{c|}{Basic Benders}& \multicolumn{6}{c|}{Benders-DRO}\\ \cline{4-15}
 					\multicolumn{1}{|c|}{} & & & {LB} & {UB} & {Gap(\%)} & {Iter.} & {Cuts} & {Time(s)}& {LB} & {UB} & {Gap(\%)} & {Iter.} & {Cuts} & {Time(s)}\\ \hline \hline
 					\multicolumn{1}{|c}{\multirow{7}{*}{15}} & \multicolumn{1}{|c}{\multirow{2}{*}{0.65}} & \multicolumn{1}{c|}{6} &13.5 &97.0 &86 &2889 &2889 &T &51.1 &51.3 &0 &250 &24975 &5977.9 \\
 					\multicolumn{1}{|c}{ } & \multicolumn{1}{|c}{ } & \multicolumn{1}{c|}{10}  &0.0 &38.9 &100 &13646 &13646 &T &5.5 &17.3 &68 &389 &27631 &T  \\\clineB{2-15}{2.5}
 					\multicolumn{1}{|c}{ } & \multicolumn{1}{|c}{\multirow{2}{*}{0.80}} & \multicolumn{1}{c|}{6} &11.1 &94.3 &88 &3231 &3231 &T &49.1 &49.2 &0 &231 &21545 &4503.7  \\
 					\multicolumn{1}{|c}{ } & \multicolumn{1}{|c}{ } & \multicolumn{1}{c|}{10} &0.0 &42.7 &100 &13755 &13755 &T &4.8 &15.8 &70 &390 &29488 &T   \\\clineB{2-15}{2.5}
 					\multicolumn{1}{|c}{ } & \multicolumn{1}{|c}{\multirow{2}{*}{0.95}} & \multicolumn{1}{c|}{6} &9.2 &93.7 &90 &3424 &3424 &T &47.1 &47.2 &0 &235 &23769 &5257.1  \\
 					\multicolumn{1}{|c}{ } & \multicolumn{1}{|c}{ } & \multicolumn{1}{c|}{10} &0.0 &44.5 &100 &13760 &13760 &T &5.6 &15.1 &63 &329 &28968 &T  \\\clineB{2-15}{2.5}
 					\multicolumn{1}{|c}{ } & \multicolumn{2}{|c|}{\textbf{Average}}  & & &\textbf{94} &\textbf{8450.8} &\textbf{8450.8} &\textbf{} & & &\textbf{34} &\textbf{304.0} &\textbf{26062.7} &\textbf{} \\\cline{1-1}\clineB{2-15}{2.5}

 					\multicolumn{1}{|c}{\multirow{7}{*}{30}} & \multicolumn{1}{|c}{\multirow{2}{*}{0.65}} & \multicolumn{1}{c|}{6} &0.0 &248.8 &100 &8063 &8063 &T &1.0 &141.0 &99 &95 &28135 &T  \\
 					\multicolumn{1}{|c}{ } & \multicolumn{1}{|c}{ } & \multicolumn{1}{c|}{10}  &0.0 &171.7 &100 &9423 &9423 &T &0.0 &79.3 &100 &422 &90730 &T  \\\clineB{2-15}{2.5}
 					\multicolumn{1}{|c}{ } & \multicolumn{1}{|c}{\multirow{2}{*}{0.80}} & \multicolumn{1}{c|}{6} &0.0 &270.8 &100 &8108 &8108 &T &0.0 &127.2 &100 &98 &28444 &T   \\
 					\multicolumn{1}{|c}{ } & \multicolumn{1}{|c}{ } & \multicolumn{1}{c|}{10} &0.0 &146.5 &100 &9249 &9249 &T &0.0 &82.4 &100 &421 &91524 &T  \\\clineB{2-15}{2.5}
 					\multicolumn{1}{|c}{ } & \multicolumn{1}{|c}{\multirow{2}{*}{0.95}} & \multicolumn{1}{c|}{6} &0.0 &276.7 &100 &8173 &8173 &T &1.1 &120.8 &99 &110 &31051 &T  \\
 					\multicolumn{1}{|c}{ } & \multicolumn{1}{|c}{ } & \multicolumn{1}{c|}{10} &0.0 &163.1 &100 &8816 &8816 &T &0.0 &79.2 &100 &400 &95203 &T  \\\clineB{2-15}{2.5}
 					\multicolumn{1}{|c}{ } & \multicolumn{2}{|c|}{\textbf{Average}}  & & &\textbf{100} &\textbf{8638.7} &\textbf{8638.7} &\textbf{} & & &\textbf{100} &\textbf{257.7} &\textbf{60847.8} &\textbf{} \\\hline
 			\end{tabular}}
 		\label{tbl:2-Stg-FLDDRO(M)}
 	\end{center}
 \end{table}

 \begin{figure}[H]
 	\bigskip
 	\begin{tikzpicture}[scale=0.55]
 		\begin{axis}[
 			width=12cm,
 			height=8cm,
 			xlabel=$ln(time)$,
 			ylabel=\% of instances sovled,
 			xmin=0,xmax=9,
 			ymin=-0.05,ymax=1.1,
 			ytick={0,0.2,...,1},
 			xtick={0,1,2,...,9},
 			grid,
 			grid style={dotted},
 			legend style={at={(0.5,1.25)},  
 				anchor=north,legend columns=4},
 			legend entries = {Basic C\&CG, C\&CG-DRO, Basic Benders, Benders-DRO}]
 			\addplot [line width=0.5mm, red] table {Demand_FLDIUCCG_UB_number.dat};
 			\addplot [line width=0.5mm, blue] table {Demand_FLDIUCCGDRO_UB_number.dat};
 			\addplot [line width=0.5mm, brown] table {Demand_FLDIUBenders_UB_number.dat};
 			\addplot [line width=0.5mm, green] table {Demand_FLDIUBendersDRO_UB_number.dat};
 			\addplot+[color = red, mark=square*, mark size=2pt] table {Demand_FLDIUCCG_UB_number.dat};
 			\addplot+[color = blue, mark=triangle*, mark size=2pt] table {Demand_FLDIUCCGDRO_UB_number.dat};
 			\addplot+[color = brown, mark=halfcircle*, mark size=2pt] table {Demand_FLDIUBenders_UB_number.dat};
 			\addplot+[color = green, mark=pentagon*, mark size=2pt] table {Demand_FLDIUBendersDRO_UB_number.dat};
 		\end{axis}
 		\node at (rel axis cs:0.5,-0.25){\scriptsize (a) \% of Instances Solved over $ln(time)$};
 	\end{tikzpicture}
 	\qquad
 	\begin{tikzpicture}[scale=0.55]
 		\begin{axis}[  
 			ybar,
 			ylabel=Relative gap,
 			xlabel=Number of sites,
 			yticklabel=\pgfmathprintnumber{\tick}\,$\%$,
 			legend style={at={(0.5,1.2)},  
 				anchor=north,legend columns=4},
 			symbolic x coords={15, 30},  
 			xtick=data,  
 			nodes near coords,  
 			nodes near coords={\pgfkeys{/pgf/fpu}\pgfmathparse{\pgfplotspointmeta}\pgfmathprintnumber{\pgfmathresult}\,\%},
 			nodes near coords style={font=\tiny},
 			width=13.8cm,
 			height=8cm,
 			bar width=17pt,
 			ymin=0, 
 			enlarge x limits=0.5, enlarge y limits={upper=0},]  
 			\addplot [red, fill=red!40!white, postaction={pattern=north east lines}, pattern color = red] coordinates {(15, 0) (30, 0)};  
 			\addplot[blue, fill=blue!40!white, postaction={pattern=north east lines}, pattern color = blue] coordinates {(15, 0) (30, 0)};  
 			\addplot [brown, fill=brown!40!white, postaction={pattern=north east lines}, pattern color = brown] coordinates {(15, 94) (30, 100)};  
 			\addplot [green, fill=green!40!white, postaction={pattern=north east lines}, pattern color = green] coordinates {(15, 34) (30, 100)};  
 			\legend{Basic C\&CG, C\&CG-DRO, Basic Benders, Benders-DRO}  
 		\end{axis}
 		
 		\node at (rel axis cs:0.75,-0.25) {\scriptsize (b) Average Gaps upon Termination};
 	\end{tikzpicture}
 	\caption{Overall performance profiles of Four Algorithms\label{fig:4_algorithms}}
 \end{figure}

	\begin{table}[H]
	\caption{Computational Results of Uncapacitated Models with $\mathcal P_{\mathbf d}^W$}
	\begin{center}
		\scalebox{0.72}{
			\begin{tabular}{c|c|cc|cccccc|cccccc|}
				\hline
				\multicolumn{1}{|c|}{\multirow{2}{*}{\begin{tabular}{@{}c@{}} $I$ \end{tabular}}} & \multicolumn{1}{c|}{\multirow{2}{*}{$\rm{r}_{\bf d}$}} & \multicolumn{1}{c}{\multirow{2}{*}{\begin{tabular}{@{}c@{}}$N$\end{tabular}}} & \multicolumn{1}{c|}{\multirow{2}{*}{$p$}}& \multicolumn{6}{c|}{Basic C\&CG}& \multicolumn{6}{c|}{C\&CG-DRO} \\ \cline{5-16}
				\multicolumn{1}{|c|}{} & & & & {LB} & {UB} & {Gap(\%)} & {Iter.} & {$|\bsm{\widehat \xi}|$} & {Time(s)}& {LB} & {UB} & {Gap(\%)} & {Iter.} & {$|\bsm{\widehat \xi}|$} & {Time(s)} \\ \hline \hline
				\multicolumn{1}{|c}{\multirow{14}{*}{15}} & \multicolumn{1}{|c|}{\multirow{7}{*}{2}} & \multicolumn{1}{c}{\multirow{2}{*}{50}} &6 &62.3 &62.6 &0  &5 &250 &114.1 &62.4 &62.4 &0  &2 &80 &26.0 \\
				\multicolumn{1}{|c}{} & \multicolumn{1}{|c|}{} & \multicolumn{1}{c}{} &10 &20.1 &20.2 &0  &3 &150 &12.9 &20.1 &20.1 &0 &1 &14 &5.1 \\\clineB{3-16}{2.5}
				\multicolumn{1}{|c}{} & \multicolumn{1}{|c|}{} & \multicolumn{1}{c}{\multirow{2}{*}{100}} &6 &62.1 &62.3 &0  &6 &600 &726.7 &62.3 &62.3 &0  &2 &139 &52.5 \\
				\multicolumn{1}{|c}{} & \multicolumn{1}{|c|}{} & \multicolumn{1}{c}{} &10 &20.0 &20.1 &0  &3 &300 &40.3 &20.0 &20.0 &0  &1 &27 &10.4  \\\clineB{3-16}{2.5}
				\multicolumn{1}{|c}{} & \multicolumn{1}{|c|}{} & \multicolumn{1}{c}{\multirow{2}{*}{200}} &6 &62.1 &62.3 &0  &5 &1000 &3077.7 &62.1 &62.1 &0  &2 &185 &98.7 \\
				\multicolumn{1}{|c}{} & \multicolumn{1}{|c|}{} & \multicolumn{1}{c}{} &10 &19.9 &20.0 &0  &3 &600 &134.5 &19.9 &19.9 &0  &1 &31 &16.7 \\\clineB{3-16}{2.5}
				\multicolumn{1}{|c}{} & \multicolumn{1}{|c|}{} & \multicolumn{2}{c|}{\textbf{Average}}  & &\textbf{41.3} &\textbf{0 } &\textbf{4.2} &\textbf{483.3} &\textbf{684.4} & &\textbf{41.2} &\textbf{0 } &\textbf{1.5} &\textbf{79.3} &\textbf{34.9} \\\cline{2-2}\clineB{3-16}{2.5}
				
				\multicolumn{1}{|c}{} & \multicolumn{1}{|c|}{\multirow{7}{*}{10}} & \multicolumn{1}{c}{\multirow{2}{*}{50}} &6 &63.7 &63.9 &0  &5 &250 &214.0 &63.8 &64.0 &0  &3 &153 &70.2 \\
				\multicolumn{1}{|c}{} & \multicolumn{1}{|c|}{} & \multicolumn{1}{c}{} &10 &21.2 &21.2 &0  &3 &150 &12.5 &21.2 &21.2 &0  &1 &51 &5.4 \\\clineB{3-16}{2.5}
				\multicolumn{1}{|c}{} & \multicolumn{1}{|c|}{} & \multicolumn{1}{c}{\multirow{2}{*}{100}} &6 &63.6 &63.8 &0  &5 &500 &770.8 &63.7 &63.8 &0  &3 &303 &179.4 \\
				\multicolumn{1}{|c}{} & \multicolumn{1}{|c|}{} & \multicolumn{1}{c}{} &10 &21.1 &21.2 &0  &3 &300 &43.9 &21.1 &21.1 &0  &1 &101 &12.0 \\\clineB{3-16}{2.5}
				\multicolumn{1}{|c}{} & \multicolumn{1}{|c|}{} & \multicolumn{1}{c}{\multirow{2}{*}{200}} &6 &63.4 &63.7 &0  &4 &800 &1719.3 &63.3 &63.6 &0  &2 &402 &190.5 \\
				\multicolumn{1}{|c}{} & \multicolumn{1}{|c|}{} & \multicolumn{1}{c}{} &10 &21.0 &21.0 &0  &3 &600 &129.5 &21.0 &21.0 &0  &1 &148 &20.9 \\\clineB{3-16}{2.5}
				\multicolumn{1}{|c}{} & \multicolumn{1}{|c|}{} & \multicolumn{2}{c|}{\textbf{Average}}  & &\textbf{42.5} &\textbf{0 } &\textbf{3.8} &\textbf{433.3} &\textbf{481.6} & &\textbf{42.5} &\textbf{0} &\textbf{1.8} &\textbf{193.0}  &\textbf{79.7} \\\cline{1-2}\clineB{3-16}{2.5}

				\multicolumn{1}{|c}{\multirow{14}{*}{30}} & \multicolumn{1}{|c|}{\multirow{7}{*}{2}} & \multicolumn{1}{c}{\multirow{2}{*}{50}} &6 &143.8 &144.4 &0  &4 &200 &1105.4 &144.2 &144.2 &0  &1 &44 &130.0 \\
				\multicolumn{1}{|c}{} & \multicolumn{1}{|c|}{} & \multicolumn{1}{c}{} &10 &84.1 &84.4 &0  &4 &200 &400.4 &84.4 &84.4 &0  &1 &51 &52.4 \\\clineB{3-16}{2.5}
				\multicolumn{1}{|c}{} & \multicolumn{1}{|c|}{} & \multicolumn{1}{c}{\multirow{2}{*}{100}} &6 &143.6 &144.1 &0  &4 &400 &4803.8 &144.0 &144.0 &0  &1 &85 &596.8 \\
				\multicolumn{1}{|c}{} & \multicolumn{1}{|c|}{} & \multicolumn{1}{c}{} &10 &83.9 &84.2 &0  &4 &400 &1561.6 &84.2 &84.2 &0  &1 &101 &222.3 \\\clineB{3-16}{2.5}
				\multicolumn{1}{|c}{} & \multicolumn{1}{|c|}{} & \multicolumn{1}{c}{\multirow{2}{*}{200}} &6 &143.0 &144.8 &1  &3 &600 &T &143.5 &143.5 &0  &1 &77 &1435.9 \\
				\multicolumn{1}{|c}{} & \multicolumn{1}{|c|}{} & \multicolumn{1}{c}{} &10 &83.7 &84.0 &0  &4 &800 &6587.7 &83.9 &83.9 &0  &1 &109 &487.9 \\\clineB{3-16}{2.5}
				\multicolumn{1}{|c}{} & \multicolumn{1}{|c|}{} & \multicolumn{2}{c|}{\textbf{Average}}  & &\textbf{114.3} &\textbf{0} &\textbf{3.8} &\textbf{433.3} &\textbf{} & &\textbf{114.0} &\textbf{0 } &\textbf{1.0} &\textbf{77.8} &\textbf{487.6} \\\cline{2-2}\clineB{3-16}{2.5}
				
				\multicolumn{1}{|c}{} & \multicolumn{1}{|c|}{\multirow{7}{*}{10}} & \multicolumn{1}{c}{\multirow{2}{*}{50}} &6 &146.6 &147.2 &0  &3 &150 &1235.2 &147.1 &147.1 &0  &1 &51 &154.6 \\
				\multicolumn{1}{|c}{} & \multicolumn{1}{|c|}{} & \multicolumn{1}{c}{} &10 &86.3 &86.3 &0  &3 &150 &500.5 &86.3 &86.3 &0  &1 &51 &49.2 \\\clineB{3-16}{2.5}
				\multicolumn{1}{|c}{} & \multicolumn{1}{|c|}{} & \multicolumn{1}{c}{\multirow{2}{*}{100}} &6 &146.3 &147.0 &0  &3 &300 &5120.4 &146.8 &146.8 &0  &1 &101 &723.6 \\
				\multicolumn{1}{|c}{} & \multicolumn{1}{|c|}{} & \multicolumn{1}{c}{} &10 &86.1 &86.2 &0  &3 &300 &1664.0 &86.2 &86.2 &0  &1 &101 &216.8 \\\clineB{3-16}{2.5}
				\multicolumn{1}{|c}{} & \multicolumn{1}{|c|}{} & \multicolumn{1}{c}{\multirow{2}{*}{200}} &6 &144.8 &146.7 &1  &2 &400 &T &146.5 &146.5 &0  &1 &201 &2303.9 \\
				\multicolumn{1}{|c}{} & \multicolumn{1}{|c|}{} & \multicolumn{1}{c}{} &10 &85.8 &86.3 &1  &2 &400 &T &86.0 &86.0 &0  &1 &201 &727.3 \\\clineB{3-16}{2.5}
				\multicolumn{1}{|c}{} & \multicolumn{1}{|c|}{} & \multicolumn{2}{c|}{\textbf{Average}}  & &\textbf{116.6} &\textbf{0} &\textbf{2.7} &\textbf{283.3} &\textbf{} & &\textbf{116.5} &\textbf{0 } &\textbf{1.0} &\textbf{117.7} &\textbf{695.9}\\\hline
		\end{tabular}}
	\end{center}\label{tbl:2-Stg-FLDDRO(W)}
	\end{table}

\subsubsection{Results for Instances of Discrete Sample Space}
Given the extremely poor performance of Benders type of algorithms, we do not consider them in the remainder of this section, unless otherwise noted. For $\mathbf{FL-DRO (\bf u)}$, the number of sites under consideration is smaller than that for $\mathbf{FL-DRO (\bf d)}$ due to its computational challenge.  Except for parameter $k$,  other parameters are kept the same as those in Tables \ref{tbl:2-Stg-FLDDRO(M)} and \ref{tbl:2-Stg-FLDDRO(W)}. Computational results for $\mathbf{FL-DRO (\bf u)}$ with $\cal P^m_{\bf u}$ and $\cal P^W_{\bf u}$ are presented in  Tables \ref{tbl:2-Stg-FLRDRO(m)} and \ref{tbl:2-Stg-FLRDRO(W)}, respectively.

From Tables \ref{tbl:2-Stg-FLRDRO(m)} and \ref{tbl:2-Stg-FLRDRO(W)}, it can be seen that both algorithms are more sensitive to parameters determining the size of instances, including $|I|$ and $N$. A few instances have not been solved to the exactness by either algorithm, indicating that $\mathbf{FL-DRO (\bf u)}$ is more complex than $\mathbf{FL-DRO (\bf d)}$. The number of iterations is generally much greater than that of its counterpart $\mathbf{FL-DRO (\bf d)}$. The same observation holds for the number of scenarios involved in. One explanation is that the non-convex structure of the discrete sample space renders the associated ambiguity sets more intricate. As a result, a large number of iterations needs to be performed to generate many scenarios that are critical to capture the worst-case distributions. Especially for Wasserstein metric-based $\cal P^W_{\bf u}$, we note that the number of scenarios generated could be larger than 2,000, which is drastically more than that for moment-based $\cal P^m_{\bf u}$. As the number of scenarios largely determines the scale of $\mathbf{MMP}$,  it deserves a deep investigation to develop efficient algorithms for $\mathbf{MMP}$ when the number of empirical samples is large. Another observation from handling $\cal P^W$ is that the numbers of iterations could become smaller with respect to larger empirical sets, which verifies the effectiveness to employ those samples for initialization.     
  Between basic C\&CG and C\&CG-DRO, they have roughly the same performance for simple instances. Yet, for more challenging instances, C\&CG-DRO demonstrates a much stronger capacity. It either can solve instances that cannot be computed by basic C\&CG, or produces solutions with significantly smaller gaps.   

	\begin{table}[H]
	\caption{Computational Results of Uncapacitated Models with $\mathcal P_{\mathbf u}^m$}
	\begin{center}
		\scalebox{0.72}{
			\begin{tabular}{c|ccc|cccccc|cccccc|}
				\hline
				\multicolumn{1}{|c|}{\multirow{2}{*}{\begin{tabular}{@{}c@{}} $|I|$\end{tabular}}} & \multicolumn{1}{c}{\multirow{2}{*}{$\mathfrak{r}$}} & \multicolumn{1}{c}{\multirow{2}{*}{$p$}} & \multicolumn{1}{c|}{\multirow{2}{*}{$k$}} & \multicolumn{6}{c|}{Basic C\&CG}& \multicolumn{6}{c|}{C\&CG-DRO}\\ \cline{5-16}
				\multicolumn{1}{|c|}{} & & & & {LB} & {UB} & {Gap(\%)} & {Iter.} & {$|\bsm{\widehat \xi}|$} & {Time(s)}& {LB} & {UB} & {Gap(\%)} & {Iter.} & {$|\bsm{\widehat \xi}|$} & {Time(s)} \\ \hline \hline
				\multicolumn{1}{|c}{\multirow{7}{*}{15}} & \multicolumn{1}{|c}{\multirow{2}{*}{0.65}} & \multicolumn{1}{c}{6} & \multicolumn{1}{c|}{3} &230.1 &230.1 &0  &120 &120 &1105.3 &229.1 &230.1 &0  &67 &99 &676.1  \\
				\multicolumn{1}{|c}{ } & \multicolumn{1}{|c}{ } & \multicolumn{1}{c}{10}  & \multicolumn{1}{c|}{5} &182.4 &182.4 &0  &164 &164 &790.3 &181.6 &182.4 &0 &64  &179 &530.5 \\\clineB{2-16}{2.5}
				\multicolumn{1}{|c}{ } & \multicolumn{1}{|c}{\multirow{2}{*}{0.80}} & \multicolumn{1}{c}{6} & \multicolumn{1}{c|}{3} &236.3 &236.3 &0  &105 &105 &640.9 &235.5 &236.3 &0  &30 &59 &124.5 \\
				\multicolumn{1}{|c}{ } & \multicolumn{1}{|c}{ } & \multicolumn{1}{c}{10}  & \multicolumn{1}{c|}{5} &195.1 &195.1 &0  &87 &87 &173.0 &194.2 &195.1 &0  &26 &81 &109.1  \\\clineB{2-16}{2.5}
				\multicolumn{1}{|c}{ } & \multicolumn{1}{|c}{\multirow{2}{*}{0.95}} & \multicolumn{1}{c}{6} & \multicolumn{1}{c|}{3} &242.6 &242.6 &0  &93 &93 &417.9 &242.2 &242.6 &0  &12 &25 &15.8 \\
				\multicolumn{1}{|c}{ } & \multicolumn{1}{|c}{ } & \multicolumn{1}{c}{10}  & \multicolumn{1}{c|}{5} &207.8 &207.8 &0  &45 &45 &51.7 &207.8 &207.8 &0  &13 &40 &35.8 \\\clineB{2-16}{2.5}
				\multicolumn{1}{|c}{ } & \multicolumn{3}{|c|}{\textbf{Average}} & &\textbf{215.7} &\textbf{0} &\textbf{102.3} &\textbf{102.3} &\textbf{529.9} & &\textbf{215.7} &\textbf{0} &\textbf{35.3} &\textbf{80.5} &\textbf{248.6} \\\cline{1-1}\clineB{2-16}{2.5}

				\multicolumn{1}{|c}{\multirow{7}{*}{20}} & \multicolumn{1}{|c}{\multirow{2}{*}{0.65}} & \multicolumn{1}{c}{6} & \multicolumn{1}{c|}{3} &215.9 &299.2 &28  &107 &107 &T &266.4 &280.3 &5 &74 &118 &T \\
				\multicolumn{1}{|c}{ } & \multicolumn{1}{|c}{ } & \multicolumn{1}{c}{10}  & \multicolumn{1}{c|}{5} &197.2 &237.6 &17 &138 &138 &T &208.1 &227.6 &9 &59 &234 &T \\\clineB{2-16}{2.5}
				\multicolumn{1}{|c}{ } & \multicolumn{1}{|c}{\multirow{2}{*}{0.80}} & \multicolumn{1}{c}{6} & \multicolumn{1}{c|}{3} &220.5 &299.2 &26 &104 &104 &T &280.3 &288.6 &3 &76 &124 &T \\
				\multicolumn{1}{|c}{ } & \multicolumn{1}{|c}{ } & \multicolumn{1}{c}{10}  & \multicolumn{1}{c|}{5} &212.2 &246.6 &14 &167 &167 &T &226.3 &233.2 &3 &59 &219 &T \\\clineB{2-16}{2.5}
				\multicolumn{1}{|c}{ } & \multicolumn{1}{|c}{\multirow{2}{*}{0.95}} & \multicolumn{1}{c}{6} & \multicolumn{1}{c|}{3} &226.2 &299.2 &24 &116 &116 &T &293.3 &296.5 &1 &92 &123 &T\\
				\multicolumn{1}{|c}{ } & \multicolumn{1}{|c}{ } & \multicolumn{1}{c}{10}  & \multicolumn{1}{c|}{5} &237.6 &254.3 &7 &165 &165 &T &239.5 &240.2 &0 &45 &151 &2318.8 \\\clineB{2-16}{2.5}
				\multicolumn{1}{|c}{ } & \multicolumn{3}{|c|}{\textbf{Average}}  & &\textbf{272.7} &\textbf{19} &\textbf{132.8} &\textbf{132.8} && &\textbf{261.1} &\textbf{3} &\textbf{67.5} &\textbf{161.5} & \\\hline
		\end{tabular}}
	\end{center}\label{tbl:2-Stg-FLRDRO(m)}
\end{table}

	\begin{table}[H]
	\caption{Computational Results of Uncapacitated Models with $\mathcal P_{\mathbf u}^W$}
	\begin{center}
		\scalebox{0.72}{
			\begin{tabular}{c|c|ccc|cccccc|cccccc|}
				\hline
				\multicolumn{1}{|c|}{\multirow{2}{*}{\begin{tabular}{@{}c@{}} $|I|$ \end{tabular}}} & \multicolumn{1}{c|}{\multirow{2}{*}{$\rm{r}_{\bf u}$}} & \multicolumn{1}{c}{\multirow{2}{*}{\begin{tabular}{@{}c@{}} $N$ \end{tabular}}} & \multicolumn{1}{c}{\multirow{2}{*}{$p$}}& \multicolumn{1}{c|}{\multirow{2}{*}{$k$}}& \multicolumn{6}{c|}{Basic C\&CG}& \multicolumn{6}{c|}{C\&CG-DRO} \\ \cline{6-17}
				\multicolumn{1}{|c|}{} & & & & & {LB} & {UB} & {Gap(\%)} & {Iter.} & {$|\bsm{\widehat \xi}|$} & {Time(s)}& {LB} & {UB} & {Gap(\%)} & {Iter.} & {$|\bsm{\widehat \xi}|$} & {Time(s)} \\ \hline \hline
				\multicolumn{1}{|c}{\multirow{14}{*}{15}} & \multicolumn{1}{|c|}{\multirow{7}{*}{0.5}} & \multicolumn{1}{c}{\multirow{2}{*}{50}} &6 &3 &201.8 &201.9 &0 &4 &150 &29.4&201.9 &201.9 &0 &3 &91 &25.8  \\
				\multicolumn{1}{|c}{} & \multicolumn{1}{|c|}{} & \multicolumn{1}{c}{} &10 &5 &181.9 &181.9 &0 &4 &138 &113.2 &181.5 &181.9 &0 &4 &104 &162.3 \\\clineB{3-17}{2.5}
				\multicolumn{1}{|c}{} & \multicolumn{1}{|c|}{} & \multicolumn{1}{c}{\multirow{2}{*}{100}} &6 &3 &194.3 &194.3 &0 &4 &292 &72.9 &194.2 &194.5 &0 &3 &251 &49.7 \\
				\multicolumn{1}{|c}{} & \multicolumn{1}{|c|}{} & \multicolumn{1}{c}{} &10 &5 &181.1 &181.9 &0 &3 &263 &176.9 &181.0 &181.1 &0 &3 &162 &271.4 \\\clineB{3-17}{2.5}
				\multicolumn{1}{|c}{} & \multicolumn{1}{|c|}{} & \multicolumn{1}{c}{\multirow{2}{*}{200}} &6 &3 &189.8 &189.8 &0 &4 &568 &204.1 &189.4 &189.8 &0 &4 &551 &173.1 \\
				\multicolumn{1}{|c}{} & \multicolumn{1}{|c|}{} & \multicolumn{1}{c}{} &10 &5 &168.7 &169.5 &0 &3 &390 &365.1 &168.9 &169.6 &0 &2 &246 &418.3 \\\clineB{3-17}{2.5}
				\multicolumn{1}{|c}{} & \multicolumn{1}{|c|}{} & \multicolumn{3}{c|}{\textbf{Average}}  & &\textbf{186.5} &\textbf{0} &\textbf{3.7} &\textbf{300.2} &\textbf{160.3} & &\textbf{186.5} &\textbf{0} &\textbf{3.2} &\textbf{234.2} &\textbf{183.4} \\\cline{2-2}\clineB{3-17}{2.5}
				
				\multicolumn{1}{|c}{} & \multicolumn{1}{|c|}{\multirow{7}{*}{2}} & \multicolumn{1}{c}{\multirow{2}{*}{50}} &6 &3 &234.8 &234.8 &0 &11 &536 &317.6 &234.1 &234.8 &0 &11 &488 &236.4 \\
				\multicolumn{1}{|c}{} & \multicolumn{1}{|c|}{} & \multicolumn{1}{c}{} &10 &5 &200.2 &200.7 &0 &6 &207 &189.9 &200.1 &200.7 &0 &7 &280 &555.9 \\\clineB{3-17}{2.5}
				\multicolumn{1}{|c}{} & \multicolumn{1}{|c|}{} & \multicolumn{1}{c}{\multirow{2}{*}{100}} &6 &3 &229.0 &229.0 &0 &9 &681 &442.0 &228.6 &229.3 &0 &7 &667 &332.2 \\
				\multicolumn{1}{|c}{} & \multicolumn{1}{|c|}{} & \multicolumn{1}{c}{} &10 &5 &198.2 &198.6 &0 &7 &522 &559.7 &198.2 &198.6 &0 &6 &541 &865.9 \\\clineB{3-17}{2.5}
				\multicolumn{1}{|c}{} & \multicolumn{1}{|c|}{} & \multicolumn{1}{c}{\multirow{2}{*}{200}} &6 &3 &226.4 &227.4 &0 &10 &1721 &2934.4 &226.3 &226.4 &0 &9 &1884 &3291.1 \\
				\multicolumn{1}{|c}{} & \multicolumn{1}{|c|}{} & \multicolumn{1}{c}{} &10 &5 &192.8 &193.1 &0 &7 &980 &1214.7 &192.6 &193.4 &0 &6 &851 &1985.6 \\\clineB{3-17}{2.5}
				\multicolumn{1}{|c}{} & \multicolumn{1}{|c|}{} & \multicolumn{3}{c|}{\textbf{Average}}  & &\textbf{213.9} &\textbf{0} &\textbf{8.3} &\textbf{774.5} &\textbf{943.4} & &\textbf{213.9} &\textbf{0} &\textbf{7.7} &\textbf{785.2} &\textbf{1211.2} \\\cline{1-2}\clineB{3-17}{2.5}

				\multicolumn{1}{|c}{\multirow{14}{*}{20}} & \multicolumn{1}{|c|}{\multirow{7}{*}{0.5}} & \multicolumn{1}{c}{\multirow{2}{*}{50}} &6 &3 &173.2 &386.1 &55 &16 &800 &T &181.2 &192.8 &6 &36 &1767 &T \\
				\multicolumn{1}{|c}{} & \multicolumn{1}{|c|}{} & \multicolumn{1}{c}{} &10 &5 &183.8 &184.7 &0 &6 &249 &334.6 &184.0 &184.7 &0 &4 &120 &303.2 \\\clineB{3-17}{2.5}
				\multicolumn{1}{|c}{} & \multicolumn{1}{|c|}{} & \multicolumn{1}{c}{\multirow{2}{*}{100}} &6 &3 &173.5 &475.8 &64 &8 &800 &T &178.2 &205.3 &13 &22 &2118 &T \\
				\multicolumn{1}{|c}{} & \multicolumn{1}{|c|}{} & \multicolumn{1}{c}{} &10 &5 &171.2 &171.3 &0  &4 &275 &407.5 &170.5 &171.3 &0 &2 &137 &233.7 \\\clineB{3-17}{2.5}
				\multicolumn{1}{|c}{} & \multicolumn{1}{|c|}{} & \multicolumn{1}{c}{\multirow{2}{*}{200}} &6 &3 &169.8 &452.5 &62 &5 &1000 &T &183.2 &318.0 &42 &8 &1471 &T \\
				\multicolumn{1}{|c}{} & \multicolumn{1}{|c|}{} & \multicolumn{1}{c}{} &10 &5 &146.4 &147.0 &0 &3 &537 &967.9 &146.5 &146.5 &0 &2 &325 &480.6 \\\clineB{3-17}{2.5}
				\multicolumn{1}{|c}{} & \multicolumn{1}{|c|}{} & \multicolumn{3}{c|}{\textbf{Average}}  & &\textbf{302.9} &\textbf{30} &\textbf{7} &\textbf{610.2} & & &\textbf{203.1} &\textbf{10} &\textbf{12.3} &\textbf{989.7} &\textbf{}  \\\cline{2-2}\clineB{3-17}{2.5}
				
				\multicolumn{1}{|c}{} & \multicolumn{1}{|c|}{\multirow{7}{*}{2}} & \multicolumn{1}{c}{\multirow{2}{*}{50}} &6 &3 &187.8 &320.8 &41 &14 &700 &T &220.6 &259.8 &15 &22 &1112 &T \\
				\multicolumn{1}{|c}{} & \multicolumn{1}{|c|}{} & \multicolumn{1}{c}{} &10 &5 &213.5 &215.9 &1 &21 &858 &T &212.8 &213.8 &0 &13 &513 &1604.5 \\\clineB{3-17}{2.5}
				\multicolumn{1}{|c}{} & \multicolumn{1}{|c|}{} & \multicolumn{1}{c}{\multirow{2}{*}{100}} &6 &3 &178.8 &521.6 &66 &8 &800 &T &189.4 &312.6 &39 &12 &1124 &T \\
				\multicolumn{1}{|c}{} & \multicolumn{1}{|c|}{} & \multicolumn{1}{c}{} &10 &5 &210.1 &210.8 &0 &9 &663 &4590.6 &210.4 &210.8 &0 &11 &624 &3920.4 \\\clineB{3-17}{2.5}
				\multicolumn{1}{|c}{} & \multicolumn{1}{|c|}{} & \multicolumn{1}{c}{\multirow{2}{*}{200}} &6 &3 &171.5 &521.6 &67 &4 &800 &T &183.2 &318.0 &42 &8 &1471 &T \\
				\multicolumn{1}{|c}{} & \multicolumn{1}{|c|}{} & \multicolumn{1}{c}{} &10 &5 &196.9 &197.4 &0 &5 &650 &2740.6 &196.9 &197.4 &0 &5 &923 &3028.5 \\\clineB{3-17}{2.5}
				\multicolumn{1}{|c}{} & \multicolumn{1}{|c|}{} & \multicolumn{3}{c|}{\textbf{Average}}  & &\textbf{331.3} &\textbf{29} &\textbf{10.2} &\textbf{745.2} & & &\textbf{252.1} &\textbf{16} &\textbf{11.8} &\textbf{961.2} & \\\hline
		\end{tabular}}
	\end{center}
\label{tbl:2-Stg-FLRDRO(W)}
\end{table}

\subsection{Computational Results of the Capacitated Case}
In this subsection, we present computational results for instances of various models with capacity considerations. Capacity parameters $F_j$ are set to the total demands from a random number of neighboring sites for $j\in J$. Parameters determining the size of instances, e.g., $|I|$, $p$ and $k$, are set to smaller values compared to those for the uncapacitated case, due to the increased computational challenge. All results are presented in Tables \ref{Capacitated_P_m} and \ref{Capacitated_P_w}, noting that column ``$|\bsm{\widehat \xi}|^o$'' reports  the number of scenarios in the optimality  set, and column ``$|\bsm{\widehat \xi}|$'' shows the number of scenarios in the complete set from solving both $\mathbf{WCEV(\it{O})}$ and $\mathbf{WCEV(\it{F})}$. Hence, the number of scenarios in the feasibility set is the difference between them. 
	
From those tables, it is clear that those instances are more difficult to solve. On average, the number of iterations, the computational time, and the optimality gap are all significantly more than those of uncapacitated instances. It is worth noting that the number of scenarios in the feasibility sets is often non-trivial. Especially for instances with a discrete sample space, it can outnumber the optimality one by an order of magnitude. This indicates that the importance of detecting and addressing the infeasibility issue in practice, although it has been overlooked in the existing literature. Compared to the uncapacitated case, we also note that it is more clear that the number of iterations is often smaller for larger $N$, which shows that empirical samples are very useful in addressing the complex structure of the capacitated case. 

\begin{table}[H]
	\caption{Computational Results of Capacitated Models on $\mathcal P^m$}
	\label{Capacitated_P_m}
	\begin{center}
		\scalebox{0.6}{
			\begin{tabular}{c|cc|ccccccc|cc|ccccccc|}
				\hline
				\multicolumn{1}{|c|}{\multirow{2}{*}{$|I|$}} & \multicolumn{1}{c}{\multirow{2}{*}{$\mathfrak{r}$}} & \multicolumn{1}{c|}{\multirow{2}{*}{$p$}} & \multicolumn{7}{c|}{$\mathcal P_{\mathbf d}^m$}  & \multicolumn{1}{c}{\multirow{2}{*}{$p$}} & \multicolumn{1}{c|}{\multirow{2}{*}{$k$}} & \multicolumn{7}{c|}{$\mathcal P_{\mathbf u}^m$}\\\cline{4-10} \cline{13-19}
				\multicolumn{1}{|c|}{}  & & & {LB} & {UB} & {Gap(\%)} & {Iter.} & {$|\bsm{\widehat \xi}^o|$} & {$|\bsm{\widehat \xi}|$} & {Time(s)} & & & {LB} & {UB} & {Gap(\%)} & {Iter.} & {$|\bsm{\widehat \xi}^o|$} & {$|\bsm{\widehat \xi}|$} & {Time(s)}\\ \hline \hline
				
				\multicolumn{1}{|c|}{\multirow{7}{*}{15}} & \multicolumn{1}{c}{\multirow{2}{*}{0.65}} & \multicolumn{1}{c|}{6} &52.8 &52.8 &0 &5 &27 &38 &23.0 & \multicolumn{1}{c}{6}& \multicolumn{1}{c|}{3} &322.3 &322.3 &0 &88 &12 &106 &754.8\\
				\multicolumn{1}{|c|}{}  & & \multicolumn{1}{c|}{10} &16.5 &16.5 &0 &4 &11 &24 &5.0 & \multicolumn{1}{c}{10}& \multicolumn{1}{c|}{5} &183.6 &184.3 &0 &53 &107 &107 &236.2 \\\clineB{2-19}{2.5}
				\multicolumn{1}{|c|}{} &\multicolumn{1}{c}{\multirow{2}{*}{0.80}}  & \multicolumn{1}{c|}{6} &51.4 &51.4 &0 &5 &24 &39 &24.8 & \multicolumn{1}{c}{6}& \multicolumn{1}{c|}{3} &331.4 &332.6 &0 &89 &9 &105 &642.1 \\
				\multicolumn{1}{|c|}{} & & \multicolumn{1}{c|}{10} &15.8 &15.8 &0 &4 &16 &30 &5.4 & \multicolumn{1}{c}{10}& \multicolumn{1}{c|}{5} &200.7 &201.7 &0 &45 &88 &88 &158.5 \\\clineB{2-19}{2.5}
				\multicolumn{1}{|c|}{} &\multicolumn{1}{c}{\multirow{2}{*}{0.95}}  & \multicolumn{1}{c|}{6} &50.8 &50.8 &0 &4 &24 &31 &27.5 & \multicolumn{1}{c}{6}& \multicolumn{1}{c|}{3} &340.3 &340.7 &0 &100 &5 &125 &676.1\\
				\multicolumn{1}{|c|}{} & & \multicolumn{1}{c|}{10} &15.1 &15.1 &0 &4 &12 &25 &6.6 & \multicolumn{1}{c}{10}& \multicolumn{1}{c|}{5} &213.3 &214.0 &0 &17 &37 &37 &29.4 \\\clineB{2-19}{2.5}
				\multicolumn{1}{|c|}{} & \multicolumn{2}{c|}{\textbf{Average}} & &\textbf{} &\textbf{0} &\textbf{4.3} &\textbf{19.0} &\textbf{31.2} &\textbf{15.4} &\multicolumn{2}{c|}{\textbf{Average}} & &\textbf{} &\textbf{0} &\textbf{65.3} &\textbf{43.0} &\textbf{94.7} &\textbf{416.2} \\\cline{1-1}\clineB{2-19}{2.5}

				\multicolumn{1}{|c|}{\multirow{7}{*}{20}} & \multicolumn{1}{c}{\multirow{2}{*}{0.65}} & \multicolumn{1}{c|}{6} &81.8 &81.9 &0 &5 &47 &64 &66.1 & \multicolumn{1}{c}{6}& \multicolumn{1}{c|}{3} &280.5 &324.9 &14 &98 &66 &125 &T\\
				\multicolumn{1}{|c|}{} & & \multicolumn{1}{c|}{10} &37.3 &37.3 &0 &4 &22 &38 &24.1 & \multicolumn{1}{c}{10}& \multicolumn{1}{c|}{5} &211.1 &228.7 &8 &70 &172 &172 &T \\\clineB{2-19}{2.5}
				\multicolumn{1}{|c|}{}  &\multicolumn{1}{c}{\multirow{2}{*}{0.80}}  & \multicolumn{1}{c|}{6} &80.7 &81.0 &0 &5 &46 &64 &76.2 & \multicolumn{1}{c}{6}& \multicolumn{1}{c|}{3} &295.7 &334.0 &11 &105 &61 &133 &T \\
				\multicolumn{1}{|c|}{} & & \multicolumn{1}{c|}{10} &35.9 &35.9 &0 &4 &25 &47 &28.5 & \multicolumn{1}{c}{10}& \multicolumn{1}{c|}{5} &230.9 &245.8 &6 &73 &178 &178 &T \\\clineB{2-19}{2.5}
				\multicolumn{1}{|c|}{}  &\multicolumn{1}{c}{\multirow{2}{*}{0.95}}  & \multicolumn{1}{c|}{6} &79.7 &79.9 &0 &5 &41 &53 &87.9 & \multicolumn{1}{c}{6}& \multicolumn{1}{c|}{3} &306.8 &339.9 &10 &115 &63 &146 &T\\
				\multicolumn{1}{|c|}{} & & \multicolumn{1}{c|}{10} &34.3 &34.4 &0 &4 &26 &40 &34.6 & \multicolumn{1}{c}{10}& \multicolumn{1}{c|}{5} &247.0 &255.8 &3 &85 &199 &199 &T \\\clineB{2-19}{2.5}
				\multicolumn{1}{|c|}{} & \multicolumn{2}{c|}{\textbf{Average}} & &\textbf{} &\textbf{0} &\textbf{4.5} &\textbf{34.5} &\textbf{51.0} &\textbf{52.9} &\multicolumn{2}{c|}{\textbf{Average}} & &\textbf{} &\textbf{9} &\textbf{91.0} &\textbf{123.2} &\textbf{158.8} &\textbf{}\\\hline
		\end{tabular}}
	\end{center}
\end{table}

\begin{table}[H]
	\caption{Computational Results of Capacitated Models on $\mathcal P^W$}
	\label{Capacitated_P_w}
	\begin{center}
		\scalebox{0.55}{
			\begin{tabular}{c|c|cc|ccccccc|c|cc|ccccccc|}
				\hline
				\multicolumn{1}{|c|}{\multirow{2}{*}{$|I|$}} & \multicolumn{1}{c|}{\multirow{2}{*}{$\rm{r}_{\mathbf d}$}} & \multicolumn{1}{c}{\multirow{2}{*}{$N$}} & \multicolumn{1}{c|}{\multirow{2}{*}{$p$}}& \multicolumn{7}{c|}{$\mathcal P^W_{\mathbf d}$} & \multicolumn{1}{c|}{\multirow{2}{*}{$\rm{r}_{\mathbf u}$}} &\multicolumn{1}{c}{\multirow{2}{*}{$p$}} & \multicolumn{1}{c|}{\multirow{2}{*}{$k$}}& \multicolumn{7}{c|}{$\mathcal P^W_{\mathbf u}$} \\ \cline{5-11}\cline{15-21}
				\multicolumn{1}{|c|}{} & \multicolumn{1}{c|}{} & &\multicolumn{1}{c|}{} & {LB} & {UB} & {Gap(\%)} & {Iter.} & {$|\widehat{\bsm \xi}^o|$} & {$|\widehat{\bsm \xi}|$} & \multicolumn{1}{c|}{Time(s)} &\multicolumn{1}{c|}{} &\multicolumn{1}{c}{}  &\multicolumn{1}{c|}{}  & {LB} & {UB} & {Gap(\%)} & {Iter.} & {$|\widehat{\bsm \xi}^o|$} & {$|\widehat{\bsm \xi}|$} & {Time(s)} \\ \hline \hline
				\multicolumn{1}{|c|}{\multirow{14}{*}{15}} & \multicolumn{1}{c|}{\multirow{7}{*}{2}} & \multicolumn{1}{c}{\multirow{2}{*}{50}} &6 &93.2 &93.2 &0 &2 &2 &4 & 114.4  &\multicolumn{1}{c|}{\multirow{7}{*}{0.5}}  &6 &3 &286.3 &286.3 &0 &18 &13 &233 &148.2  \\
				\multicolumn{1}{|c|}{} & \multicolumn{1}{c|}{} & &10 &20.1 &20.1 &0 &1 &14 &14 &194.5 & &10 &5 &184.0 &184.4 &0 &1 &13 &13 &26.8 \\\clineB{3-11}{2.5} \clineB{13-21}{2.5}
				\multicolumn{1}{|c|}{} & \multicolumn{1}{c|}{} &\multicolumn{1}{c}{\multirow{2}{*}{100}} &6 &92.8 &92.8 &0 &2 &5 &9 &272.5 & &6 &3 &255.3 &255.3 &0 &6 &33 &133 &117.6   \\
				\multicolumn{1}{|c|}{} & \multicolumn{1}{c|}{} & &10 &20.0 &20.0 &0 &1 &27 &27 &371.6 & &10 &5 &183.2 &183.4 &0 &3 &50 &50 &138.5   \\\clineB{3-11}{2.5} \clineB{13-21}{2.5}
				\multicolumn{1}{|c|}{} & \multicolumn{1}{c|}{} & \multicolumn{1}{c}{\multirow{2}{*}{200}} &6 &92.3 &92.3 &0 &2 &8 &15 &445.3 & &6 &3 &241.5 &241.5 &0 &3 &75 &125 &206.1  \\
				\multicolumn{1}{|c|}{} & \multicolumn{1}{c|}{} & \multicolumn{1}{c}{} &10 &19.9 &19.9 &0 &1 &31 &31 &585.3 & &10 &5 &171.3 &171.4 &0 &2 &72 &72 &343.5   \\\clineB{3-11}{2.5} \clineB{13-21}{2.5}
				\multicolumn{1}{|c|}{} & \multicolumn{1}{c|}{} & \multicolumn{2}{c|}{\textbf{Average}}  & & &\textbf{0} &\textbf{1.5} &\textbf{14.5} &\textbf{16.7} &\textbf{330.6} &  &\multicolumn{2}{c|}{\textbf{Average}} & & &\textbf{0} &\textbf{5.5} &\textbf{42.7} &\textbf{104.3} &\textbf{163.4} \\\cline{2-2}\clineB{3-11}{2.5} \cline{12-12} \clineB{13-21}{2.5}
				
				\multicolumn{1}{|c|}{} & \multicolumn{1}{c|}{\multirow{7}{*}{10}} & \multicolumn{1}{c}{\multirow{2}{*}{50}} &6 &98.2 &98.2 &0 &2 &10 &18 &118.7 &\multicolumn{1}{c|}{\multirow{7}{*}{2}}  &6 &3 &331.8 &331.8 &0 &45 &71 &1341 &1640.4  \\
				\multicolumn{1}{|c|}{} & \multicolumn{1}{c|}{} & &10 &21.2 &21.2 &0 &1 &51 &51 &191.1 & &10 &5 &204.9 &205.7 &0 &6 &208 &208 &209.0 \\\clineB{3-11}{2.5} \clineB{13-21}{2.5}
				\multicolumn{1}{|c|}{} & \multicolumn{1}{c|}{} &\multicolumn{1}{c}{\multirow{2}{*}{100}} &6 &97.9 &97.9 &0 &2 &20 &35 &262.6 & &6 &3 &316.5 &316.7 &0 &22 &125 &1334 &1083.4   \\
				\multicolumn{1}{|c|}{} & \multicolumn{1}{c|}{} & &10 &21.1 &21.1 &0 &1 &101 &101 &362.2 & &10 &5 &202.1 &203.0 &0 &5 &308 &308 &358.7   \\\clineB{3-11}{2.5} \clineB{13-21}{2.5}
				\multicolumn{1}{|c|}{} & \multicolumn{1}{c|}{} & \multicolumn{1}{c}{\multirow{2}{*}{200}} &6 &97.4 &97.4 &0 &2 &38 &68 &447.7 & &6 &3 &309.6 &310.2 &0 &13 &233 &1404 &1346.4  \\
				\multicolumn{1}{|c|}{} & \multicolumn{1}{c|}{} & \multicolumn{1}{c}{} &10 &21.0 &21.0 &0 &1 &148 &148 &589.3 & &10 &5 &196.1 &197.0 &0 &6 &595 &595 &1001.4   \\\clineB{3-11}{2.5} \clineB{13-21}{2.5}
				\multicolumn{1}{|c|}{} & \multicolumn{1}{c|}{} & \multicolumn{2}{c|}{\textbf{Average}}  & & &\textbf{0}  &\textbf{1.5} &\textbf{61.3} &\textbf{70.2} &\textbf{328.6} &  &\multicolumn{2}{c|}{\textbf{Average}} & & &\textbf{0} &\textbf{16.2} &\textbf{256.7} &\textbf{865.0} &\textbf{939.9} \\\cline{1-2}\clineB{3-11}{2.5} \cline{12-12} \clineB{13-21}{2.5}

				\multicolumn{1}{|c|}{\multirow{14}{*}{20}} & \multicolumn{1}{c|}{\multirow{7}{*}{2}} & \multicolumn{1}{c}{\multirow{2}{*}{50}} &6 &114.3 &114.3 &0 &1 &51 &51 &331.5 &\multicolumn{1}{c|}{\multirow{7}{*}{0.5}}  &6 &3 &197.0 &197.5 &0 &30 &173 &213 &5465.2  \\
				\multicolumn{1}{|c|}{} & \multicolumn{1}{c|}{} & &10 &44.1 &44.1 &0 &1 &43 &43 &594.6 & &10 &5 &184.5 &185.2 &0 &3 &25 &25 &99.5 \\\clineB{3-11}{2.5} \clineB{13-21}{2.5}
				\multicolumn{1}{|c|}{} & \multicolumn{1}{c|}{} &\multicolumn{1}{c}{\multirow{2}{*}{100}} &6 &114.0 &114.0 &0 &1 &95 &95 &717.2 & &6 &3 &190.0 &199.5 &5 &18 &243 &269 &T   \\
				\multicolumn{1}{|c|}{} & \multicolumn{1}{c|}{} & &10 &44.0 &44.0 &0 &1 &60 &60 &1037.5 & &10 &5 &170.9 &171.7 &0 &2 &29 &29 &158.4   \\\clineB{3-11}{2.5} \clineB{13-21}{2.5}
				\multicolumn{1}{|c|}{} & \multicolumn{1}{c|}{} & \multicolumn{1}{c}{\multirow{2}{*}{200}} &6 &113.7 &113.7 &0 &1 &98 &98 &1485.1 & &6 &3 &180.8 &210.2 &14 &9 &464 &464 &T  \\
				\multicolumn{1}{|c|}{} & \multicolumn{1}{c|}{} & \multicolumn{1}{c}{} &10 &43.8 &43.8 &0 &1 &41 &41 &1556.9 & &10 &5 &147.3 &147.6 &0 &1 &50 &50 &198.8   \\\clineB{3-11}{2.5} \clineB{13-21}{2.5}
				\multicolumn{1}{|c|}{} & \multicolumn{1}{c|}{} & \multicolumn{2}{c|}{\textbf{Average}}  & & &\textbf{0} &\textbf{1.0} &\textbf{64.7} &\textbf{64.7} &\textbf{953.8} & &\multicolumn{2}{c|}{\textbf{Average}} & & &\textbf{3} &\textbf{10.5} &\textbf{164.0} &\textbf{175.0} &\textbf{}  \\\cline{2-2}\clineB{3-11}{2.5} \cline{12-12} \clineB{13-21}{2.5}
				
				\multicolumn{1}{|c|}{} & \multicolumn{1}{c|}{\multirow{7}{*}{10}} & \multicolumn{1}{c}{\multirow{2}{*}{50}} &6 &118.4 &118.4 &0 &1 &51 &51 &316.0 &\multicolumn{1}{c|}{\multirow{7}{*}{2}}  &6 &3 &207.3 &277.2 &25 &17 &397 &397 &T  \\
				\multicolumn{1}{|c|}{} & \multicolumn{1}{c|}{} & &10 &45.5 &45.5 &0 &1 &51 &51 &443.3 & &10 &5 &213.4 &214.2 &0 &10 &188 &188 &723.7 \\\clineB{3-11}{2.5} \clineB{13-21}{2.5}
				\multicolumn{1}{|c|}{} & \multicolumn{1}{c|}{} &\multicolumn{1}{c}{\multirow{2}{*}{100}} &6 &118.1 &118.1 &0 &1 &101 &101 &695.8 & &6 &3 &191.2 &290.1 &34 &9 &432 &432 &T   \\
				\multicolumn{1}{|c|}{} & \multicolumn{1}{c|}{} & &10 &45.4 &45.4 &0 &1 &101 &101 &781.1 & &10 &5 &212.0 &212.8 &0 &7 &321 &321 &1925.9   \\\clineB{3-11}{2.5} \clineB{13-21}{2.5}
				\multicolumn{1}{|c|}{} & \multicolumn{1}{c|}{} & \multicolumn{1}{c}{\multirow{2}{*}{200}} &6 &117.8 &117.8 &0 &1 &150 &150 &2047.9& &6 &3 &186.9 &357.5 &48 &6 &697 &697 &T  \\
				\multicolumn{1}{|c|}{} & \multicolumn{1}{c|}{} & \multicolumn{1}{c}{} &10 &45.4 &45.4 &0 &1 &191 &191 &1901.9 & &10 &5 &198.7 &199.4 &0 &5 &568 &568 &3802.5   \\\clineB{3-11}{2.5} \clineB{13-21}{2.5}
				\multicolumn{1}{|c|}{} & \multicolumn{1}{c|}{} & \multicolumn{2}{c|}{\textbf{Average}}  & & &\textbf{0} &\textbf{1.0} &\textbf{107.5} &\textbf{107.5} &\textbf{1031.0} & &\multicolumn{2}{c|}{\textbf{Average}} & & &\textbf{18} &\textbf{9.0} &\textbf{433.8} &\textbf{433.8} &\textbf{} \\\hline
		\end{tabular}}
	\end{center}
\end{table}

\subsection{Computational Results of Complex Ambiguity Sets}
\label{subsect_compu_complexAMB}
In this subsection, we consider  uncapacitated $\mathbf{FL-DRO}$ defined on more sophisticated ambiguity sets. Note that our purpose here is to demonstrate C\&CG-DRO's (with \textit{Oracle-2}) capacity to handle complex problems, rather than to provide a comprehensive evaluation.
One type of ambiguity sets is Wasserstein metric-based set defined by $L_2$ norm, denoted by $\cal P^{W_2}_{\bf d}$. Note that it has a second-order conic (SOC) structure and requires SOC programming  solver within \textit{Oracle-2} to compute $\mathbf{PSP}$ problems. Another type is mixed integer ambiguity sets extending $\cal P^m_{\bf u}$ as in the following, which is denoted by $\mathcal P_{\mathbf u}^{mI}$.  
\begin{align*}
	\begin{split}
		\mathcal P_{\mathbf u}^{mI} = \Big\{P \in \mathcal{M}(\mathcal{U},\mathcal{F}) : \underline{\bm k}-\theta\bf z\leq E_{P} [\mathbf u] \leq \tilde{\bm k}-\theta\bf z, \ \sum_i z_i\geq z^0, \ z_i\in \{0,1\} \ \ \forall i\in I
		\Big\}.
	\end{split}
\end{align*}
In our numerical study, $\underline{\bm k}=0.4*\bf 1$, $\tilde{\bm k}=0.8*\bf 1$, $\theta$ is set to 0.4 and $z^0$ to 2. All results are presented in Table \ref{tbl:SOCP_MIP}. Compared to the results in Tables \ref{tbl:2-Stg-FLDDRO(W)} and \ref{tbl:2-Stg-FLRDRO(m)}, there is no significant difference in the computational performance, including time and the number of iterations. The C\&CG-DRO method remains effective for handling $L_2$ norm-defined Wasserstein metric-based ambiguity sets, while it is still sensitive to the instance scale when the sample space is discrete. However, we note that this observation is based on this small-scale study and is not conclusive. 

	\begin{table}[H]
	\caption{Computational Results of Uncapacitated Models with $\mathcal P_{\mathbf d}^{W_2}$ and $\mathcal P_{\mathbf u}^{mI}$}
	\begin{center}
		\scalebox{0.6}{
			\begin{tabular}{c|ccc|cccccc|ccc|cccccc|}
				\hline
				\multicolumn{1}{|c|}{\multirow{2}{*}{$|I|$}} &\multicolumn{1}{c}{\multirow{2}{*}{$\rm{r}^2_{\bf d}$}} & \multicolumn{1}{c}{\multirow{2}{*}{$N$}} & \multicolumn{1}{c|}{\multirow{2}{*}{$p$}} &  \multicolumn{6}{c|}{SOC Ambiguity Set ($\mathcal P_{\mathbf d}^{W_2}$)} & \multicolumn{1}{c}{\multirow{2}{*}{$\mathfrak{r}$}} & \multicolumn{1}{c}{\multirow{2}{*}{$p$}} & \multicolumn{1}{c|}{\multirow{2}{*}{$k$}} & \multicolumn{6}{c|}{MIP Ambiguity Set ($\mathcal P_{\mathbf u}^{mI}$)}\\ \cline{5-10} \cline{14-19}
				\multicolumn{1}{|c|}{} & && & {LB} & {UB} & {Gap(\%)} & {Iter.} & {$|\widehat{\bsm \xi}|$} & {Time(s)}& & & & {LB} & {UB} & {Gap(\%)} & {Iter.} & {$|\widehat{\bsm \xi}|$} & {Time(s)} \\ \hline \hline
				\multicolumn{1}{|c}{\multirow{7}{*}{15}} &\multicolumn{1}{|c}{\multirow{3}{*}{2}} &50 & 6 &62.4 &62.7 &0 &1 &51 &31.6 & \multicolumn{1}{c}{\multirow{2}{*}{0.65}}  &\multicolumn{1}{c}{6} & \multicolumn{1}{c|}{3} &235.6 &236.4 &0 &34 &69 &216.4 \\
				\multicolumn{1}{|c}{} & \multicolumn{1}{|c}{} &100 & 6 &62.3 &62.5 &0 &1 &101 &58.5  & &10 & \multicolumn{1}{c|}{5} &193.3 &194.2 &0 &30 &76 &103.0 \\\clineB{11-19}{2.5}
				\multicolumn{1}{|c}{} &\multicolumn{1}{|c}{}  &200 &6 &62.1 &62.4 &0 &1 &201 &135.2 &\multicolumn{1}{c}{\multirow{2}{*}{0.80}}&\multicolumn{1}{c}{6} & \multicolumn{1}{c|}{3}  &235.4 &236.4 &0 &34 &64 &195.3 \\ \clineB{2-10}{2.5}
				\multicolumn{1}{|c}{} &\multicolumn{1}{|c}{\multirow{3}{*}{30}}  &50 &6 &64.4 &64.7 &0 &2 &98 &72.3 & &\multicolumn{1}{c}{10} & \multicolumn{1}{c|}{5} &194.7 &195.1 &0 &27 &58 &83.4 \\\clineB{11-19}{2.5}
				\multicolumn{1}{|c}{} &\multicolumn{1}{|c}{} &100 & \multicolumn{1}{c|}{6} &64.3 &64.6 &0 &2 &198 &150.9 &\multicolumn{1}{c}{\multirow{2}{*}{0.95}} &\multicolumn{1}{c}{6} & \multicolumn{1}{c|}{3} &235.4 &236.4 &0 &33 &63 &191.1 \\
				\multicolumn{1}{|c}{} & \multicolumn{1}{|c}{} &200 & \multicolumn{1}{c|}{6} &64.2 &64.4 &0 &2 &397 &333.3 &  &\multicolumn{1}{c}{10} & \multicolumn{1}{c|}{5} &195.0 &195.2 &0 &27 &60 &68.1  \\\clineB{2-19}{2.5}
				\multicolumn{1}{|c}{ } & \multicolumn{3}{|c|}{\textbf{Average}} & &\textbf{} &\textbf{0} &\textbf{1.5} &\textbf{174.3} &\textbf{130.3}  & \multicolumn{3}{c|}{\textbf{Average}} & &\textbf{} &\textbf{0} &\textbf{30.8} &\textbf{65.0} &\textbf{68.1} \\\cline{1-1}\clineB{2-19}{2.5}

				\multicolumn{1}{|c}{\multirow{7}{*}{20}} &\multicolumn{1}{|c}{\multirow{3}{*}{2}} &50 & 6 &93.6 &93.6 &0 &1 &51 &72.0 & \multicolumn{1}{c}{\multirow{2}{*}{0.65}}  &\multicolumn{1}{c}{6} & \multicolumn{1}{c|}{3} &276.6 &289.3 &4 &70 &130 &T \\
				\multicolumn{1}{|c}{} & \multicolumn{1}{|c}{} &100 & 6 &93.4 &93.4 &0 &1 &101 &217.3  & &10 & \multicolumn{1}{c|}{5} &212.5 &225.2 &6 &69 &165 &T  \\\clineB{11-19}{2.5}
				\multicolumn{1}{|c}{} &\multicolumn{1}{|c}{}  &200 &6 &93.2 &93.2 &0 &1 &201 &776.9 &\multicolumn{1}{c}{\multirow{2}{*}{0.80}}&\multicolumn{1}{c}{6} & \multicolumn{1}{c|}{3}  &279.7 &290.5 &4 &74 &122 &T \\ \clineB{2-10}{2.5}
				\multicolumn{1}{|c}{} &\multicolumn{1}{|c}{\multirow{3}{*}{30}}  &50 &6 &96.2 &96.6 &0 &1 &51 &69.2 & &\multicolumn{1}{c}{10} & \multicolumn{1}{c|}{5} &225.4 &233.2 &3 &65 &155 &T \\\clineB{11-19}{2.5}
				\multicolumn{1}{|c}{} &\multicolumn{1}{|c}{} &100 & \multicolumn{1}{c|}{6} &96.1 &96.4 &0 &1 &101 &222.4 &\multicolumn{1}{c}{\multirow{2}{*}{0.95}} &\multicolumn{1}{c}{6} & \multicolumn{1}{c|}{3} &280.5 &293.2 &4 &65 &133 &T \\
				\multicolumn{1}{|c}{} & \multicolumn{1}{|c}{} &200 & \multicolumn{1}{c|}{6} &96.0 &96.3 &0 &1 &201 &687.3 &  &\multicolumn{1}{c}{10} & \multicolumn{1}{c|}{5} &224.9 &233.2 &4 &66 &170 &T \\\clineB{2-19}{2.5}
				\multicolumn{1}{|c}{ } & \multicolumn{3}{|c|}{\textbf{Average}} & &\textbf{} &\textbf{0} &\textbf{1.0} &\textbf{117.7} &\textbf{340.9}  & \multicolumn{3}{c|}{\textbf{Average}} & &\textbf{} &\textbf{4} &\textbf{68.2} &\textbf{145.8} & \\\hline
		\end{tabular}}
	\end{center}\label{tbl:SOCP_MIP}
\end{table}

\section{Conclusions}
\label{sect_conclusion}

In this paper, rather than following the dual perspective that is popular in the literature, we present a new study on two-stage DRO by taking the primal perspective. This perspective allows us to gain a  deeper and more  intuitive  understanding on DRO, to develop a general and fast decomposition algorithm (and its variants) by leveraging existing powerful solution methods, and to address a couple of unsolved issues underlying two-stage DRO. Theoretical analyses regarding the strength, convergence, and iteration complexity of the developed algorithm are also presented. A systematic numerical study on the distributionally robust facility location problem has been conducted, taking into account multiple critical factors.. Results clearly demonstrate that our new solution algorithm (and its variants) is generally applicable and achieves remarkable superiority over existing methods, often solving problems up to several orders of magnitude faster. 

Regarding future research directions, it would be interesting to  extend our primal perspective to consider other types of risk measures beyond the expected value under the DRO framework. Also, enhancing the developed algorithm (and its variants) and promoting it in solving large-scale and complex data-driven problems will be carried out to support various real-world systems.

\medskip

\bibliographystyle{unsrt}
\bibliography{paper}

\begin{thebibliography}{10}

\bibitem{scarf1958min}
Herbert Scarf.
\newblock A min max solution of an inventory problem.
\newblock {\em Studies in the Mathematical Theory of Inventory and Production},
  1958.

\bibitem{rahimian2022frameworks}
Hamed Rahimian and Sanjay Mehrotra.
\newblock Frameworks and results in distributionally robust optimization.
\newblock {\em Open Journal of Mathematical Optimization}, 3:1--85, 2022.

\bibitem{delage2010distributionally}
Erick Delage and Yinyu Ye.
\newblock Distributionally robust optimization under moment uncertainty with
  application to data-driven problems.
\newblock {\em Operations Research}, 58(3):595--612, 2010.

\bibitem{wiesemann2014distributionally}
Wolfram Wiesemann, Daniel Kuhn, and Melvyn Sim.
\newblock Distributionally robust convex optimization.
\newblock {\em Operations Research}, 62(6):1358--1376, 2014.

\bibitem{ben2013robust}
Aharon Ben-Tal, Dick Den~Hertog, Anja De~Waegenaere, Bertrand Melenberg, and
  Gijs Rennen.
\newblock Robust solutions of optimization problems affected by uncertain
  probabilities.
\newblock {\em Management Science}, 59(2):341--357, 2013.

\bibitem{lam2019recovering}
Henry Lam.
\newblock Recovering best statistical guarantees via the empirical
  divergence-based distributionally robust optimization.
\newblock {\em Operations Research}, 67(4):1090--1105, 2019.

\bibitem{gao2023distributionally}
Rui Gao and Anton Kleywegt.
\newblock Distributionally robust stochastic optimization with wasserstein
  distance.
\newblock {\em Mathematics of Operations Research}, 48(2):603--655, 2023.

\bibitem{mohajerin2018data}
Peyman Mohajerin~Esfahani and Daniel Kuhn.
\newblock Data-driven distributionally robust optimization using the
  wasserstein metric: Performance guarantees and tractable reformulations.
\newblock {\em Mathematical Programming}, 171(1):115--166, 2018.

\bibitem{zhao2018data}
Chaoyue Zhao and Yongpei Guan.
\newblock Data-driven risk-averse stochastic optimization with wasserstein
  metric.
\newblock {\em Operations Research Letters}, 46(2):262--267, 2018.

\bibitem{ben2004adjustable}
Aharon Ben-Tal, Alexander Goryashko, Elana Guslitzer, and Arkadi Nemirovski.
\newblock Adjustable robust solutions of uncertain linear programs.
\newblock {\em Mathematical programming}, 99(2):351--376, 2004.

\bibitem{bertsimas2010models}
Dimitris Bertsimas, Xuan~Vinh Doan, Karthik Natarajan, and Chung-Piaw Teo.
\newblock Models for minimax stochastic linear optimization problems with risk
  aversion.
\newblock {\em Mathematics of Operations Research}, 35(3):580--602, 2010.

\bibitem{hanasusanto2018conic}
Grani~A Hanasusanto and Daniel Kuhn.
\newblock Conic programming reformulations of two-stage distributionally robust
  linear programs over wasserstein balls.
\newblock {\em Operations Research}, 66(3):849--869, 2018.

\bibitem{xie2020tractable}
Weijun Xie.
\newblock Tractable reformulations of two-stage distributionally robust linear
  programs over the type-$\infty$ wasserstein ball.
\newblock {\em Operations Research Letters}, 48(4):513--523, 2020.

\bibitem{bagheri2018resilient}
Ali Bagheri, Chaoyue Zhao, Feng Qiu, and Jianhui Wang.
\newblock Resilient transmission hardening planning in a high renewable
  penetration era.
\newblock {\em IEEE Transactions on Power Systems}, 34(2):873--882, 2018.

\bibitem{bansal2018decomposition}
Manish Bansal, Kuo-Ling Huang, and Sanjay Mehrotra.
\newblock Decomposition algorithms for two-stage distributionally robust mixed
  binary programs.
\newblock {\em SIAM Journal on Optimization}, 28(3):2360--2383, 2018.

\bibitem{saif2021data}
Ahmed Saif and Erick Delage.
\newblock Data-driven distributionally robust capacitated facility location
  problem.
\newblock {\em European Journal of Operational Research}, 291(3):995--1007,
  2021.

\bibitem{gamboa2021decomposition}
Carlos~Andr{\'e}s Gamboa, Davi~Michel Vallad{\~a}o, Alexandre Street, and Tito
  Homem-de Mello.
\newblock Decomposition methods for wasserstein-based data-driven
  distributionally robust problems.
\newblock {\em Operations Research Letters}, 49(5):696--702, 2021.

\bibitem{gangammanavar2022stochastic}
Harsha Gangammanavar and Manish Bansal.
\newblock Stochastic decomposition method for two-stage distributionally robust
  linear optimization.
\newblock {\em SIAM Journal on Optimization}, 32(3):1901--1930, 2022.

\bibitem{duque2022distributionally}
Daniel Duque, Sanjay Mehrotra, and David~P Morton.
\newblock Distributionally robust two-stage stochastic programming.
\newblock {\em SIAM Journal on Optimization}, 32(3):1499--1522, 2022.

\bibitem{shapiro2001duality}
Alexander Shapiro.
\newblock {\em On Duality Theory of Conic Linear Problems}, pages 135--165.
\newblock Springer US, Boston, MA, 2001.

\bibitem{blanchet2019quantifying}
Jose Blanchet and Karthyek Murthy.
\newblock Quantifying distributional model risk via optimal transport.
\newblock {\em Mathematics of Operations Research}, 44(2):565--600, 2019.

\bibitem{ZENG2013457}
Bo~Zeng and Long Zhao.
\newblock Solving two-stage robust optimization problems using a
  column-and-constraint generation method.
\newblock {\em Operations Research Letters}, 41(5):457--461, 2013.

\bibitem{luo2022decomposition}
Fengqiao Luo and Sanjay Mehrotra.
\newblock A decomposition method for distributionally-robust two-stage
  stochastic mixed-integer conic programs.
\newblock {\em Mathematical Programming}, 196(1):673--717, 2022.

\bibitem{el2024distributionally}
Mohamed El~Tonbari, George Nemhauser, and Alejandro Toriello.
\newblock Distributionally robust disaster relief planning under the
  wasserstein set.
\newblock {\em Computers \& Operations Research}, page 106689, 2024.

\bibitem{long2024supermodularity}
Daniel~Zhuoyu Long, Jin Qi, and Aiqi Zhang.
\newblock Supermodularity in two-stage distributionally robust optimization.
\newblock {\em Management Science}, 70(3):1394--1409, 2024.

\bibitem{jiang2018risk}
Ruiwei Jiang and Yongpei Guan.
\newblock Risk-averse two-stage stochastic program with distributional
  ambiguity.
\newblock {\em Operations Research}, 66(5):1390--1405, 2018.

\bibitem{bertsimas2019adaptive}
Dimitris Bertsimas, Melvyn Sim, and Meilin Zhang.
\newblock Adaptive distributionally robust optimization.
\newblock {\em Management Science}, 65(2):604--618, 2019.

\bibitem{georghiou2021optimality}
Angelos Georghiou, Angelos Tsoukalas, and Wolfram Wiesemann.
\newblock On the optimality of affine decision rules in robust and
  distributionally robust optimization.
\newblock {\em Available at Optimization Online}, 2021.

\bibitem{zhang2022simple}
Luhao Zhang, Jincheng Yang, and Rui Gao.
\newblock A simple and general duality proof for wasserstein distributionally
  robust optimization.
\newblock {\em arXiv preprint arXiv:2205.00362}, 2022.

\bibitem{KanZal.1951}
L.V. Kantorovic and V.A. Zalgaller.
\newblock {\em Rational cutting of industrial materials}.

\bibitem{ford1958suggested}
Lester~Randolph Ford~Jr and Delbert~R Fulkerson.
\newblock A suggested computation for maximal multi-commodity network flows.
\newblock {\em Management Science}, 5(1):97--101, 1958.

\bibitem{dantzig1960decomposition}
George~B Dantzig and Philip Wolfe.
\newblock Decomposition principle for linear programs.
\newblock {\em Operations research}, 8(1):101--111, 1960.

\bibitem{gilmore1961linear}
Paul~C Gilmore and Ralph~E Gomory.
\newblock A linear programming approach to the cutting stock problem.
\newblock {\em Operations research}, 9(6):849--859, 1961.

\bibitem{zeng2022two}
Bo~Zeng and Wei Wang.
\newblock Two-stage robust optimization with decision dependent uncertainty.
\newblock {\em arXiv preprint arXiv:2203.16484}, 2022.

\bibitem{barnhart1998branch}
Cynthia Barnhart, Ellis~L Johnson, George~L Nemhauser, Martin~WP Savelsbergh,
  and Pamela~H Vance.
\newblock Branch-and-price: Column generation for solving huge integer
  programs.
\newblock {\em Operations research}, 46(3):316--329, 1998.

\bibitem{zeng2014solving}
Bo~Zeng and Yu~An.
\newblock Solving bilevel mixed integer program by reformulations and
  decomposition.
\newblock {\em Optimization online}, pages 1--34, 2014.

\bibitem{snyder2005reliability}
Lawrence~V Snyder and Mark~S Daskin.
\newblock Reliability models for facility location: the expected failure cost
  case.
\newblock {\em Transportation science}, 39(3):400--416, 2005.

\end{thebibliography}

\end{document}